\documentclass[11pt,leqno,british]{amsart}
\usepackage{ifpdf}

\usepackage[T1]{fontenc}

\usepackage{eucal}
\usepackage{url}
\usepackage[pagebackref]{hyperref}
\usepackage{amsfonts}
\usepackage{amsthm}
\usepackage{amsmath}
\usepackage{amscd}
\usepackage{amssymb}
\usepackage{graphicx}
\usepackage{rotating}
\usepackage{wasysym}
\usepackage[small,nohug,heads=LaTeX]{diagrams}
\diagramstyle[labelstyle=\scriptstyle]
\usepackage{ulem}
\usepackage{textcomp}

\usepackage{pgf}
\usepackage{tikz}
\usetikzlibrary{matrix,arrows,decorations.pathmorphing}

\newcommand*{\DashedArrow}[1][]{\mathbin{\tikz [baseline=-0.25ex,-latex, dashed,#1] \draw [#1] (0pt,0.5ex) -- (1.3em,0.5ex);}}%

\def\CC {{\mathbb C}}     
\def\HH {{\mathbb H}}     
\def\NN {{\mathbb N}}     
\def\RR {{\mathbb R}}     
\def\SS {{\mathbb S}}     
\def\ZZ {{\mathbb Z}}     

\def\ring#1{\ifmmode \mathaccent'027 #1\else \rm\accent'027 #1\fi}

\newcommand{\ri}{{\mathrm i}}

\def\ul  {\underline}

\def\wt  {\widetilde}

\def\im {\mathfrak{im}}

\def\mc {\mathcal}
\def\mk {\mathfrak}

\def\Hom {\mathrm{Hom}}

\def\st {\mathrm{Stab}}

\def \bd {\begin{diagram}}
\def \ed {\end{diagram}}
\def\be  {\begin{eqnarray}}
\def\ee  {\end{eqnarray}}
\def\ben {\begin{eqnarray*}}
\def\een {\end{eqnarray*}}

\def\bpr {\begin{proof}[Proof]}
\def\epr {\end{proof}}
\def\bsp {\begin{split}}
\def\esp {\end{split}}
\def\bprr {\begin{proof}[solution]}
\def\bpru {\begin{proof}[hint]}
\def\bpro {\begin{proof}[answer]}
\def\bcd {\begin{CD}}
\def\ecd {\end{CD}}

\newcommand{\abs}[1]{\left\vert#1\right\vert}
\newcommand{\scal}[1]{\left\langle#1\right\rangle}

\newtheorem{theorem}{Theorem}[section]
\newtheorem{lemma}[theorem]{Lemma}
\newtheorem{prop}[theorem]{Proposition}
\newtheorem{coro}[theorem]{Corollary}
\newtheorem{remark}[theorem]{Remark}
\newtheorem{df}[theorem]{Definition}

\begin{document}

\title[Non-semistable exceptional objects in hereditary categories]%
{Non-semistable exceptional objects in hereditary categories}
\author{George Dimitrov}
\address[Dimitrov]{Universität Wien\\
Oskar-Morgenstern-Platz 1, 1090 Wien\\
Österreich}
\email{gkid@abv.bg}

\author{Ludmil Katzarkov}
\address[Katzarkov]{Universität Wien\\
Oskar-Morgenstern-Platz 1, 1090 Wien\\
Österreich}
\email{lkatzark@math.uci.edu}

\begin{abstract}
For a given stability condition $\sigma$ on a triangulated category
   we  define  a $\sigma$-exceptional collection  as an Ext-exceptional collection, whose elements
 are $\sigma$-semistable with phases contained in an open interval of length one. If there  exists a full $\sigma$-exceptional collection,
 then    $\sigma$ is  generated by this collection in a procedure described by  E. Macr\`i.

Constructing  $\sigma$-exceptional collections of length at least three  in  $D^b(\mc A)$ from a non-semistable exceptional object,
where  $\mc A$ is  a hereditary $\hom$-finite abelian category, we introduce certain conditions on   the   Ext-nontrivial couples (couples of
exceptional objects $X,Y\in \mc A$  with ${\rm Ext}^1(X,Y)\neq 0$, ${\rm Ext}^1(Y,X)\neq 0$).

     After a   detailed study of  the exceptional objects  of the quivers  $Q_1=\begin{diagram}[size=0.8em]
        &       &  \circ  &       &    \\
        & \ruTo &         & \luTo &       \\
 \circ  & \rTo  &         &       &  \circ
\end{diagram}$,  $Q_2= \begin{diagram}[size=0.8em]
 \circ &  \rTo  &  \circ    \\
  \uTo &        & \uTo     \\
 \circ   & \rTo  &    \circ
 \end{diagram}$
we observe that the needed conditions   do hold  in   $Rep_k(Q_1)$,  $Rep_k(Q_2)$.

   Combining these findings, we   prove that for each $\sigma \in \st(D^b(Q_1))$ there exists a full $\sigma$-exceptional collection. It follows that
   $\st(D^b(Q_1))$ is connected. \end{abstract}

\maketitle
\setcounter{tocdepth}{2}
\tableofcontents

\section{Introduction}

   T. Bridgeland introduced in his seminal work \cite{Bridg1} the definition of a locally finite stability condition  on a triangulated category $\mc T$,
   motivated by the work of Douglas on $\Pi$-stability for Dirichlet branes. He
       proved that the set of these stability conditions is a complex manifold, denoted by $\st(\mc T)$.

  Bridgeland's axioms imply\footnote{For a subset $S \subset Ob(\mc T)$ we denote by $\langle  S \rangle  \subset \mc T$  the triangulated subcategory of $\mc T$
   generated by $S$.} that  $\st(\langle E \rangle)=\CC$ for an exceptional object $E$ in $\mc T$.
   The guiding motivation of this paper is the study of $\st(\langle E_1,E_2,\dots, E_n \rangle )$, where  $(E_1,\dots,E_n)$ is an exceptional collection in $\mc T$
    and $n\geq 2$.  This study was initiated by E. Macr\`i in \cite{Macri}. Here, we  proceed further.

    Collins and Polishchuk     defined and studied  in \cite{CP} a gluing procedure for Bridgeland stability conditions in the situation when $\mc T$ has a
     semiorthogonal decomposition $\mc T=\langle \mc A_1, \mc A_2 \rangle $.

   \vspace{3mm}

\noindent 1.1. T. Bridgeland constructed  a stability condition $\sigma \in \st(\mc T)$  from a bounded t-structure  $\mc A \subset \mc T$  and a stability
function\footnote{I.e. $Z$ is homomorphism  $\bd K(\mc A)&\rTo^{Z}&\CC \ed$, s. t. $Z(X)\in \HH=\{ r \exp(\ri \pi t): r>0 \ \ \mbox{and} \ \ 0< t \leq 1 \}$ for
$X \in \mc A\setminus \{0\}$.}  $Z:K(\mc A) \rightarrow \CC$ satisfying certain restrictions. Keeping $\mc A$  fixed and varying $Z$ produces  a family      of
 stability conditions,  which we denote by $\HH^{\mc A}\subset \st(\mc T)$.
E. Macr\`i  proved in \cite[Lemma 3.14]{Macri},   using results of \cite{BBD},   that the extension closure  ${\mc A}_{\mc E}$  of  a full Ext-exceptional
collection\footnote{An exceptional  collection ${\mc E}=(E_0,E_1,\dots,E_n)$ is said to be Ext-exceptional if $\forall i\neq j$ $\Hom^{\leq 0}(E_i, E_j)=0$.}
    ${\mc E}=(E_0,E_1,\dots,E_n)$  in $\mc T$  is a heart of a bounded t-structure, and for  each $\sigma \in \HH^{{\mc A}_{\mc E}}$ the objects
    $E_0,E_1,\dots,E_n$ are $\sigma$-stable with phases in $ (0,1]$. Motivated by this result, for a given $\sigma \in \st(\mc T)$ we define a
    \textit{$\sigma$-exceptional collection}(Definition \ref{sigma exceptional collection})  as an  Ext-exceptional collection ${\mc E } = (E_0,E_1,\dots,E_n)$,
     s. t. the
     objects $\{E_i\}_{i=0}^n$ are $\sigma$-semistable, and  $\{\phi(E_i)\}_{i=0}^n \subset (t,t+1)$
 for some $t\in \RR$ .  It follows easily from \cite[Lemmas 3.14, 3.16]{Macri} that for any  full Ext-exceptional collection $ \mc E $  the set
 $\{ \sigma \in \st(\mc T)$:  $ \mc E $ is $\sigma$-exceptional$\}$ coincides with\footnote{Recall that $\st(\mc T)$ carries a  right action
   by $\wt{GL}^+(2,\RR)$.} $\HH^{{\mc A}_{\mc E}}\cdot\widetilde{GL}^+(2,\RR)$(Corollary \ref{coro of Macri}).

E. Macr\`i, studying   $ \st(D^b(K(l))$ in  \cite{Macri},     gave  an idea for producing a $\sigma$-exceptional pair in  $D^b(K(l))$ from
 a non-semistable exceptional object, where $K(l)$ is the $l$-Kronecker quiver.

 Throughout sections  \ref{non-stable exc obj in...}, \ref{some terminilogy}, \ref{no bad after good}, \ref{sequence}, \ref{final}, \ref{constructing}
   we develop tools   for constructing $\sigma$-exceptional collections of length at least three in  $D^b(\mc A)$, where  $\mc A$ is  a hereditary
   $\hom$-finite abelian category.      Combining them with  the findings of Section \ref{two examples} about $Rep_k(Q_1)$ we    prove  in
   Section \ref{main theo for Q_1}  the following theorem:

 \begin{theorem} \label{main theorem for Q_1 in intro} Let $Q_1$ be the quiver  $\begin{diagram}[size=1em]
   &       &  \circ  &       &    \\
   & \ruTo &    & \luTo &       \\
\circ  & \rTo  &    &       &  \circ
\end{diagram}$. Let $k$ be an algebraically closed field.
 For each $\sigma \in \st(D^b(Rep_k(Q_1)))$ there exists a full
$\sigma$-exceptional collection.
\end{theorem}

Theorem \ref{main theorem for Q_1 in intro} is one  novelty of this paper.  In particular, it implies
 that $\st(D^b(Rep_k(Q_1)))$ is connected (Corollary \ref{connectedness}).

   The $K(l)$-analogue of Theorem \ref{main theorem for Q_1 in intro}(Lemma \ref{main theorem for K(l)}) is already  treated by E. Macr\`i
   in \cite[Lemma 4.2 on p.10]{Macri}. For the sake of completeness, we  add a proof of this analogue in Appendix \ref{exceptional pairs}.

   The proof of  Theorem \ref{main theorem for Q_1 in intro}  is   more complicated than of its $K(l)$-analogue  not only because the full collections are triples
   instead of pairs, but also due to the presence of
     Ext-nontrivial   couples\footnote{These are couples of exceptional objects $X,Y$ with ${\rm Ext}^1(X,Y)\neq 0$, ${\rm Ext}^1(Y,X)\neq 0$(Definition \ref{Ext-nontrivial couple}).} in $Rep_k(Q_1)$.
     We circumvent this difficulty  by observing  remarkable patterns, which  the Ext-nontrivial couples obey. These patterns and the notion of
      \textit{regularity-preserving hereditary category}, which  they imply,   are other novelties of the paper.

\vspace{3mm}

\noindent 1.2. We explain now the organization of the paper and give details about the intermediate results.

Here, by  $\mc A$ we denote  a $k$-linear  $\hom$-finite hereditary  abelian category, where  $k$ is an algebraically closed field, and we denote $D^b(\mc A)$  by  $\mc T$.

In  Section \ref{non-stable exc obj in...} we analyze  the following data: an  exceptional
   object $E \in D^b(\mc A)$, which is not $\sigma$-semistable for a given stability condition $\sigma\in \st(\mc T)$.  Macr\`i initiated such an analysis in \cite[p.10]{Macri}. \\   
  We end up  in  Section \ref{non-stable exc obj in...} with   a distinguished triangle,   denoted by $\mk{alg}(E)$, which satisfies one of five possible lists of properties,
  named \textbf{C1,C2,C3,B1,B2}.  If the resulting list is one of \textbf{C1,C2} or \textbf{C3},  then we say that the object $E$ is  \textit{$\sigma$-regular},
  otherwise -  \textit{$\sigma$-irregular}.  The    triangle  $\mk{alg}(R)=\begin{diagram}[size=1.0em]
U & \rTo      &     &       &   R \\
  & \luDashto &     & \ldTo &       \\
  &           & V &       &
\end{diagram}  $   of   a $\sigma$-regular $R$  has the feature that for any indecomposable components $S$ and $E$ of $V$  and $U$, respectively,
   the  pair $(S,E)$ is exceptional with semistable  first element $S$.   We denote this relation between a $\sigma$-regular object  $R$ and the exceptional pair
   $(S,E)$ by $\bd  R & \rDotsto^{X} & (S,E)  \ed $, where $X$ contains further information as explained in Section \ref{some terminilogy}.
   This feature   is not   available in the irregular cases \textbf{B1} and \textbf{B2}, and       the obstruction  to  obtaining  it      are the  Ext-nontrivial
   couples.  Such couples exist in $Rep_k(Q_1)$ and $Rep_k(Q_2)$,
 as shown  in Section \ref{two examples}.
Essential part of our efforts concerns the Ext-nontrivial couples.
 It follows from   \cite[Lemma 4.1]{Macri} that there are not
 such couples  in $Rep_k(K(l))$ (Appendix \ref{there are no Ext-nontrivial...}).

Thus, in  Sections \ref{non-stable exc obj in...},   \ref{some
terminilogy} from each $\sigma$-regular exceptional  object $R$ we
obtain  at least one exceptional pair $(S,E)$  with $\bd  R &
\rDotsto^{X} & (S,E)  \ed $.  The first component   $S$ in such a
pair is always semistable.
If the second component $E$ is not semistable, which is possible
iff $R$ is non-final as defined in Definition \ref{def final
good}, then it is natural to ask:  Is  $E$ a  $\sigma$-regular
exceptional object?

 Motivated by this question, we introduce   in Section \ref{no bad after good}  certain conditions on  the Ext-nontrivial couples of $\mc A$, which we call
   \textit{RP property 1 and RP property 2} (Subsection \ref{RP properties}), and using them we give a positive answer. We say that $\mc A$ is a
   \textit{regularity-preserving category}(Definition \ref{RP category}), when the answer is positive. RP properties 1, 2 themselves are not important for the
   rest of the paper, but that  $\mc A$ is regularity-preserving, which follows from them.

   Whence, in regularity-preserving category $\mc A$ the relation  $\bd   & \rDotsto &   \ed $ circumvents the irregular objects, and  each non-final
   $\sigma$-regular object $R$ generates  a long  sequence\footnote{By ``long'' we mean that it has at least two steps. This sequence  is not uniquely determined
   by $R$.} of the form:  \begin{gather} \label{sequence of cases in intro}  \bd[height=1.5em] R & \rDotsto^{X_1} & (S_1,E_1) & \rMapsto^{proj_2}& E_1 & \rDotsto^{X_2} & (S_2,E_2)  & \rMapsto^{proj_2}& E_2 & \rDotsto^{X_3} & (S_3,E_3)& \rMapsto^{proj_2}&  \dots \\
  &   & \dMapsto^{proj_1} &    &  &    & \dMapsto^{proj_1} &  &  &  & \dMapsto^{proj_1} &  &  & \\
   &  &  S_1              &    &  &    &  S_2              &  &  &  & S_3               &  &  &  \ed. \end{gather}  In such a sequence, which  we call an  \textit{$R$-sequence},  the exceptional objects    $S_1, S_2$, $\dots$ are all semistable, and furthermore,  if    $E_{n}$ is final for some $n$, then,
   by the very definition of a final object(Definition
   \ref{def final good}), the pair $(S_{n+1},E_{n+1})$ is semistable and exceptional. 
   
    In Section \ref{sequence} we proceed  further in direction   $\sigma$-exceptional collections   by refining on   the phases and
   the degrees of $\{S_i \}$, and  showing various situations, in which  the  vanishings   $\Hom^*(S_i,S_1)=\Hom^*(E_i,S_1)=0$ hold for $i>1$. However,
    these vanishings do not hold in each $R$-sequence.
       Nevertheless, we show  that  starting from any  $\sigma$-regular $R$ through any $R$-sequence  we  reach a final  $\sigma$-regular object
        $E_n$ for some $n\geq 1$.

       After a  careful   examination of  the final  $\sigma$-regular objects,  in Section \ref{final},  we   find that    an exceptional  pair $(S,E)$
        produced from such an object is not only semistable, but also $(S,E[-i])$ is a $\sigma$-exceptional  pair for some $i\geq 0$ (e.g., a   situation
         as: $\phi(S)=\phi(E)$, $\Hom(S,E)\neq 0$ cannot happen).

The  proofs in Sections \ref{sequence} and \ref{final} are
facilitated by the use  of a   function $\theta_\sigma :
Ob({\mathcal T}) \rightarrow \NN^{(\sigma^{ss}_{ind}/\cong)}$,
introduced in subsection \ref{comments on stab cond}. For an
object $X\in Ob(\mc T)$ the function $\theta_\sigma(X):
\sigma^{ss}_{ind}/\cong \rightarrow \NN$ indicates (with
multiplicities)  the indecomposable components of the
Harder-Narasimhan  factors of $X$.    The  relation $\bd  R &
\rDotsto^{X} & (S,E)  \ed $  implies
$\theta_\sigma(E)<\theta_\sigma(R)$ and $\theta_\sigma(R)(S)>
0$. This feature gives an upper bound  of the   lengths of all
$R$-sequences  with a fixed $R$.  It also  plays a role in
avoiding some   situations  as the mentioned in the end of the
previous paragraph.

  In Section \ref{two examples}  we classify   exceptional objects of the categories $Rep_k(Q_1)$, $Rep_k(Q_2)$.  After that we obtain tables with
    dimensions of $\Hom(X,Y)$, ${\rm Ext}^1(X,Y)$ for any two exceptional objects $X,Y$, and observe that  one of these  always vanishes.
      RP property 1 and RP property 2 follow  by a careful analysis of these tables.
 For the Ext-nontrivial couples of the  quiver $Q_1$ we observe  an additional pattern: Corollary \ref{additional RP property}, which  helps
  us further to avoid the irregular cases.  We refer to it  as the additional RP property.  It does not hold in $Q_2$. 
 In the end of Subsection \ref{sketch of proof} we obtain the  lists  of all  exceptional pairs and triples  in $Rep_k(Q_1)$.

The results before  Section \ref{constructing} contain  the
implications (the first is due to regularity-preserving): \\
$\sigma$-regular  object $\Rightarrow$ final  $\sigma$-regular
object  $\Rightarrow$  $\sigma$-exceptional pair (Corollary
\ref{sigma pairs} and Remark \ref{if there are not Ext-nontrivila
pairs}).

 In Section \ref{constructing}  we  develop      various  criteria  for existence of $\sigma$-exceptional triples in $D^b(\mc A)$,
 assuming that the  exceptional objects of $\mc A$ obey the global properties  observed for  $Rep_k(Q_1)$  in Section \ref{two examples}.\footnote{The precise assumptions are specified after Lemma \ref{another prop of exc in quiver}.}
   It is shown that  any   non-final  \textbf{C2} or \textbf{C3} object induces such a triple. Thus, if $R$ is a  \textbf{C2} or \textbf{C3} object,
    then  any  $R$-sequence of length two produces a  $\sigma$-exceptional triple. If  $R$ is a  \textbf{C1} object, then our results imply  that any
    $R$-sequence of length three is enough, but for length less or equal to  two - only under special  circumstances (Lemmas \ref{coro for C3 after C1},
    \ref{when C1 is final}, Corollary \ref{after C1}).      \\
  If $R$ is a final $\sigma$-regular object, then  we have no long $R$-sequences, they are all of length one and each of them induces  a $\sigma$-exceptional pair.
   To obtain a $\sigma$-triple in this case we apply  two ideas. The first   is   to combine the pairs coming from different $R$-sequences, which leads
   to the result that a final $\sigma$-regular object $R$ whose  Harder-Narasimhan filtration differs from  $\mk{alg}(R)$  induces  a $\sigma$-exceptional triple.
    The other idea is to  utilize  the infimum $\phi_{min}$  and the supremum $\phi_{max}$  of the set of phases of semistable exceptional objects in $\mc A$.
     More precisely, we show that a  relation  $\bd R &  \rDotsto & (S[1],E)\ed$ with a final  \textbf{C3} object $R\in\mc A$ and $\phi(S)>\phi_{min}$ induces
     a $\sigma$-triple(Corollary   \ref{from C3 to min}). There is an analogous criterion using a final \textbf{C2} object $R\in\mc A$ and $\phi_{min}$,
      shown in Corollary   \ref{from C2 to max}, but there is not  an analogue  for  final \textbf{C1} objects (Lemma \ref{when C1 is final} uses a
      non-final \textbf{C1} object and in different setting).
When $\phi_{max}-\phi_{min}>1$, we show that, if  $(S_{min},E,S_{max})$  is an exceptional triple in $\mc A$ with $S_{min}\in \mc P(\phi_{min})$ and  $S_{max}\in \mc P(\phi_{max})$,  then non-semistability of $E$ (no
matter regular or irregular) implies a $\sigma$-exceptional
triple. The last is    widely used in Subsection \ref{>1}.

   The criteria obtained in Section \ref{constructing} combined with   the lists  of the exceptional pairs and the exceptional triples of $Rep_k(Q_1)$ at our
   disposal (due to Section \ref{two examples}) turn out to be enough for the  proof of the main Theorem \ref{main theorem for Q_1 in intro}, which is demonstrated
   in Section \ref{main theo for Q_1}. The locally finiteness of the stability condition $\sigma \in \st(\mc T)$ plays an important role as well. The proof is
   divided into two  steps: $\phi_{max}-\phi_{min}>1$ and $\phi_{max}-\phi_{min}\leq 1$.

\vspace{3mm}

\noindent 1.3.  The  following  three  statements  are  proved in   \cite{DHKK}. In  the first and the second statement, $Q$ is  an acyclic Euclidean quiver:
\begin{itemize}
    \item[(a)] \cite[Corollary 3.15]{DHKK}: For each   $\sigma \in \st(D^b(Q))$  the set of semistable phases  is either finite or has  two limit points in $\SS^1$.
    \item[(b)] \cite[Corollary 3.31]{DHKK}:   For any exceptional pair $(A,B)$ in $D^b(Q)$ and any $i\in \ZZ$ holds the inequality $\hom^i(A,B)\leq 2$.
    \item[(c)] \cite[Proposition 3.32]{DHKK}:   Any connected quiver Q, which is neither Euclidean nor Dynkin, has a family of stability
conditions with phases which are dense in an arc. The proof of this fact  relies on extendability, as defined in \cite[Definition 3.25]{DHKK}, of certain stability
conditions on a subcategory of $D^b(Q)$ to the entire $D^b(Q)$ (the precise setting is described right after Theorem 3.27 in \cite{DHKK}).
\end{itemize}

    We   construct  in Subsection \ref{two limit points} stability conditions $\sigma \in \st(D^b(Q_1))$ with two limit points in $\SS^1$,    concerning (a).

 In  Remark \ref{dimension 2} we point out   exceptional pairs $(A,B)$ in $D^b(Q_1)$ with  $\hom^i(A,B)= 2$, concerning (b).

 In  Subsection \ref{subsection about macri}  we comment on the stability conditions constructed by E. Macr\`i \cite{Macri} via exceptional collections.   By slightly modifying
  the statement of  \cite[Proposition 3.17]{Macri} and  refining its proof  is obtained Proposition \ref{projection}, which provides the extendability needed in (c).

\vspace{3mm}

\noindent 1.4.   It is known  \cite{WCB1} that the Braid group  acts transitively on the exceptional collections of $Rep_k(Q_1)$. The list of these collections
shows that this action is not free (Remark \ref{braid}).
 \vspace{3mm}

\noindent 1.5.  This paper gives  a few answers, and poses many questions.    We expect that there is a proof of Theorem \ref{main theorem for Q_1 in intro},
governed by a  general principle. The notion of regularity-preserving hereditary category (Definition \ref{RP category})  should be related to this principle.
     RP property 1 and  RP property 2 are our method to prove  regualrity-preserving.
 The fact that they hold not only in $Rep_k(Q_1)$, but also in $Rep_k(Q_2)$ (Corollary \ref{RP property 1,2 and.. for Q1}) seems to be a trace of   a larger
 unexplored  picture.  We expect  that there are further non-trivial examples of  regularity-preserving  categories.

 We do not give an answer to the question: is there  a $\sigma$-exceptional quadruple for each $\sigma \in \st(D^b(Q_2))$ (the $Q_2$-analogue of  Theorem \ref{main theorem for Q_1 in intro}). We show  that $Rep_k(Q_2)$ is regularity-preserving,   and the results of  Sections \ref{sequence}, \ref{final}, and Subsection \ref{without the additional RP}  hold for  $Rep_k(Q_2)$ entirely.  These  are clues for a positive answer(see especially  Corollary \ref{between min max 1}).  In section \ref{two examples} we give  the dimensions of $\Hom(X,Y)$, ${\rm Ext}^1(X,Y)$ for any two exceptional objects $X,Y$ in $Q_2$ as well. This lays  a  ground for working on the  $Q_2$-analogue of  Theorem \ref{main theorem for Q_1 in intro}.

We expect that the results in Section \ref{two examples} and  Theorem \ref{main theorem for Q_1 in intro} can be  used  for the  study of the topology of $\st(D^b(Rep_k(Q_1)))$ further (e. g. to check its simply-connectivity).

\vspace{3mm}

 \textit{\textbf{Some notations.}} In these notes the letters ${\mathcal T}$ and $\mc A$ denote always  a triangulated category and an abelian category, respectively, linear over a field\footnote{in some sections algebraically closedness of $k$ is not important, but overall this feature  is necessary.} $k$, the shift functor  in ${\mathcal T}$ is designated by $[1]$.   We write $\Hom^i(X,Y)$ for  $\Hom(X,Y[i])$ and  $\hom^i(X,Y)$ for  $\dim_k(\Hom(X,Y[i]))$, where $X,Y\in \mc T$.  For $X,Y\in\mc  A$,  writing $\Hom^i(X,Y)$, we consider $X,Y$ as elements in   $\mc T=D^b(\mc A)$, i.e.  $\Hom^i(X,Y)={\rm Ext}^i(X,Y)$.

   We denote by ${\mc A}_{exc}$, resp. $D^b(\mc A)_{exc}$,  the set of all
    exceptional objects of  $\mc A$, resp. of  $D^b(\mc A)$.

 An abelian category $\mc A$ is said to be hereditary, if ${\rm Ext}^i(X,Y)=0$ for any two $X,Y \in \mc A$ and each $i\geq 2$.

\textit{ For an object  $X \in D^b(\mc A)$ of the form $X\cong X'[j]$, where $X'\in \mc A$ and $j \in \ZZ$, we write $\deg(X)=j$.}

  For any quiver $Q$ we write $D^b(Q)$ for $D^b(Rep_k(Q))$.
\vspace{3mm}

\textit{{\bf Acknowledgements:}}
The authors wish to express their gratitude to Igor Dolgachev,  M. Umut Isik, Maxim Kontsevich,  Alexander  Kuznetsov,   Tony Pantev
    for their  interest in this paper.

The first author wishes to express his thanks to Matthew Ballard, Dragos Deliu, David Favero, Sergey Galkin, Fabian Haiden, M. Umut Isik, Gabriel Kerr, Alexander Noll, Pranav Pandit, Victor  Przyjalkowski  for helpful educational discussions.

The authors were funded by NSF DMS 0854977 FRG, NSF DMS 0600800, NSF DMS 0652633
FRG, NSF DMS 0854977, NSF DMS 0901330, FWF P 24572 N25, by FWF P20778 and by an
ERC Grant.

\section{On the Ext-nontrivial couples of some hereditary categories} \label{two examples}
In  Sections \ref{no bad after good}, \ref{sequence}, \ref{final}, \ref{constructing}   we treat hereditary abelian categories whose exceptional objects are supposed to obey  specific pairwise relations.
    In  this section  we give examples of such categories.

\subsection{The categories} For any finite quiver $Q$ and an algebraically closed field $k$ we denote   the category of $k$-representations
of $Q$ by $Rep_{k}(Q)$.  It is well known that $Rep_{k}(Q)$ is a $\hom$-finite hereditary $k$-linear abelian category (see e. g. \cite{WCB2}).

In this section we classify  the exceptional objects of the categories of representations of the following quivers:
\be \label{Q1} Q_1= \begin{diagram}[1.5em]
   &       &  \circ  &       &    \\
   & \ruTo &    & \luTo &       \\
\circ  & \rTo  &    &       &  \circ
\end{diagram}  \ \qquad \qquad \ Q_2= \begin{diagram}[1.5em]
 \circ &  \rTo  &  \circ    \\
  \uTo &        & \uTo     \\
\circ   & \rTo  &    \circ
\end{diagram}. \ee
After that we compute  the dimensions of $\Hom(X,Y)$, ${\rm Ext}^1(X,Y)$ for any two exceptional objects $X,Y$.  The obtained information reveals
 some patterns, which are of   importance for the rest of the paper.

  More precisely, Corollary \ref{RP property 1,2 and.. for Q1} \textbf{(a)} claims  that $Rep_k(Q_1)$ and  $Rep_k(Q_2)$ have RP property 1 and RP property 2(see subsection \ref{RP properties} for definition). These properties ensure that $Rep_k(Q_1)$ and  $Rep_k(Q_2)$ are regularity-preserving(Definition \ref{RP category}, Proposition \ref{prop no bad after good}), which  is of primary importance for  Sections  \ref{sequence}, \ref{constructing}, \ref{main theo for Q_1}.

In the end of Section \ref{sequence} and in Section
\ref{constructing},  the property that for any two exceptional
objects $X,Y$ at most one of the spaces $\Hom(X,Y)$, ${\rm
Ext}^1(X,Y)$ is nonzero  plays an important role.    Corollary
\ref{RP property 1,2 and.. for Q1} \textbf{(b)} asserts that  this
property holds for both the quivers  $Q_1$, $Q_2$.

For   $Q_1$ we observe  the additional RP property(see Corollary
\ref{additional RP property}),  used in Subsection \ref{with the
additional RP property}. In the end we obtain  the lists of
exceptional pairs and exceptional triples in $Rep_k(Q_1)$, which
are  widely used in Section \ref{main theo for Q_1}.

We give now more details.
\subsection{The dimensions \texorpdfstring{ $\hom(X,Y),\hom^1(X,Y)$}{\space} for \texorpdfstring{ $X,Y\in Rep_k(Q_i)_{exc}$ and $i\in \{1,2\}$ }{\space} } \label{sketch of proof} \mbox{} \\

For a representation $\rho= \begin{diagram}[size=1.2em]
 k^{\alpha_+} &  \rTo  & k^{\alpha_e}   \\
  \uTo &        & \uTo     \\
k^{\alpha_b}  & \rTo  &   k^{\alpha_-}
\end{diagram}$ $\in Rep_k(Q_2)$, where $\alpha_b,\alpha_-,\alpha_+,\alpha_e \in \NN$, we denote its dimension vector by $\ul{\dim}(\rho)=(\alpha_b,\alpha_-,\alpha_+,\alpha_e)$ and for a representation $\begin{diagram}[size=1em]
   &       &  k^{\alpha_e} &       &    \\
   & \ruTo &    & \luTo &       \\
k^{\alpha_b}  & \rTo  &    &       &  k^{\alpha_{mid}}
\end{diagram} = \rho \in Rep_k(Q_1)$ we denote $\ul{\dim}(\rho)=(\alpha_b,\alpha_{mid},\alpha_e)$.
The Euler forms of $Q_1, Q_2$  are:
\begin{gather} \left\langle(\alpha_b,\alpha_{mid},\alpha_e), \  (\alpha_b',\alpha_{mid}',\alpha_e') \right \rangle = \alpha_b \alpha_b' + \alpha_{mid} \alpha_{mid}' + \alpha_e \alpha_e' - \alpha_b \alpha_e' - \alpha_b \alpha_{mid}' - \alpha_{mid} \alpha_e' ,  \nonumber \\  \left\langle (\alpha_b,\alpha_-,\alpha_+,\alpha_e), \ (\alpha_b',\alpha_-',\alpha_+',\alpha_e')\right \rangle =\begin{array}{c} \alpha_+ \alpha_+'+\alpha_- \alpha_-'+\alpha_b \alpha_b'+\alpha_e \alpha_e' \\ -\alpha_b \alpha_+'-\alpha_b \alpha_-'-\alpha_+ \alpha_e'-\alpha_- \alpha_e' \end{array}. \nonumber  \end{gather}
Recall(see page 8 in \cite{WCB2}) that for  any $\rho,\rho'\in Rep_k(Q)$ we have the formula
\begin{gather}\label{euler} \hom(\rho,\rho')-\hom^1(\rho,\rho')=\scal{\ul{\dim}(\rho),\ul{\dim}(\rho')}.  \end{gather}
In particular, it follows  that if $\rho \in Rep_k(Q)$ is an exceptional object, then $\scal{\ul{\dim}(\rho),\ul{\dim}(\rho)}=1$. The vectors satisfying this equality are called real roots(see \cite[p. 17]{WCB2}).    For example, one can show that the real roots of $Q_1$ are $(m+1,m,m)$,$(m,m+1,m+1)$, $(m,m,m+1)$, $(m
 +1,m+1,m)$, $(m+1,m,m+1)$, $(m,m+1,m)$, $m\geq 0$. The imaginary roots\footnote{Imaginary root is a  vector $\rho$ with $\scal{\ul{dim}(\rho),\ul{dim}(\rho)}\leq 0$.}  of $Q_1$,  are $(m,m,m)$, $m\geq 1$. Not every real root is a dimension vector of an exceptional representation. More precisely:
\begin{lemma} \label{roots no exceptional} Let $m\geq 1$. If  $(\alpha_b,\alpha_{mid},\alpha_e) \in \{ (m+1,m,m+1), (m,m+1,m) \}_{m\in \NN}$, then   $(\alpha_b,\alpha_{mid},\alpha_e)$ is not dimension vector of any exceptional representation in $Rep_k(Q_1)$. If
 $(\alpha_b,\alpha_-,\alpha_+, \alpha_e) \in \{ (m,m+1,m,m)$, $(m,m,m+1,m)$, $(m+1,m,m+1,m+1)$, $(m+1,m+1,m,m+1) \}_{m\in \NN}$, then  $(\alpha_b,\alpha_-,\alpha_+, \alpha_e)$ is not dimension vector of any exceptional representation in $Rep_k(Q_2)$.
\end{lemma}
\textit{Sketch of proof.} For the proof of this lemma one can use   (see \cite[Lemma 1 on page 13]{WCB2}) that  a representation  $\rho \in Rep_k(Q_i)$ is without self-extensions iff
$\dim({\mc O}_{\rho})=\dim(Rep_k(Q_i))$, where ${\mc O}_{\rho}$ is the orbit of $\rho$ in $Rep_k(Q_i)$ as defined in  \cite[page 11,12]{WCB2}.
Using this argument, it can be shown  that any representation without self-extensions with dimension vector among the listed in the lemma is decomposable.
\qed

Now we classify the exceptional objects on $Rep_k(Q_1)$,
$Rep_k(Q_2)$(Propositions \ref{exceptional objects in Q1} and
\ref{exceptional objects in Q2}). In these propositions we use the
following notations  for any $m\geq 1$:
\begin{gather} \nonumber \pi_+^m: k^{m+1} \rightarrow k^{m}, \quad  \pi_-^m: k^{m+1} \rightarrow k^{m}, \quad j_+^m: k^{m} \rightarrow k^{m+1}, \quad  j_-^m: k^{m} \rightarrow k^{m+1}  \\
\nonumber  \pi_+^m(a_1,a_2,\dots, a_m, a_{m+1}) =(a_1,a_2,\dots, a_m) \qquad   \pi_-^m(a_1,a_2,\dots, a_m, a_{m+1})=(a_2,\dots, a_m, a_{m+1}) \\
\nonumber  j_+^m(a_1,a_2,\dots, a_m) =(a_1,a_2,\dots, a_m,0)  \qquad
  j_-^m(a_1,a_2,\dots, a_m)=(0,a_1,\dots, a_m).
\end{gather}

\begin{prop} \label{exceptional objects in Q1} The exceptional objects up to isomorphism in  $ Rep_{k}(Q_1) $ are ($m=0,1,2,\dots$)
\begin{gather} E_1^m = \begin{diagram}[1.5em]
   &       &  k^m &       &    \\
   & \ruTo^{\pi_+^m} &    & \luTo^{Id} &       \\
k^{m+1}  & \rTo^{\pi_-^m}  &    &       &  k^m
\end{diagram} \ \ \ \  E_2^m = \begin{diagram}[1.5em]
   &       &  k^{m+1} &       &    \\
   & \ruTo^{j_+^m} &    & \luTo^{Id} &       \\
k^{m}  & \rTo^{j_-^m}  &    &       &  k^{m+1}
\end{diagram} \ \ \ \  E_3^m = \begin{diagram}[1.5em]
   &       &  k^{m+1} &       &    \\
   & \ruTo^{j_+^m} &    & \luTo^{j_-^m} &       \\
k^{m}  & \rTo^{Id}  &    &       &  k^{m}
\end{diagram} \nonumber  \\
E_4^m = \begin{diagram}[1.5em]
   &       &  k^m &       &    \\
   & \ruTo^{\pi_+^m} &    & \luTo^{\pi_-^m} &       \\
k^{m+1}  & \rTo^{Id}  &    &       &  k^{m+1}
\end{diagram} \ \ \ \ M = \begin{diagram}[1.5em]
   &       &  0 &       &    \\
   & \ruTo &    & \luTo &       \\
0  & \rTo  &    &       &  k
\end{diagram} \ \ \ \  M'= \begin{diagram}[1.5em]
   &       &  k &       &    \\
   & \ruTo^{Id} &   &  \luTo  &       \\
k  &  \rTo &    &       &  0
\end{diagram}.\nonumber \end{gather}
\end{prop}
\textit{Sketch of proof.} We showed that the dimension vectors of the exceptional representations are real roots. The list of real roots is given before Lemma \ref{roots no exceptional} and  some of them are excluded in Lemma \ref{roots no exceptional}.   Moreover,  there is at most one representation without self-extensions of a given dimension vector up to isomorphism \cite[p. 13]{WCB2}. Taking into account these arguments, the proposition follows   by   showing  that the endomorphism space of each of the listed representations is $k$ (recall also \eqref{euler}). The computations, which we skip, are  reduced to  table \eqref{vect space table} in Appendix \ref{table with matrices}. \qed

\begin{prop} \label{exceptional objects in Q2} The exceptional objects up to isomorphism in  $ Rep_{k}(Q_2) $ are($m=0,1,2,\dots$)
\begin{gather} E_1^m = \begin{diagram}[1.5em]
 k^m &  \rTo^{Id}  &  k^m   \\
  \uTo^{\pi_+^m} &        & \uTo^{Id}     \\
k^{m+1}  & \rTo^{\pi_-^m}  &    k^m
\end{diagram} \ \ \ \  E_2^m = \begin{diagram}[1.5em]
 k^{m+1} &  \rTo^{Id}  &  k^{m+1}   \\
  \uTo^{j_+^m} &        & \uTo^{Id}     \\
k^{m}  & \rTo^{j_-^m}   &    k^{m+1}
\end{diagram} \ \ \ \  E_3^m = \begin{diagram}[1.5em]
 k^{m} &  \rTo^{j_+^m}  &  k^{m+1}   \\
  \uTo^{Id} &        & \uTo^{j_-^m}     \\
k^{m}  & \rTo^{Id}   &    k^{m}
\end{diagram}\ \ \ \  E_4^m = \begin{diagram}[1.5em]
 k^{m+1} &  \rTo^{\pi_+^m}  &  k^{m}   \\
  \uTo^{Id} &        & \uTo^{\pi_-^m}     \\
k^{m+1}  & \rTo^{Id}   &    k^{m+1}
\end{diagram}\nonumber  \end{gather}
\begin{gather}
E_5^m = \begin{diagram}[1.5em]
 k^m &  \rTo^{j_+^m}  &  k^{m+1}   \\
  \uTo^{Id} &        & \uTo^{Id}     \\
k^{m}  & \rTo^{j_-^m}  &    k^{m+1}
\end{diagram} \ \ \ \  E_6^m = \begin{diagram}[1.5em]
 k^{m+1} &  \rTo^{\pi_+^m}  &  k^{m}   \\
  \uTo^{Id} &        & \uTo^{Id}     \\
k^{m+1}  & \rTo^{\pi_-^m}   &    k^{m}
\end{diagram} \ \ \ \  E_7^m = \begin{diagram}[1.5em]
 k^{m} &  \rTo^{Id}  &  k^{m}   \\
  \uTo^{\pi_+^m} &        & \uTo^{\pi_-^m}     \\
k^{m+1}  & \rTo^{Id}   &    k^{m+1}
\end{diagram}\ \ \ \  E_8^m = \begin{diagram}[1.5em]
 k^{m+1} &  \rTo^{Id}  &  k^{m+1}   \\
  \uTo^{j_+^m} &        & \uTo^{j_-^m}     \\
k^{m}  & \rTo^{Id}   &    k^{m}
\end{diagram}\nonumber \end{gather}

\begin{gather} F_+ = \begin{diagram}[1.5em]
 k &  \rTo  &  0   \\
  \uTo &        & \uTo     \\
0  & \rTo  &    0
\end{diagram} \ \ \ \  F_- = \begin{diagram}[1.5em]
 0 &  \rTo  &  0   \\
  \uTo &        & \uTo     \\
0  & \rTo   &    k
\end{diagram} \ \ \ \  G_+ = \begin{diagram}[1.5em]
 k &  \rTo^{Id}  &  k   \\
  \uTo^{Id} &        & \uTo    \\
k  & \rTo   &    0
\end{diagram}\ \ \ \  G_- = \begin{diagram}[1.5em]
 0 &  \rTo  &  k   \\
  \uTo &        & \uTo^{Id}     \\
k  & \rTo^{Id}   &    k
\end{diagram}.\nonumber \end{gather}
\end{prop}
\textit{Sketch of proof.} The same as  Proposition \ref{exceptional objects in Q1}. \qed

Now we  compute  $\hom(\rho,\rho')$, $\hom^1(\rho,\rho')$ with  $\rho,\rho'$ varying throughout the obtained lists.
\begin{prop} \label{prop for Q1 table}
The dimensions of the vector spaces $\Hom(X,Y) $ and $\Hom^1(X,Y) $ for any pair of exceptional objects $X,Y \in Rep_{k}(Q_1)$ are contained in the following table:
\begin{gather} \label{Q1 table}\begin{array}{| c | c | c | c | c | c | c |}
  \hline
                      &                &   \hom  &     \hom^1 &               & \hom    &   \hom^1       \\ \hline
  0\leq m < n        & (E_1^m,E_1^n)  &   0    & n-m-1   & (E_1^n,E_1^m) & 1+n-m  & 0           \\ \hline
  0\leq n < m        & (E_2^m,E_2^n)  &   0    & m-n-1   & (E_2^n,E_2^m) & 1+m-n  & 0           \\ \hline
  0\leq n < m        & (E_3^m,E_3^n)  &   0    & m-n-1   & (E_3^n,E_3^m) & 1+m-n  & 0           \\ \hline
  0\leq m < n        & (E_4^m,E_4^n)  &   0    & n-m-1   & (E_4^n,E_4^m) & 1+n-m  & 0           \\ \hline
    m\geq 0, n \geq 0  & (E_1^m,E_2^n)  &   0    & n+m+2   & (E_2^n,E_1^m) & n+m    & 0           \\ \hline
  m\geq 0, n \geq 0  & (E_1^m,E_3^n)  &   0    & n+m+1   & (E_3^n,E_1^m) & n+m    & 0           \\ \hline
  0\leq m \leq n     & (E_1^m,E_4^n)  &   0    & n-  m   & (E_4^n,E_1^m) & 1+n-m  & 0           \\ \hline
  0\leq n <    m     & (E_1^m,E_4^n)  &   m-n  & 0       & (E_4^n,E_1^m) & 0      & m-n-1       \\ \hline
  0\leq n \leq m     & (E_2^m,E_3^n)  &   0    & m-n     & (E_3^n,E_2^m) & 1+m-n  & 0           \\ \hline
  0\leq m <n         & (E_2^m,E_3^n)  &   n-m  & 0       & (E_3^n,E_2^m) & 0      & n-m-1       \\ \hline
  m \geq 0, n \geq 0 & (E_2^m,E_4^n)  & 1+ n+m & 0       & (E_4^n,E_2^m) & 0      & n+m+2       \\ \hline
  m \geq 0, n \geq 0 & (E_3^m,E_4^n)  &    n+m & 0       & (E_4^n,E_3^m) & 0      & n+m+2       \\ \hline
  m \geq 0           & (M,E_1^m    )  & 0      & 0       & (E_1^m,M)     & 0      & 1           \\ \hline
  m \geq 0           & (M,E_2^m    )  & 0      & 0       & (E_2^m,M)     & 1      & 0           \\ \hline
  m \geq 0           & (M,E_3^m    )  & 0      & 1       & (E_3^m,M)     & 0      & 0           \\ \hline
  m \geq 0           & (M,E_4^m    )  & 1      & 0       & (E_4^m,M)     & 0      & 0           \\ \hline
  m \geq 0           & (M',E_1^m   )  & 1      & 0       & (E_1^m,M')    & 0      & 0           \\ \hline
  m \geq 0           & (M',E_2^m   )  & 0      & 1       & (E_2^m,M')    & 0      & 0           \\ \hline
  m \geq 0           & (M',E_3^m   )  & 0      & 0       & (E_3^m,M')    & 1      & 0           \\ \hline
  m \geq 0           & (M',E_4^m   )  & 0      & 0       & (E_4^m,M')    & 0      & 1           \\ \hline
                     & (M,M'   )      & 0      & 1       & (M',M)        & 0      & 1           \\ \hline
  \end{array}
  \end{gather}
\end{prop}
\textit{Sketch of proof.} Via    computations, which we do not write out here, we obtain   $\hom(\rho,\rho')$ for any two representations $\rho,\rho'$ taken from  Proposition \ref{exceptional objects in Q1}.
The computations are  reduced to determining  the dimensions of
some vector spaces of matrices. These spaces and their dimensions
are listed   in  Appendix \ref{table with matrices}, table \eqref{vect space table}. Having $\hom(\rho,\rho')$,
 the dimension $\hom^1(\rho,\rho')$ is  computed by \eqref{euler}.
\qed
\begin{prop} \label{prop for Q2 table}
The dimensions $\hom(X,Y) $ and $\hom^1(X,Y) $ for any pair of exceptional objects $X,Y \in Rep_{k}(Q_2)$ are contained in the following table:
\tiny
\begin{gather}  \begin{array}{| c | c | c | c | c | c | c |}
  \hline
                     &                &   \hom  &     \hom^1 &               & \hom    &   \hom^1       \\ \hline
  0\leq n < m        & (E_1^m,E_1^n)  &   1+m-n    & 0   & (E_1^n,E_1^m) & 0  & m-n-1           \\ \hline
   0\leq m < n        & (E_2^m,E_2^n)  &   1+n-m    & 0   & (E_2^n,E_2^m) & 0  & n-m-1           \\ \hline
   0\leq m < n        & (E_3^m,E_3^n)  &   1+n-m    & 0   & (E_3^n,E_3^m) & 0  & n-m-1           \\ \hline
   0\leq n < m        & (E_4^m,E_4^n)  &   1+m-n    & 0   & (E_4^n,E_4^m) & 0  & m-n-1           \\ \hline
      0\leq m < n        & (E_5^m,E_5^n)  &   1+n-m    & 0   & (E_5^n,E_5^m) & 0  & n-m-1           \\ \hline
    0\leq n < m        & (E_6^m,E_6^n)  &   1+m-n    & 0   & (E_6^n,E_6^m) & 0  & m-n-1           \\ \hline
   0\leq n < m        & (E_7^m,E_7^n)  &   1+m-n    & 0   & (E_7^n,E_7^m) & 0  & m-n-1           \\ \hline
    0\leq m < n        & (E_8^m,E_8^n)  &   1+n-m    & 0   & (E_8^n,E_8^m) & 0  & n-m-1           \\ \hline
0\leq m, 0 \leq n     & (E_1^m,E_2^n)  &   0    & 2+n+m   & (E_2^n,E_1^m) & m+n  & 0             \\ \hline
0\leq m, 0 \leq n     & (E_1^m,E_3^n)  &   0    &   n+m   & (E_3^n,E_1^m) & m+n  & 0             \\ \hline
      0\leq n<m          & (E_1^m,E_4^n)  &   m-n-1 &   0    & (E_4^n,E_1^m) & 0    & m-n-1          \\ \hline
    0\leq m\leq n      & (E_1^m,E_4^n)  &   0     &  n-m+1  & (E_4^n,E_1^m)& n-m+1   & 0          \\ \hline
       0\leq m, 0 \leq n     & (E_1^m,E_5^n)  &   0     &  n+m+1  & (E_5^n,E_1^m)& n+m     & 0          \\ \hline
         0\leq n< m        & (E_1^m,E_6^n)  &   m-n     &  0    & (E_6^n,E_1^m)& 0    & m-n-1          \\ \hline
       0\leq m\leq n     & (E_1^m,E_6^n)  &   0     &  n-m    & (E_6^n,E_1^m)& n-m+1    &   0        \\ \hline
     0\leq n<m         & (E_1^m,E_7^n)  &   m-n     &  0    & (E_7^n,E_1^m)& 0    &   m-n-1      \\ \hline
     0\leq m\leq n      & (E_1^m,E_7^n)  &   0      &  n-m  & (E_7^n,E_1^m)& n-m+1    &   0         \\ \hline
 0\leq m, 0 \leq n     & (E_1^m,E_8^n)  &   0     &  n+m+1  & (E_8^n,E_1^m)& n+m     & 0          \\ \hline
 \end{array} \nonumber
  \end{gather}

\begin{gather}  \begin{array}{| c | c | c | c | c | c | c |}
  \hline
                       &                &   \hom  &     \hom^1 &               & \hom    &   \hom^1       \\ \hline
    0\leq n\leq m      & (E_2^m,E_3^n)  &   0      & m-n+1  & (E_3^n,E_2^m)& m-n+1    &   0         \\ \hline
             0\leq m< n        & (E_2^m,E_3^n)  &   n-m-1      & 0  & (E_3^n,E_2^m)& 0    &   n-m-1         \\ \hline
     0\leq m, 0 \leq n     & (E_2^m,E_4^n)  &   2+m+n     &  0  & (E_4^n,E_2^m)& 0     & n+m+2         \\ \hline
0\leq m< n     & (E_2^m,E_5^n)  &   n-m     &  0  & (E_5^n,E_2^m)& 0     & n-m-1         \\ \hline
0\leq n\leq m     & (E_2^m,E_5^n)  &   0    &  m-n   & (E_5^n,E_2^m)& m-n+1      &    0     \\ \hline
0\leq m, 0 \leq n     & (E_2^m,E_6^n)  &   1+m+n     &  0  & (E_6^n,E_2^m)& 0     & n+m+2         \\ \hline
0\leq m, 0 \leq n     & (E_2^m,E_7^n)  &   1+m+n     &  0  & (E_7^n,E_2^m)& 0     & n+m+2         \\ \hline
   0\leq m< n     & (E_2^m,E_8^n)  &   n-m    &  0   & (E_8^n,E_2^m)& 0      &    n-m-1     \\ \hline
   0\leq n\leq m      & (E_2^m,E_8^n)  &   0      & m-n    & (E_8^n,E_2^m)& m-n+1    &   0         \\ \hline
0\leq m, 0 \leq n     & (E_3^m,E_4^n)  &   m+n     &  0  & (E_4^n,E_3^m)& 0     & n+m+2         \\ \hline
   0\leq m\leq n      & (E_3^m,E_5^n)  &   n-m+1      &  0  & (E_5^n,E_3^m)& 0    &   n-m         \\ \hline
   0\leq n<m         & (E_3^m,E_5^n)  &   0     &  m-n-1    & (E_5^n,E_3^m)& m-n    &   0         \\ \hline
   0\leq m, 0 \leq n     & (E_3^m,E_6^n)  &   m+n     &  0  & (E_6^n,E_3^m)& 0     & n+m+1         \\ \hline
0\leq m, 0 \leq n     & (E_3^m,E_7^n)  &   m+n     &  0  & (E_7^n,E_3^m)& 0     & n+m+1         \\ \hline
   0\leq m\leq n      & (E_3^m,E_8^n)  &   n-m+1      &  0  & (E_8^n,E_3^m)& 0    &   n-m         \\ \hline
   0\leq n<m         & (E_3^m,E_8^n)  &   0     &  m-n-1    & (E_8^n,E_3^m)& m-n    &   0         \\ \hline
0\leq m, 0 \leq n     & (E_4^m,E_5^n)  &   0     &  2+m+n  & (E_5^n,E_4^m)& 1+m+n     & 0         \\ \hline
  0\leq m< n     & (E_4^m,E_6^n)  &   0    &  n-m-1   & (E_6^n,E_4^m)& n-m      &     0    \\ \hline
   0\leq n\leq m      & (E_4^m,E_6^n)  &  1+ m-n      &   0  & (E_6^n,E_4^m)& 0    &   m-n         \\ \hline
  0\leq m< n     & (E_4^m,E_7^n)  &   0    &  n-m-1   & (E_7^n,E_4^m)& n-m      &     0    \\ \hline
   0\leq n\leq m     & (E_4^m,E_7^n)  &   m-n+1    &   0  & (E_7^n,E_4^m)& 0      &     m-n    \\ \hline
     0\leq m, 0 \leq n     & (E_4^m,E_8^n)  &   0     &  2+m+n  & (E_8^n,E_4^m)& 1+m+n     & 0         \\ \hline
0\leq m, 0 \leq n     & (E_5^m,E_6^n)  &   m+n     &  0  & (E_6^n,E_5^m)& 0     & 2+m+n         \\ \hline
0\leq m, 0 \leq n     & (E_5^m,E_7^n)  &   1+m+n     &  0  & (E_7^n,E_5^m)& 0     & 1+m+n         \\ \hline
 0\leq m\leq n      & (E_5^m,E_8^n)  &   n-m       &  0  & (E_8^n,E_5^m)& 0    &   n-m         \\ \hline
   0\leq n\leq m     & (E_5^m,E_8^n)  &   0     &   m-n  & (E_8^n,E_5^m)& m-n      &     0    \\ \hline
   0\leq n\leq m     & (E_6^m,E_7^n)  &   m-n     &   0  & (E_7^n,E_6^m)& 0     &     m-n     \\ \hline
  0\leq m\leq n     & (E_6^m,E_7^n)  &   0     &   n-m  & (E_7^n,E_6^m)& n-m     &     0     \\ \hline
0\leq m, 0 \leq n     & (E_6^m,E_8^n)  &   0    &  1+m+n   & (E_8^n,E_6^m)& 1+m+n     & 0         \\ \hline
0\leq m, 0 \leq n     & (E_7^m,E_8^n)  &   0    &  2+m+n   & (E_8^n,E_7^m)&   m+n     & 0         \\ \hline
 0\leq m     & (F_+,E_1^m)  &   0    &   0   & (E_1^m,F_+) & 0  & 1             \\ \hline
0\leq m     & (F_-,E_1^m)  &   0    &   0   & (E_1^m,F_-) & 0  & 1             \\ \hline
0\leq m     & (F_+,E_2^m)  &   0    &   0   & (E_2^m,F_+) & 1  & 0             \\ \hline
0\leq m     & (F_-,E_2^m)  &   0    &   0   & (E_2^m,F_-) & 1  & 0             \\ \hline
0\leq m     & (F_+,E_3^m)  &   0    &   1   & (E_3^m,F_+) & 0  & 0             \\ \hline
0\leq m     & (F_-,E_3^m)  &   0    &   1   & (E_3^m,F_-) & 0  & 0             \\ \hline
0\leq m     & (F_+,E_4^m)  &   1    &   0   & (E_4^m,F_+) & 0  & 0             \\ \hline
0\leq m     & (F_-,E_4^m)  &   1    &   0   & (E_4^m,F_-) & 0  & 0             \\ \hline
0\leq m     & (F_+,E_5^m)  &   0    &   1   & (E_5^m,F_+) & 0  & 0             \\ \hline
0\leq m     & (F_-,E_5^m)  &   0    &   0   & (E_5^m,F_-) & 1  & 0             \\ \hline
0\leq m     & (F_+,E_6^m)  &   1    &   0   & (E_6^m,F_+) & 0  & 0             \\ \hline
0\leq m     & (F_-,E_6^m)  &   0    &   0   & (E_6^m,F_-) & 0  & 1             \\ \hline
0\leq m     & (F_+,E_7^m)  &   0    &   0   & (E_7^m,F_+) & 0  & 1             \\ \hline
0\leq m     & (F_-,E_7^m)  &   1    &   0   & (E_7^m,F_-) & 0  & 0             \\ \hline
0\leq m     & (F_+,E_8^m)  &   0    &   0   & (E_8^m,F_+) & 1  & 0             \\ \hline
0\leq m     & (F_-,E_8^m)  &   0    &   1   & (E_8^m,F_-) & 0  & 0             \\ \hline
0\leq m     & (G_\pm,E_1^m)  &   1    &   0   & (E_1^m,G_\pm) & 0  & 0             \\ \hline
0\leq m     & (G_\pm,E_2^m)  &   0    &   1   & (E_2^m,G_\pm) & 0  & 0             \\ \hline
0\leq m     & (G_\pm,E_3^m)  &   0    &   0   & (E_3^m,G_\pm) & 1  & 0             \\ \hline
0\leq m     & (G_\pm,E_4^m)  &   0    &   0   & (E_4^m,G_\pm) & 0  & 1             \\ \hline
0\leq m     & (G_+,E_5^m)  &   0    &   1   & (E_5^m,G_+) & 0  & 0             \\ \hline
0\leq m     & (G_-,E_5^m)  &   0    &   0   & (E_5^m,G_-) & 1  & 0             \\ \hline
0\leq m     & (G_+,E_6^m)  &   1    &   0   & (E_6^m,G_+) & 0  & 0             \\ \hline
0\leq m     & (G_-,E_6^m)  &   0    &   0   & (E_6^m,G_-) & 0  & 1             \\ \hline
0\leq m     & (G_+,E_7^m)  &   0    &   0   & (E_7^m,G_+) & 0  & 1             \\ \hline
0\leq m     & (G_-,E_7^m)  &   1    &   0   & (E_7^m,G_-) & 0  & 0             \\ \hline
0\leq m     & (G_+,E_8^m)  &   0    &   0   & (E_8^m,G_+) & 1  & 0             \\ \hline
0\leq m     & (G_-,E_8^m)  &   0    &   1   & (E_8^m,G_-) & 0  & 0             \\ \hline
           & (F_+,F_-)  &   0    &   0   & (F_-,F_+) & 0  & 0             \\ \hline
            & (F_+,G_+)  &   0    &   0   & (G_+,F_+) & 0  & 0             \\ \hline
            & (F_+,G_-)  &   0    &   1   & (G_-,F_+) & 0  & 1             \\ \hline
          & (F_-,G_+)  &   0    &   1   & (G_+,F_-) & 0  & 1             \\ \hline
            & (F_-,G_-)  &   0    &   0   & (G_-,F_-) & 0  & 0             \\ \hline
            & (G_+,G_-)  &   0    &   0   & (G_-,G_+) & 0  & 0             \\ \hline
            \end{array} \nonumber
  \end{gather}
\normalsize
\end{prop}
\textit{Sketch of proof.}
The table for  $Rep_k(Q_2)$  is obtained by the same method as  for $Rep_k(Q_1)$.
\qed

The next subsection contains corollaries of the  obtained tables. 
\subsection{The Ext-nontrivial couples and their properties}
From the table in   Proposition \ref{prop for Q1 table}  we see
that the only couple  $\{X,Y\}$ of exceptional objects in
$Rep_k(Q_1)$ satisfying $\hom^1(X,Y)\neq 0$ and $\hom^1(Y,X)\neq
0$ is $\{M,M'\}$ . We call such a couple an \textit{Ext-nontrivial
couple} (see Definition \ref{Ext-nontrivial couple}). By
Proposition \ref{prop for Q2 table} we see that  the
Ext-nontrivial couples in $Rep_k(Q_2)$ are  $\{F_+,G_-\}$,
$\{F_-,G_+\}$.

Corollary \ref{RP property 1,2 and.. for Q1} concerns both
$Rep_k(Q_1)$ and $Rep_k(Q_2)$.
\begin{coro} \label{RP property 1,2 and.. for Q1} The categories $Rep_{k}(Q_1)$, $Rep_{k}(Q_2)$ satisfy the following properties:
\begin{itemize}
    \item[\textbf{(a)}] RP property 1, RP property 2 (see subsection \ref{RP properties} for description).
    \item[\textbf{(b)}] For any two exceptional objects $X, Y \in Rep_{k}(Q_i)$ at most one  degree in $\{ \hom^p(X,Y) \}_{p\in \ZZ}$ is nonzero, where  $i\in \{1,2\}$.
    \end{itemize}
\end{coro}
\bpr It follows by a careful  case by case check, using the tables in Propositions \ref{prop for Q1 table}, \ref{prop for Q2 table}. \epr
The following four corollaries concern only $Rep_{k}(Q_1)$  and are contained in table \eqref{Q1 table}.
\begin{coro} \label{additional RP property}
     If $\{\Gamma_1,\Gamma_2\}$ is an Ext-nontrivial couple  in $Rep_{k}(Q_1)$(see Definition \ref{Ext-nontrivial couple}), then for each exceptional object $X\in  Rep_{k}(Q_1)$  we have $\hom^p(\Gamma_i,X)\neq 0$ for some $i\in \{1,2\}$, $p\in \ZZ$   and $ \hom^q(X,\Gamma_j) $ for some $j\in \{1,2\}$, $q\in \ZZ$.
\end{coro}

\begin{coro}\label{coro Q1 exceptional pairs} The exceptional pairs $(X,Y)$  in $Rep_{k}(Q_1)$  are ($m \in \NN$):  \begin{gather} (E_1^{m+1}, E_1^m)\ \  (E_2^m, E_2^{m+1})\ \ (E_3^m, E_3^{m+1})\ \ (E_4^{m+1}, E_4^m)\ \ (E_1^0, E_2^0)\ \ (E_1^0, E_3^0) \nonumber \\
\label{Q1 exceptional pairs}  (E_4^{m}, E_1^m)\ \ (E_1^{m+1}, E_4^m)\ \ (E_3^{m}, E_2^m)\ \ (E_2^{m}, E_3^{m+1})\ \ (E_4^{0}, E_3^0)\ \ (E_1^{m}, M) \\
 (E_2^{m}, M)\ \  (M, E_3^{m})\ \ (M, E_4^{m})\ \ (M', E_1^{m})\ \ (M', E_2^{m})\ \ ( E_3^{m},M')\ \ ( E_4^{m},M'). \nonumber \end{gather}

\end{coro}
Using this corollary  we  obtain the list of the exceptional
triples of $Rep_{k}(Q_1)$, which by \cite{WCB1} are  the full
exceptional collections.
\begin{coro} \label{exceptional colleections}  The full  exceptional collections in $Rep_{k}(Q_1)$ up to isomorphism are ($m \in \NN$):
\ben \begin{array}{c  c  c} (E_1^{m+1},E_1^m,M) & (E_1^{m+1}, E_4^m,E_1^m) & (E_1^{m+1} , M, E_4^m) \\
(E_1^0,E_2^0, M) & (E_1^0,E_3^0,E_2^0) & (E_1^{0},M,E_3^0) \\
(E_2^m,E_2^{m+1},M) & (E_2^m,E_3^{m+1},E_2^{m+1}) & (E_2^{m},M,E_3^{m+1})\\
(E_3^m,E_2^{m},E_3^{m+1}) & (E_3^m,E_3^{m+1},M') & (E_3^{m},M',E_2^{m})\\
(E_4^{m+1},E_4^{m},M') & (E_4^{m+1},E_1^{m+1},E_4^{m}) & (E_4^{m+1},M',E_1^{m+1})\\
(E_4^{0},E_1^{0},E_3^0) & (E_4^{0},E_3^{0},M') & (E_4^{0},M',E_1^{0})\\
(M,E_3^{m},E_3^{m+1}) & (M,E_4^{m+1},E_4^m) & (M,E_4^{0},E_3^{0}) \\
(M',E_1^{m+1},E_1^{m}) & (M',E_2^{m},E_2^{m+1}) & (M',E_1^{0},E_2^{0}).\end{array}
\een
\end{coro}
The following corollary is a special case of a result in \cite{WCB1}.  It also follows   from Corollary \ref{exceptional colleections}.

\begin{coro} \label{coro for isom triples} Let $(A_0,A_1,A_2)$,$(A_0',A_1',A_2')$ be two exceptional triples in $ Rep_{k}(Q_1)$.
  If  $A_i \cong A_i'$, $A_j \cong A_j'$ for two different $i,j\in \{0,1,2\}$, then $A_k \cong A_k'$ for the third $k\in \{0,1,2\}$.
 \end{coro}

 \begin{remark} \label{dimension 2} In \cite{DHKK} is shown that   any  exceptional pair $(A,B)$  in $D^b(Q)$ for  an  acyclic Euclidean quiver $Q$ satisfies $\hom^i(A,B)\leq 2$.   Among the pairs of $Rep_k(Q_1)$ listed in Corollary \ref{coro Q1 exceptional pairs} equality is
  attained   in the following cases: $2=\hom(E_1^{m+1},E_1^m)$ $=\hom(E_2^{m},E_2^{m+1})=$ $\hom(E_3^{m},E_3^{m+1})=$
   $\hom(E_4^{m+1},E_4^m)=$ $\hom^1(E_1^{0},E_2^0)=$ $\hom^1(E_4^{0},E_3^0)$.
 \end{remark}

\begin{remark} \label{braid} From Corollaries \ref{exceptional colleections} and  \ref{coro for isom triples} we see that the action of  the Braid group $B_3$ on the exceptional collections of $Rep_k(Q_1)$ is not free.  We give an example here.\end{remark}

\noindent \textbf{Example of  fixed triples by  a Braid group element.}  \textit{For any exceptional triple $(A,B,C)$ we  denote here   the triple\footnote{ Recall that for any exceptional pair $(A,B)$  the exceptional objects $L_A(B)$ and $R_B(A)$ are determined by the triangles $\bd L_A(B)&\rTo &\Hom^*(A,B)\otimes A & \rTo^{ev^*_{A,B}} & B \ed$; \  $\bd A &\rTo^{coev^*_{A,B}} &\Hom^*(A,B)^{\check{}}\otimes B &\rTo & R_B(A) \ed$ and that $(L_A(B),A)$, $(B,R_B(A))$ are exceptional pairs.} $(A,L_B(C),B)$   by  $L_1(A,B,C)$. We keep in mind also Corollary \ref{coro for isom triples} and that each exceptional object in $D^b(Q_1)$ is a shift of an exceptional object in $Rep_k(Q_1)$.}

\textit{ The first row in the list of  Corollary \ref{exceptional colleections}  shows that, up to shifts, we have the equalities \\ $L_1(E_1^{m+1},M,E_4^m)=(E_1^{m+1},E_1^m,M)$; $L_1(E_1^{m+1},E_1^m,M)=(E_1^{m+1},E_4^m,E_1^m)$; $L_1(E_1^{m+1},E_4^m,E_1^m)=(E_1^{m+1},M,E_4^m)$. Hence, the triple $(E_1^{m+1},M,E_4^m)$ is fixed by  $(L_1)^3$. The element  $(L_1)^3$ is not trivial in the braid group $B_3$, since $B_3$ is torsion free.}

\textit{
  Acting with $L_1$ on each of the rest  rows, except   the last two rows, we find the same
 behavior.}

\section{Preliminaries}
Here we comment on Bridgeland's stability conditions  and on Macr\`i's construction of stability  conditions via exceptional collections.

In Subsection \ref{comments on stab cond} for a  a  Krull-Schmidt category ${\mathcal T}$, we introduce  a function  $\bd Ob({\mathcal T})& \rTo^{\theta_\sigma}& \NN^{(\sigma^{ss}_{ind}/\cong)} \ed $, depending on a stability condition $\sigma \in  \st(\mc T)$. It helps us later to encode  useful  features of the relation $\bd R & \rDotsto & (S,E)\ed$ in the simple   expressions $\theta_\sigma(R)>\theta_\sigma(E)$, $\theta_\sigma(R)(S)>0$(see Section \ref{some terminilogy}).
Lemma \ref{finiteness coro}, based on the locally finiteness of the elements in $\st(\mc T)$, has an important  role in Section \ref{main theo for Q_1}.
The simple fact observed in Lemma \ref{lemma for hom leq 1(X,S)}, used throughout  Sections  \ref{no bad after good},..., \ref{main theo for Q_1}, is helpful in our study of long $R$-sequences.

Applying some results  of Section \ref{two examples}, we obtain in Subsection \ref{two limit points} stability conditions on $D^b(Q_1)$ with two limit points in $\SS^1$, which concerns \cite[table (1)]{DHKK}.

After having recalled Macr\`i's construction in Subsection \ref{subsection about macri},  we define   in the final Subsection \ref{def of sigma-exceptional collection} the notion of a $\sigma$-exceptional collection.

\subsection{Krull-Schmidt property. The function \texorpdfstring{$\theta: Ob({\mathcal C}) \rightarrow \NN^{({\mathcal C}_{ind}/\cong )}$}{\space}} \label{section KSA} \mbox{}\\
Let  $\mathcal C$ be  an additive category. We denote by ${\mathcal C}_{ind}$ the set of all indecomposable objects in $\mathcal C$.\footnote{the set ${\mathcal C}_{ind}$  does not contain zero objects.}
  We discuss here the  well known Krull Schmidt property.

\begin{df} \label{def of kr-schm} We say that an additive category ${\mathcal C}$ has \uline{Krull-Schmidt property} if for each
$X \in Ob( { \mathcal C } )\setminus \{0\}$ there exists unique up to isomorphism and permutation  sequence  $ \{X_1, X_2,\dots X_n\}$ in  $ {\mathcal C}_{ind}$  with  $X \cong \bigoplus_{i=1}^n X_i  $.

For $X \in Ob( { \mathcal C } )\setminus \{0\}$ with a
decomposition   $X \cong \bigoplus_{i=1}^n X_i  $ as above we denote by $Ind(X)$ the set $\{Y \in Ob({\mc C}) :
  \  Y \cong X_i \ \mbox{for some} \ i=1,2,\dots,n\}$. If $X$ is a zero
object, then $Ind(X)=\emptyset$.
\end{df}
We will use two simple observations related to this property.
\begin{lemma} \label{from Hered to derived with KSA} Let  $\mathcal A$ be  a hereditary abelian
category. If $\mathcal A$ has  Krull-Schmidt property,
then $D^b({\mathcal A})$  has Krull-Schmidt property.
\end{lemma}
\bpr Recall that   any object $X \in D^b({\mathcal A})$ decomposes as follows
$X \cong \bigoplus_{i} H^i(X)[-i]$ and if $X \cong \bigoplus_{i}
X_i[-i]$ for some collection $\{X_i\} \subset {\mathcal A}$, then
$ X_i \cong H^i(X) $ for all $i$. In particular $\mc A$ is a thick subcategory of $D^b(\mc A)$. Now the lemma follows. \epr
\begin{lemma} \label{defintion of theta}
Let  ${\mathcal C}$  have   Krull-Schmidt property.  There exists unique function
$ \bd  Ob({\mathcal C}) & \rTo^{\theta} & \NN^{({\mathcal C}_{ind}/\cong )} \ed$  satisfying:\footnote{By $\NN^{({\mathcal C}_{ind}/\cong )}$ we denote  the set of functions  from ${\mathcal C}_{ind}/\cong$ to $\NN$ with finite support.}
\begin{itemize}
    \item[\textbf{(a)}] If $Y \cong \bigoplus_{i=1}^m Y_i$ in  ${\mathcal C}$, then $\theta(Y)=\sum_{i=1}^n \theta(Y_i)$.
    \item[\textbf{(b)}]     For  any $X \in {\mathcal C}_{ind}$  the function $ \bd {\mathcal C}_{ind}/\cong & \rTo^{\theta(X)}& \NN \ed$  assigns one to the equivalence class containing $X$, and zero elsewhere.
\end{itemize}
\end{lemma}
\bpr For an object $X \in Ob({\mathcal C}) $ with a decomposition $X \cong \bigoplus_{i=1}^n X_i$  as in Definition \ref{def of kr-schm}    the function $ \bd {\mathcal C}_{ind}/\cong & \rTo^{\theta(X)}& \NN \ed $ assigns  to each  $ u \in {\mathcal C}_{ind}/\cong $  the number $\#\{i : X_i \in u \}$. \epr

\subsection{Comments on stability conditions. The family \texorpdfstring{$\{  \theta_\sigma : Ob({\mathcal T}) \rightarrow \NN^{(\sigma^{ss}_{ind}/\cong)} \}_{\sigma \in \st(\mc T)}$}{\space}} \label{comments on stab cond} \mbox{}\\

 Recall
that if $\sigma = ({\mathcal P}, Z)$ is  a locally finite
stability condition on a triangulated category ${\mathcal T}$,
then for each $ t \in \RR $ the subcategory  ${\mathcal P}(t)$ is
an
 abelian category of finite length (see  \cite[p. 6]{Bridg2}). Furthermore \cite{Bridg1}, the short exact sequences in $\mc P(t)$ are exactly these sequences
  $\bd A & \rTo ^{\alpha}& B & \rTo^{\beta} & C \ed $ with $A,B,C \in \mc P(t)$, s. t.  for some $\gamma : C \rightarrow A[1]$ the sequence
$\bd A & \rTo ^{\alpha}& B & \rTo^{\beta} & C & \rTo^{\gamma} A[1]\ed $ is a  triangle in $\mc T$.
 The first lemma in this subsection, used in Section \ref{main theo for Q_1}, follows from  \textbf{locally finiteness}.
\begin{lemma}\label{finiteness coro}  Let $\sigma=({\mc P}, Z) \in \st(\mc T)$, $t \in \RR$, $A\in {\mc P}(t)$. For any object $X \in \mc T$ denote
by $[X] \in K(\mc T)$ the corresponding equivalence class in the Grothendieck group $K(\mc T)$. Then  the set
\be \label{finiteness coro formula} \{[X]\in K(\mc T) :   X \in {\mc P}(t) \  \mbox{and there exists a monic arrow} \  X \rightarrow A   \mbox{ in} \ {\mc P}(t) \}\ee
is finite.
\end{lemma}
\bpr Since $\mc P(t)$ is abelian category of finite length,
 we have a Jordan-Holder filtration for the given $A \in \mc P(t)$
\begin{diagram}[size=1em] \label{JH filtration}
0 & \rTo      &     &       &   E_1 & \rTo      &     &        & E_2 & \rTo & \dots & \rTo & E_{n-1} & \rTo      &     &       & E_n=A \\
  &           &     & \ldTo &       &           &     &  \ldTo &     &      &       &      &         &           &     & \ldTo &       \\
  &           & S_1 &       &       &           & S_2 &        &     &      &       &      &         &           & S_n &
\end{diagram}
where $E_i \rightarrow E_{i+1} \rightarrow S_{i+1}$ are  short exact sequences in $\mc P(t)$ and $S_1, S_2, \dots, S_n$ are simple objects in $\mc P(t)$.
We will show that the set \eqref{finiteness coro formula} is finite by showing that it  is a subset of:
$$\left \{\sum_{i=1}^m [S_{\xi(i)}] : \bd  \{1,2,\dots,m\} & \rTo^\xi & \{1,2,\dots,n\}\ed \ \mbox {is injective} \right \}.$$

   For any monic arrow $X \rightarrow A$ in $\mc P(t)$ we have a Jordan-Holder  filtration of $X$
\be \label{JH filtration5} \begin{diagram}[size=1em]
0 & \rTo      &     &       &   E_1' & \rTo      &     &        & E_2' & \rTo & \dots & \rTo & E_{n-1}' & \rTo      &     &       & E_m'=X \\
  &           &     & \ldTo &       &           &     &  \ldTo &     &      &       &      &         &           &     & \ldTo &       \\
  &           & S_1' &       &       &           & S_2' &        &     &      &       &      &         &           & S_m' &
\end{diagram} \ee
where $S_1', S_2', \dots, S_m'$ are
simple objects in $\mc P(t)$, s. t. $S_i' \cong S_{\xi(i)},
i=1,\dots,m$ for some injection $\xi: \{1,2,\dots,m\} \rightarrow
\{1,2,\dots,n\}$.   Since $E_i' \rightarrow E_{i+1}' \rightarrow
S_{i+1}'$ is a short exact sequences in $\mc P(t)$, it is also a
part of a triangle $E_i' \rightarrow E_{i+1}' \rightarrow S_{i+1}'
\rightarrow E_i'[1]$ in $\mc T$. Hence  by \eqref{JH filtration5} it follows $[X]= \sum_{i=1}^m [S_i']=\sum_{i=1}^m [S_{\xi(i)}]$.
\epr
Recall  that one of Bridgeland's axioms   \cite{Bridg1}  is: for any nonzero $X \in Ob({\mathcal T})$ there exists  a
diagram of triangles,\footnote{Throughout  the whole text the word triangle means distinguished triangle.} called \textbf{Harder-
Narasimhan filtration}:
\be
\begin{diagram}[size=1em] \label{HN filtration}
0 & \rTo      &     &       &   E_1 & \rTo      &     &        & E_2 & \rTo & \dots & \rTo & E_{n-1} & \rTo      &     &       & E_n=X \\
  & \luDashto &     & \ldTo &       &\luDashto  &     &  \ldTo &     &      &       &      &         & \luDashto &     & \ldTo &       \\
  &           & A_1 &       &       &           & A_2 &        &     &      &       &      &         &           & A_n &
\end{diagram}
\ee where $\{ A_i  \in {\mathcal P}(t_i) \}_{i=1}^n $, $t_1 > t_2
> \dots > t_n $ and $A_i$ is non-zero object for any $i=1,\dots,n$
(the non-vanishing  condition makes the factors $\{ A_i  \in
{\mathcal P}(t_i) \}_{i=1}^n$ unique up to isomorphism). In
\cite{Bridg1} is used the notation $\phi^\sigma_-(X):=t_n$,
$\phi^\sigma_+(X):=t_1$, and the phase of a semistable object $A
\in {\mathcal P}(t)\setminus \{0\}$
 is denoted by $\phi^{\sigma}(A):=t$, we also use these notations. \textit{The objects $\{ A_i \}_{i=1}^n$ will be called   HN factors of $X$} (HN for Harder-
Narasimhan). It is useful to give a name of the minimal  HN factor $A_n$.
\begin{df}\label{def of sigma_-} For any  $X\in \mc T\setminus \{0\}$ we choose\footnote{by the axiom of choice} a Harder-Narasimhan filtration as in   \eqref{HN filtration}. Having this diagram, we denote the semistable HN factor of minimal phase $A_n$  by $\sigma_-(X)$, and    the last triangle  $\bd E_{n-1} & \rTo      & X & \rTo & A_n & \rTo & E_{n-1}[1] \ed$ by  ${\rm H\! N}_-(X)$. In particular, $\phi(\sigma_-(X))=\phi_-(X)$.
 \end{df}
In the next Lemma \ref{lemma for hom leq 1(X,S)} we  treat $\sigma_-(X)$. We recall first (another axiom of Bridgeland \cite{Bridg1}) that  \textbf{from $\phi(A)>\phi(B)$ with semistable $A$, $B$ it follows $\hom(A,B)=0$}. This axiom implies that from   $\phi_-(X)>\phi_+(Y)$ it follows  $\hom^{\leq 0}(X,Y)=\hom^{\leq 0}(\sigma_-(X),Y)= 0$. We get $\hom^{\leq 1}(\sigma_-(X),Y)= 0$ in the following situation:
\begin{lemma} \label{lemma for hom leq 1(X,S)} If $\phi_-(X) \geq \phi_+(Y)$ and
  $\hom^{\leq 1}(X,Y)=0$, then   $\hom^{\leq 1}(\sigma_-(X), Y)=0$. \end{lemma}
\bpr Let $ {\rm HN}_-(X)= \bd[1em]
Z & \rTo      &  X  & \rTo & \sigma_-(X) & \rTo & Z[1]
\ed$. Then $\phi_-(Z) >\phi(\sigma_-(X))=\phi_-(X) \geq \phi_+(Y)$.
Hence $\Hom^{\leq 0}(Z,Y)=0$. We apply $\Hom(\_,Y[i])$ with
$i\leq 1$ to this triangle  and obtain:\\ $ 0=\Hom(Z[1], Y[i])
\rightarrow \Hom(\sigma_-(X), Y[i]) \rightarrow \Hom(X, Y[i])=0.$ The lemma follows. \epr

 In \cite{Bridg1} for a slicing $\mc P$ of $\mc T$ and an interval $I\subset \RR$ by $\mc P(I)$ is denoted the extension closure of $\{\mc P(t)\}_{t\in I}$, and     ${\mathcal P}([t, t+1)), {\mathcal P}((t, t+1])$ are shown to be hearts of bounded t-structures for any $t \in \RR$.
If $\mc P$ is a part of a stability condition $({\mathcal P}, Z) \in \st(\mathcal T)$, then    ${\mathcal P}(t)$ is shown to be abelian. 
The nonzero objects in the subcategory $\mc P(I)$ are exactly those $X\in \mc T\setminus \{0\}$, which  satisfy $\phi_\pm(X)\in I$.

From these facts it follows  that $\mc P(I)$ is a thick subcategory for any interval $I\subset \RR$:

\begin{lemma}\label{P(t) is a thick subcategory} For any slicing $\mathcal P$ of a triangulated category ${\mathcal T}$ and any interval $I \subset \RR$ the category
${\mathcal P}(I)$ is a thick subcategory of ${\mathcal T}$. In particular, if $\mc T$ has  Krull-Schmidt property, then ${\mathcal P}(I)$ has it. \end{lemma}
\bpr   In \cite{GM} t-structures are defined   as pairs of subcategories.    For any slicing $\mathcal P$ and any $t\in \RR$ the hearts  $ {\mathcal P}((t, t+1]), {\mathcal P}([t, t+1))$ come from  the pairs
$({\mathcal P}((t,+\infty)), {\mathcal P}((-\infty,t+1]))$,\\ $({\mathcal P}([t,+\infty)), {\mathcal P}((-\infty,t+1)))$, respectively, which are bounded t-structures. Let us consider for example the t-structure
$({\mathcal P}((t,+\infty)), {\mathcal P}((-\infty,t+1]))$. In terms of the notations used in \cite{GM} we denote
${\mathcal T}^{\leq 0}={\mathcal P}((t,+\infty))$, ${\mathcal T}^{\geq 0}={\mathcal P}((-\infty,t+1])$. From the properties of t-structures
we know that \ben X \in {\mathcal T}^{\leq 0} \iff \forall Y \in {\mathcal T}^{\geq 1} \ \  \hom(X,Y)=0; \qquad
 X \in {\mathcal T}^{\geq 0} \iff \forall Y \in {\mathcal T}^{\leq -1}\ \  \hom(Y,X)=0.\een
Hence ${\mathcal T}^{\leq 0}={\mathcal P}((t,+\infty))$, ${\mathcal T}^{\geq 0}={\mathcal P}((-\infty,t+1])$ are thick subcategories. Similarly
${\mathcal P}([t,+\infty))$, ${\mathcal P}((-\infty,t+1))$ are thick. Since for any interval
 $I\subset \RR$ the subcategory ${\mathcal P}(I)$ is an intersection of two subcategories of the considered types, the lemma follows.
\epr

\begin{coro} \label{ineq for phi_-} Let $X,A, B \in {\mathcal T}$ and $X \cong A \oplus B$, then for any slicing ${\mathcal P}$ of  ${\mathcal T}$ we have
$\phi_-(X) \leq \phi_-(A) \leq \phi_+(A) \leq \phi_+(X) $.
\end{coro}
\bpr We have $X \in {\mathcal P}([\phi_-(X), \phi_+(X)])$. From the previous lemma $A,B \in {\mathcal P}([\phi_-(X), \phi_+(X)])$ and the statement follows.
\epr

Thus, if $\mc T$ has Krull-Schmidt property, then all  $\{\mc P(t)\}_{t \in \RR}$ have it(Lemma \ref{P(t) is a thick subcategory}). From Lemma \ref{defintion of theta} we obtain a family of functions $ \{ {\mathcal P}(t)
\rightarrow \NN^{({\mathcal P}(t)_{ind}/\cong)} \}_{t \in \RR} $.   In Definition \ref{theta_sogma} below we build a single function on $Ob(\mc T)$ from this family of functions,  using the HN filtrations.
We need first some notations.

For $\sigma = ({\mathcal P}, Z) \in \st(\mathcal T)$ we  denote by
$\sigma^{ss}$ the  set of  $\sigma$-semistable
objects, i. e. \be
\label{sigma^{ss}} \sigma^{ss}=\cup_{t \in \RR} {\mathcal P}(t)\setminus \{0\}.
\ee  By $\sigma^{ss}_{ind}$ we denote the set of all indecomposable
semistable objects, i. e.\footnote{Recall that  $\mc P(t)$ is thick in $\mc T$ (Lemma \ref{P(t) is a thick subcategory}), hence $\mc P(t)_{ind}=\mc P(t) \cap \mc T_{ind}$.} \be
\label{sigma^ss_ind} \sigma^{ss}_{ind}=\cup_{t \in \RR} {\mathcal
P}(t)_{ind}=\sigma^{ss}\cap {\mc T}_{ind}. \ee

 In \textbf{(a)} of Definition \ref{theta_sogma} we  consider $\NN^{({\mathcal P}(t)_{ind}/\cong)}$ as a subset of  $\NN^{(\sigma^{ss}_{ind}/\cong)}$, which is reasonable since the family  $\{{\mathcal P}(t)_{ind} \}_{t\in \RR}$ is pairwise disjoint.

\begin{df}\label{theta_sogma} Let ${\mathcal T}$ have
Krull-Schmidt property. Let $\sigma = ({\mathcal P}, Z) \in
\st(\mathcal T)$.

We define   $
\theta_\sigma : Ob({\mathcal T}) \rightarrow
\NN^{(\sigma^{ss}_{ind}/\cong)} $  as the unique function satisfying the following:

 \textbf{(a)}  For each $t\in \RR$  the restriction of $\theta_\sigma$ to  ${\mathcal P}(t)$ coincides with the    function $ {\mathcal P}(t)
\rightarrow \NN^{({\mathcal P}(t)_{ind}/\cong)} $, given by Lemmas \ref{defintion of theta},  \ref{P(t) is a thick subcategory}.

 \textbf{(b)} For   any non-zero
$X \in Ob({\mathcal T})$ with a HN filtration\footnote{Recall that  the collection $\{ A_i \}_{i=1}^n $
of the HN factors is determined by $X$ up to isomorphism.}
\eqref{HN filtration} holds the equality $  \theta_\sigma (X)= \sum_{i=1}^n
\theta_\sigma(A_i). $
\end{df}
We use freely that $X\cong Y$ implies $\theta_\sigma(X)=\theta_\sigma(Y)$, $X\neq 0$ implies $\theta_\sigma(X)\neq 0$, and $\theta_\sigma(X)\leq \theta_\sigma(Y)$ implies $\phi_-(Y)\leq \phi_-(X)\leq \phi_+(X)\leq \phi_+(Y)$.  Another property of $\theta_\sigma$, to which we refer later, is:
\begin{lemma} \label{additive property of theta_sigma 1} Let $\phi_-(X_1) > \phi_+(X_2)$.  For any  triangle   $ X_1 \rightarrow X \rightarrow X_2 \rightarrow X_1[1] $
we have $\theta_{\sigma}(X )=\theta_{\sigma}(X_1)+\theta_{\sigma}(X_2)$.
\end{lemma}
\bpr If  the HN factors of $X_1$ and $X_2$  are $A_1, A_2, \dots,
A_n$ and  $B_1,B_2, \dots, B_m$, respectively,  then, using the
octahedral axiom, one can show that the HN factors of $X$ are
$A_1, A_2, \dots, A_n,$  $B_1,B_2, \dots, B_m$. Now the lemma
follows from \textbf{(b)} in Definition \ref{theta_sogma}. \epr
 The property  $\theta_{\sigma}(X \oplus Y)=\theta_{\sigma}(X)+\theta_{\sigma}(Y)$ for $X, Y \in \mc P(t)$ follows from    \textbf{(a)} in Lemma \ref{defintion of theta}. To show this additive property for    any two objects $X$, $Y \in \mc T$  we note first:
\begin{lemma} \label{lemma for without constraint for non zero} For any diagram of the type (composed of distinguished triangles):
\begin{diagram}[size=1em]
0 & \rTo      &     &       &   B_1 & \rTo      &     &        & B_2 & \rTo & \dots & \rTo & B_{n-1} & \rTo      &     &       & B_n=X,\\
  & \luDashto &     & \ldTo &       &\luDashto  &     &  \ldTo &     &      &       &      &         & \luDashto &     & \ldTo &       \\
  &           & A_1 &       &       &           & A_2 &        &     &      &       &      &         &           & A_n &
\end{diagram}
where $\{ A_i  \in {\mathcal P}(t_i) \}_{i=1}^n $, $t_1 > t_2 > \dots > t_n $, without the constraint that $A_1, A_2, \dots, A_n $ are non-zero objects,
we have $\theta_\sigma (X)= \sum_{i=1}^n \theta_\sigma(A_i)$.
\end{lemma}
\bpr We can remove all triangles where  $A_i$ is zero and in the
end we obtain the HN filtration of $X$, then the equality follows
from \textbf{(b)} in Definition \ref{theta_sogma}  and
$\theta_\sigma(A_i)=0$ if $A_i$ is a zero object. \epr

Given two non-zero objects $X_1, X_2 \in
Ob({\mathcal T})$, then  after
inserting triangles of the form $\begin{diagram}[size=1em]
E & \rTo^{Id}      &     &       &   E \\
  & \luDashto &     & \ldTo &       \\
  &           & 0 &       &
\end{diagram}$ to their HN filtrations we can obtain two($i=1,2$) equally long diagrams
with distinguished triangles
\begin{diagram}[size=1em]
0 & \rTo      &      &       &   B_1^i & \rTo      &     &        & B_2^i & \rTo & \dots & \rTo & B_{n-1}^i & \rTo      &     &       & B_n^i=X_i,  \\
  & \luDashto &      & \ldTo &       &\luDashto  &     &  \ldTo &     &      &       &      &         & \luDashto &     & \ldTo &       \\
  &           & A_1^i &       &       &           & A_2^i &        &     &      &       &      &         &           & A_n^i &
\end{diagram}
where $\{ A_j^i \in {\mathcal P}(t_j) \}_{j=1}^n$, $i=1,2$ and $t_1>t_2 > \dots > t_n$. Hence,
we get a diagram of triangles:
\begin{diagram}[size=1em]
0 & \rTo      &      &       &   B_1^1 \oplus B_1^2& \rTo      &     &        & B_2^1 \oplus B_2^2& \rTo & \dots & B_{n-1}^1 \oplus B_{n-1}^2 & \rTo      &     &       & X_1\oplus X_2.  \\
  & \luDashto &      & \ldTo &       &\luDashto  &     &  \ldTo &     &      &       &               & \luDashto &     & \ldTo &       \\
  &           & A_1^1 \oplus A_1^2&       &       &           & A_2^1\oplus A_2^2&        &     &      &             &         &           & A_n^1 \oplus A_n^2&
\end{diagram}
 We have $\{ A_j^1 \oplus A_j^2 \in {\mathcal P}(t_j)\}_{j=1}^n$ by the additivity of $ {\mathcal P}(t_j)$. Using Lemma \ref{lemma for without constraint for non zero} we obtain:
$ \theta_{\sigma}(X_1\oplus X_2)= \sum_{j=1}^n \theta_\sigma(A_j^1\oplus A_j^2)= \sum_{j=1}^n \theta_\sigma(A_j^1) + \sum_{j=1}^n \theta_\sigma(A_j^2)=\theta_{\sigma}(X_1)+\theta_{\sigma}(X_2), $ i. e. we proved:
\begin{lemma} \label{additive property of theta_sigma 2} For any  pair of objects $X_1,X_2$   in $\mathcal T$ we have:
 $ \theta_{\sigma}(X_1\oplus X_2)=\theta_{\sigma}(X_1)+\theta_{\sigma}(X_2).$
\end{lemma}
In the end of this subsection we recall the remaining    axioms of Bridgeland \cite{Bridg1}.   \textbf{A stability condition $\sigma=(\mc P, Z)\in \st(\mc T)$ has the properties:  ${\mathcal P}(t)[1]={\mathcal P}(t+1)$ for each $t\in \RR$, and} \be
\label{phase formula}  X\in \sigma^{ss} \qquad
\Rightarrow \qquad  Z(X)=r(X)\ \exp(\ri \pi \phi(X)), \ r(X) >0.
\ee
\subsection{Application: Stability conditions on \texorpdfstring{$D^b(Q_1)$}{\space} with two limit points in \texorpdfstring{$\SS^1$}{\space}}  \label{two limit points} \mbox{}\\
We start by recalling  a result in \cite{Bridg1}:
\begin{prop}[Proposition 5.3 in \cite{Bridg1}] \label{from abelian to triangulated} Let $\mc A \subset \mc T$ be a bounded t-structure in a triangulated category $\mc T$ and $\bd K(\mc A)& \rTo^{Z}& \CC \ed $ be a
stability function on $\mc A$  with HN property.\footnote{HN property for $\bd K(\mc A)& \rTo^{Z} & \CC \ed $ is defined in \cite[Definition 2.3]{Bridg1}. If $\mc A$ is an abelian category of finite length, then any  stability function  $Z$ on $\mc A$ satisfies the HN property \cite[Proposition 2.4]{Bridg1}. } Then there exists unique stability condition\footnote{If $\mc A$ has  finite length  and   finitely many simple objects, then the obtained  stability condition $\sigma$ is locally finite.} $\sigma=(\mc P, Z_e)$ on $\mc T$ satisfying:
\begin{itemize}
    \item[\textbf{(a)}]$Z_e(X)=Z(X)$  for $X \in \mc A$;
    \item[\textbf{(b)}] For $t\in (0,1]$ the objects of ${\mc P}(t)$ are:\footnote{For $u\in \HH$ we denote by $\arg(u)$ the number satisfying $\arg(u)\in (0,1]$, $u=\exp(i \arg(u))$. We set $\arg(0)=-\infty $.} \\ $Ob({\mc P}(t))=\{ X\in {\mc A} :   \ \mbox{ for each} \  \mc A\mbox{-monic } \     X' \rightarrow X \ \  \arg Z(X')\leq  \arg Z(X)= \pi t \} $  .
\end{itemize}
Conversely, for each stability condition $\sigma=(\mc P, Z_e)$ on
$\mc T$ the subcategory $\mc P((0,1]) = \mc A$ is a heart of a
bounded t-structure of $\mc T$, the restriction  $Z =Z_e\circ
(K(\mc A) \rightarrow  K(\mc T)) $ of $Z_e$ to $K(\mc A)$ is a
stability function on $\mc A$ with HN property and for $t\in
(0,1]$ the set of objects of $\mc P(t)$ is the same as in
\textbf{(b)}.
\end{prop}
\begin{df} \label{H^A} We   denote by $\HH^{\mc A}$
the   family   of stability conditions on $\mc T$ obtained by
\textbf{(a)}, \textbf{(b)} above  keeping $\mc A$ fixed and
varying $Z$ in the set of all stability functions on $\mc A$
with  HN property.
\end{df}

 Let ${\mc A} = Rep_k(Q_1) \subset D^b(Rep_k(Q_1))$ be the standard bounded $t$-structure, where  $k$ is an algebraically closed field.
  A stability function $\bd K(\mc A)& \rTo^{Z}& \CC \ed$ is uniquely determined by  $Z(E_1^0)$, $Z(M),$ $ Z(E_3^0) \in \HH$ .   Here we  choose $Z(E_1^0), Z(M), Z(E_3^0)$ as follows:

\vspace{-30mm}
\setlength{\unitlength}{6mm}
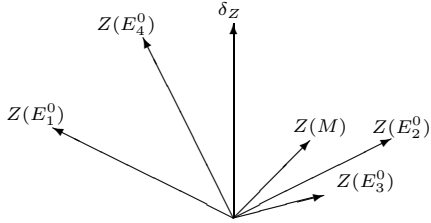
\begin{picture}(10,10)(-5,0)
\put(-4.3,2.3){\makebox(0,0){\tiny $ Z(E_1^0)  $ \normalsize}}
\put(0,0){\vector(-2,1){4}}
\put(2,2){\makebox(0,0){\tiny $Z(M)$ \normalsize}}
\put(0,0){\vector(1,1){1.7}}
\put(3,0.8){\makebox(0,0){\tiny $Z(E_3^0)$ \normalsize}}
\put(0,0){\vector(4,1){2}}
\put(3.8,2){\makebox(0,0){\tiny $Z(E_2^0)$ \normalsize}}
\put(0,0){\vector(2,1){3.5}}
\put(-2.3,4.3){\makebox(0,0){\tiny $Z(E_4^0)$ \normalsize}}
\put(0,0){\vector(-1,2){2}}
\put(0,4.6){\makebox(0,0){\tiny $\delta_Z$ \normalsize}}
\put(0,0){\vector(0,1){4.3}}
\end{picture}
where $ \begin{array}{c} \delta_Z=Z(M)+Z(E_1^0)+ Z(E_3^0),  \\ Z(E_4^0)=Z(E_1^0)+Z(M), Z(E_2^0)=Z(E_3^0)+Z(M). \end{array}$

In this subsection we show that:
     \begin{lemma} \label{P_sigma} Let  $\sigma = (\mc P, Z) \in \HH^{\mc A}\subset \st(D^b(\mc A))$ be the stability condition, uniquely  determined by  the chosen stability function $\bd[1em] K(\mc A)& \rTo^{Z}& \CC \ed $ and Bridgeland's Proposition \ref{from abelian to triangulated}. Then the set of stable phases   $P_\sigma$,  defined by\footnote{In \cite{DHKK} is shown that  if  $Q$ is  an Euclidean quiver without oriented cycles,
 then for any $\sigma \in D^b(Q)$ the set of phases $P_\sigma$ is either finite or has two limit points.}
$  P_\sigma= \exp\left ( \ri \pi \{\phi_{\sigma}(X): X \in \sigma^{ss}, X \neq 0 \} \right ), $ has  two limit points in $\SS^1$. \end{lemma}
\bpr The  values of $Z$ on $\mc A_{exc}$(see Proposition \ref{exceptional objects in Q1}) are: $Z(M)$, $Z(M')=Z(E_1^0)+Z(E_3^0)$, and
\begin{gather}  \label{Z(E_i^m)}  Z(E_j^m) = m \delta_Z +  Z(E_j^0), \ \ m \in \NN, j=1,2,3,4. \end{gather}
By Kac's theorem,\footnote{saying that the dimension vectors of the indecomposables are the same as the roots} proved in \cite{Kac}, and the description of the roots before Lemma \ref{roots no exceptional} it follows that the values of $Z$ on all indecomposable objects are the already given above and   the following:
\begin{gather}\label{non except}\{ m\delta_Z \}_{m\geq 1}, \ \ \{ m\delta_Z + Z(M), \ \ m\delta_Z + Z(E_1^0)+Z(E_3^0) \}_{m\in \NN }. \end{gather}
We will show below that,  due to the choice of   $Z(E_1^0),
Z(M),Z(E_3^0)$,  for each $m\in \NN$ the object $E_4^m\in
Rep_k(Q_1)$ satisfies the conditions  in Corollary \ref{from
abelian to triangulated} \textbf{(b)}, so  $E_4^m \in {\mc
P}(\frac{1}{\pi} \arg(Z(E_4^m)))$. Hence   for any $m \in \NN$, $i
\in \ZZ$ we obtain\footnote{Recall that ${\mathcal
P}(t)[1]={\mathcal P}(t+1)$ for each $t\in \RR$.} $\mc
P(\frac{1}{\pi} \arg(Z(E_4^m))+i) \neq 0$. By \eqref{phase
formula} we can write $\{\pm Z(E_4^m)/\abs{Z(E_4^m)} \}_{m \in
\NN} \subset P_\sigma$. From \eqref{Z(E_i^m)} we have  $
\lim_{m\rightarrow \infty }\frac{ Z(E_4^m)}{\abs{Z(E_4^m)}} =
\lim_{m\rightarrow \infty } \frac{ m \delta_Z +Z(E_4^0)}{\abs{m
\delta_Z +Z(E_4^0)}}= \pm \delta_Z/\abs{ \delta_Z }$. Since
$\delta_Z$ is not collinear with $Z(E_4^0)$, it follows that $\pm
\delta_Z/\abs{ \delta_Z }$ are limit points of $P_\sigma$.

It remains to show that $E_4^m$ satisfies the conditions  in
Corollary \ref{from abelian to triangulated} \textbf{(b)} for each
$m\in \NN$.
 Since  $\mc A$ has Krull-Schimidt property, it is enough to show that   any monic $X \rightarrow E_4^m$ with $X\in {\mc A}_{ind}$ satisfies
$\arg(Z(X))\leq \arg(Z(E_4^m))$.  In \eqref{Z(E_i^m)}, \eqref{non except} are given all  the values $\{Z(X) : X \in {\mc A}_{ind}\}$. From the picture we see that if
 $u\in \{Z(E_2^j), Z(E_3^j), j \delta_Z, j \delta_Z+Z(M)\}_{j \in \NN} \cup \{Z(E_4^j)\}_{j \geq m}$, then $\arg(u)\leq \arg(Z(E_4^m))$. Hence,  it remains to show that
  any monic $X \rightarrow E_4^m$ with $X\in {\mc A}_{ind}$ and $\ul{\dim}(X) \in \{(j+1,j,j), (j+1,j,j+1)\}_{j \in \NN} \cup \{(j+1,j+1,j)\}_{j< m}$ satisfies
  $\arg(X)\leq \arg(Z(E_4^m))$. We consider separately two options.

  If either $\ul{\dim}(X) =(j+1,j,j)$ or $\ul{\dim}(X) =(j+1,j,j+1)$, then a morphism $X\rightarrow E_4^m$ consists of three vector space
  morphisms $(f_b,f_{mid},f_e):(k^{j+1}, k^j, k^{x}) \rightarrow (k^{m+1},k^{m+1},k^m)$, where $x\in \{j,j+1\}$, satisfying three relations. One of these relations is of the form \\  $( \bd k^{j+1}& \rTo & k^j & \rTo^{f_{mid}} & k^{m+1}  \ed )=( \bd k^{j+1}& \rTo^{f_b} & k^{m+1} & \rTo^{Id} & k^{m+1} \ed) =f_b$,  therefore $f_b$ cannot be injective vector space morphism. This implies that
  $X \rightarrow E_4^m$ cannot be monic in $\mc A$.

 If  $\ul{\dim}(X)= (j+1,j+1,j)$ with $j<m$, then, as far as $(j+1,j+1,j)$ is a real root, by Kac's theorem there exists unique up to isomorphism such $X \in {\mc A}_{ind}$. From Proposition \ref{exceptional objects in Q1}
 it follows $X\cong E_4^j$. However from table \eqref{Q1 table} we see $\hom(E_4^j,E_4^m)=0$ for $j<m$, therefore again we have not any monics $X \rightarrow E_4^m$.  The lemma is proved.\epr
\subsection{On the stability conditions constructed by E. Macr\`i via exceptional collections} \label{subsection about macri} \mbox{} \\

\textbf{E. Macr\`i proved in \cite[Lemma 3.14]{Macri} that the extension closure  ${\mc A}_{\mc
E}$   of a full Ext-exceptional collection ${\mc
E}=(E_0,E_1,\dots,E_n)$ in $\mc T$ is
  a heart of a bounded t-structure}. Furthermore,
 ${\mc A}_{\mc E}$ is of finite length and $E_0, E_1,\dots, E_n$ are the simple objects in it.  
  Bridgeland's Proposition  \ref{from abelian to triangulated} produces
  a family $\HH^{\mc A_{\mc E}}\subset \st(\mc T)$ (see Definition \ref{H^A}).

\begin{df} \label{H mc E} Let ${\mc E}$ be a full Ext-exceptional collection and let ${\mc A}_{\mc E}$ be its extension closure.
     We write $\HH^{\mc E}$ for  $\HH^{\mc A_{\mc E}}$ and  denote  by $\Theta_{\mc E}'\subset \st(\mc T)$  the set obtained by acting on
      $\HH^{\mc E}$  with $\widetilde{GL}^+(2,\RR)$.
 \end{df}

 If $\mc T$ is of finite type, then starting with any full exceptional collection ${\mc E}=(E_0,E_1,\dots,E_n)$
 the collection ${\mc E}[p]=(E_0[p_0],E_1[p_1],\dots,E_n[p_n])$ is Ext for some
 integer vector $p=(p_0,p_1,\dots, p_n)\in \ZZ^{n+1}$ and  to each such vector
  corresponds a  subset $\Theta_{{\mc E}[p]}' \subset \st(\mc T)$. E. Macr\`i denotes  the union of these open subsets   by $\Theta_{{\mc E}}$,
  and  the union of the subsets $\{ \Theta_{{\mc M}}: {\mc M} \ \mbox{is a mutation of} \ {\mc E} \}$ by $\Sigma_{\mc E}$, i. e. \begin{gather}\label{theta_{mc E}} \Theta_{{\mc E}}= \bigcup_{\{p \in \ZZ^{n+1} : {\mc E}[p] \ \mbox{is Ext} \}} \Theta_{{\mc E}[p]}'\subset \st(\mc T); \qquad \qquad \Sigma_{\mc E} = \bigcup_{ \{ \Theta_{{\mc M}}: {\mc M} \ \mbox{is a mutation of} \ {\mc E} \}}\Theta_{{\mc M}}. \end{gather}
 Lemma 3.19 in \cite{Macri}  says that $\Theta_{{\mc E}}$ is an open, connected and simply connected subset of $\st (\mc T)$, which implies  (see  \cite[ Corollary 3.20]{Macri}) that, if  all  iterated mutations of ${\mc E}$  are regular,\footnote{Here  \textit{regular}  means that for $0\leq i \leq n-1$ at most one degree in   $\{ \Hom^{p}(E_i, E_{i+1})=0\}_{p \in \ZZ}$ does not vanish.}
   then   $\Sigma_{\mc E}$ is an open, connected subset of $\st(\mc T)$.

The following proposition  ensures   extendability of certain stability conditions used in \cite[Section 3]{DHKK}. The statement of Proposition \ref{projection}  is  a slight modification of  the first part of  \cite[Proposition 3.17]{Macri}.   The  difference is  that in the statement of  \cite[Proposition 3.17]{Macri} is claimed that one must take ${\mc E}_{ij} = (E_i,E_j)$, whereas  we take ${\mc E}_{ij} = (E_i, E_{i+1}, \dots ,E_j)$.  For the sake of clarity, we give a proof of Proposition \ref{projection} here.

\begin{prop} \label{projection} Let $\mc E = (E_0, E_1, \dots ,E_n)$ be a full Ext-exceptional collection in $\mc T$. Let $0\leq i< j \leq n$ and denote   ${\mc E}_{ij} = (E_i, E_{i+1}, \dots ,E_j)$, ${\mc T}_{ij}=\left \langle {\mc E}_{ij}\right \rangle \subset \mc T$. Let $\HH^{{\mc E}_{ij}} \subset \st({\mc T}_{ij})$, $\HH^{{\mc E}} \subset \st({\mc T})$  be the corresponding families  as in Definition \ref{H mc E}.

Then the map $\pi_{i j}: \HH^{{\mc
E}} \rightarrow \HH^{{\mc E}_{ij}}$, which assigns to $(\mc P, Z)
\in \HH^{{\mc E}}$ the unique $({\mc P}', Z') \in \HH^{{\mc
E}_{ij}}$ with $\{ Z'(E_k)=Z(E_k) \}_{k=i}^j$, is surjective. For any  $(\mc P, Z) \in \HH^{{\mc E}}$ and  $({\mc P}', Z') \in \HH^{{\mc E}_{ij}}$ holds the implication  \be \pi_{ij}(\mc P, Z)= ({\mc P}', Z') \ \ \Rightarrow \
\ \{{\mc P}'(t)=\mc P(t)\cap \mc       T_{ij}\}_{t\in\RR}. \ee
\end{prop}
\bpr
Using the definition of $\HH^{{\mc E}}$, $\HH^{{\mc E}_{ij}}$ (Definitions  \ref{H^A}, \ref{H mc E}), one easily reduces the proof of this proposition to the following   lemma   (compare with   the proof of   \cite[Proposition 3.17, p.7]{Macri}). \epr
\begin{lemma}  Let $\mc E $, $\mc E_{ij}$ be as in Proposition \ref{projection}. Let us denote  by $\mc A$, $\mc A_{i j}$ the extension closures of $\mc E $ and $\mc E_{ij}$ in $\mc T$. Then $\mc A_{i j}$ is an exact Serre subcategory of $\mc A$. In particular the embedding functor induces an embedding $K(\mc A_{i j}) \rightarrow K(\mc A)$.
\end{lemma}
\bpr    Since both $\mc A, \mc A_{ij}$ are abelian categories (\cite[Lemma 3.14]{Macri} ), if  $\mc A_{ij}$ is a Serre
subcategory of $\mc A$ it follows that  $\mc A_{ij}$ is an exact subcategory. Whence, it is enough to show that $\mc A_{ij}$ is a
Serre subcategory. Let $0\rightarrow B_1 \rightarrow S \rightarrow B_2 \rightarrow 0$ be any short exact sequence  in $\mc A$.

Assume that $B_1,B_2 \in \mc A_{ij}$.  Since $\mc A$ is a heart of  bounded t-structure\footnote{Recall that    the short exact sequences in a  heart of a t-structure $\mc A$  are exactly those sequences
  $\bd[1em] A & \rTo ^{\alpha}& B & \rTo^{\beta} & C \ed $ with $A,B,C \in  \mc A$, s. t. for some $\gamma : C \rightarrow A[1]$ the triangle
$\bd[1em] A & \rTo ^{\alpha}& B & \rTo^{\beta} & C & \rTo^{\gamma} A[1]\ed $ is  distinguished in $\mc T$.} in $\mc T$, the given short exact sequence is part of a  triangle in $\mc T$. Since  $\mc A_{ij}$ is extension closed in $\mc T$, it follows   $S\in \mc A_{ij}$.

Next, assume  that  $S \in \mc A_{i j}$. We have to show that   $B_1,B_2 \in \mc A_{i j}$. By $B_1,B_2  \in \mc A$ and   the definition of $\mc A$,  we have   diagrams of short exact sequences in $\mc A$ for $l=1,2$(the superscript is a power of $E_i$):
\be\label{diagrams in proof of Serre} \begin{diagram}[size=1em]
0 & \rTo      &     &       &   U_{l,n} & \rTo      &     &        & U_{l,n-1} & \rTo & \dots & \rTo & U_{l,1} & \rTo      &     &       & U_{l,0}=B_l \\
  &           &     & \ldTo &       &           &     &  \ldTo &     &      &       &      &         &           &     & \ldTo &       \\
  &           & E_n^{p_{l, n}} &       &       &           & E_{n-1}^{p_{l, n-1}} &        &     &      &       &      &         &           & E_0^{p_{l, 0}} &
\end{diagram} \qquad   l=1,2. \ee
From $S \in \mc A_{ij}$  it follows $\Hom^*(S,E_l)=0$ for $l<i$ and  $\Hom^*(E_l,S)=0$ for $l>j$. Since we have  $\mc A$-epic arrows $S\rightarrow B_2$, $B_2\rightarrow E_0^{p_{2,0}}$ and  $\mc A$-monic arrows $E_0^{p_{1, n}} \rightarrow B_1$, $B_1\rightarrow S$,  it follows that   $p_{2,0}=0$, if $0<i$ and  $p_{1,n}=0$, if $n>j$. Now by induction  it follows:
\be \label{serre first ineqs} p_{2,k}=0 \  \mbox{for}  \  k<i ,  \qquad  p_{1,k}=0 \  \mbox{for}  \ k>j. \ee

 We show bellow that $\Hom(E_k,B_2)=0$ for $k>j$   and  $\Hom(B_1,E_k)=0$ for $k<i$. Since there exist  $\mc A$-monic $E_n^{p_{2, n}} \rightarrow B_2$ and  $\mc A$-epic $B_1\rightarrow E_0^{p_{1,0}}$,  by  the diagrams  \eqref{diagrams in proof of Serre}  and  induction  we obtain   $  p_{2,k}=0 \  \mbox{for}  \  k>j ,  \quad  p_{1,k}=0 \  \mbox{for}  \ k<i$.   These vanishings together  with \eqref{serre first ineqs} imply the lemma.

 Having  \eqref{diagrams in proof of Serre} and \eqref{serre first ineqs} we can write  $B_2 \in \left \langle E_i,E_{i+1},\dots,E_n \right \rangle$ and  $B_1 \in \left \langle E_0,E_1,\dots,E_j \right \rangle$,
  hence \begin{gather}\label{serre step} \Hom^*(B_2,E_k)=\Hom^*(S,E_k)=0 \ \mbox{for} \  k<i, \Hom^*(E_k,B_1)=\Hom^*(E_k,S)=0 \ \mbox{for} \  k>j. \end{gather} From the short exact sequence $0\rightarrow B_1 \rightarrow S \rightarrow B_2 \rightarrow 0$ in $\mc A$  we get a distinguished triangle  $ B_1 \rightarrow S \rightarrow B_2 \rightarrow B_1[1]$ in $\mc T$.  Since we have \eqref{serre step}, applying to this triangle $\Hom(E_k,\_)$ and $\Hom(\_,E_k)$   we obtain the desired $\Hom(E_k,B_2)=0$ for $k>j$, $\Hom(B_1,E_k)=0$ for $k<i$. \epr

\subsection{\texorpdfstring{$\sigma$}{\space}-exceptional collections} \label{def of sigma-exceptional collection}
Motivated by the work of E. Macr\`i,  discussed in the introduction and in the  previous Subsection \ref{subsection about macri}, we define:
\begin{df}\label{sigma exceptional collection} Let $\sigma =({\mc P}, Z) \in \st(\mc T)$. We call an exceptional collection ${\mc E }=(E_0,E_1,\dots,E_n)$ \\
\underline{$\sigma$-exceptional collection} if the following properties hold:
\begin{itemize}
    \item $\mc E$ is semistable w. r. to $\sigma$ (i. e. all $E_i$ are semistable).
    \item $\forall i\neq j$ $\hom^{\leq 0}(E_i, E_j)=0$ (i. e. this is an Ext-exceptional collection).
    \item There exists $t \in \RR$, s. t. $\{\phi(E_i)\}_{i=0}^n \subset (t,t+1]$.
\end{itemize}
\end{df}

The set stability conditions for which $\mc E$ is $\sigma$-exceptional coincides with   $\Theta_{\mc E}'=\HH^{\mc E}\cdot\widetilde{GL}^+(2,\RR)$ (Definition \ref{H mc E}). More precisely, we have:
\begin{coro}[of Lemmas 3.14,  3.16 in \cite{Macri}]  \label{coro of Macri} Let $\sigma =({\mc P}, Z) \in \st(\mc T)$.  Let $\mc E$ be a full Ext-exceptional collection
in $\mc T$.
Then we have the equivalences:

$\sigma \in \Theta_{\mc E}'\ \ \ \ \iff$  $\ \ \ \ \mc E \subset \mc P(t,t+1] $
for some $t \in \RR\ \ \ \ $  $\iff\ \ \ \ $  $\mc E$ is a
$\sigma$-exceptional collection.
\end{coro}
\bpr First, note  \cite[Lemma 3.16]{Macri} that  from $\{ E_i \}_{i=0}^n \subset  {\mc P}((t,t+1])$   it follows $ {\mc  A}_{\mc E}={\mc P}((t,t+1])$, and then  all $\{ E_i \}_{i=0}^n$ are stable in $\sigma$, because they are simple in  $ {\mc  A}_{\mc E}={\mc P}((t,t+1])$.  Indeed,  ${\mc  A}_{\mc E}$ and  ${\mc P}((t,t+1])$ are both bounded t-structures, therefore the inclusion $ {\mc  A}_{\mc E} \subset {\mc P}((t,t+1])$ implies equality $ {\mc  A}_{\mc E}={\mc P}((t,t+1])$. Whence, if $\{ E_i \}_{i=0}^n \subset  {\mc P}((t,t+1])$, then  ${\mc E}$ is $\sigma$-exceptional (see Definition \ref{sigma exceptional collection}).

Now the corollary follows from the last part of Bridgeland's Proposition
\ref{from abelian to triangulated} and the   following comments on
the action of $\widetilde{GL}^+(2,\RR)$. If $(\wt{{\mc P}}, \wt{Z}
)$ is obtained by the action with $\widetilde{GL}^+(2,\RR)$ on
$({\mc P}, Z )$, then $\{ \wt{{\mc P}}(\psi(t))={\mc P}(t)\}_{t
\in \RR}$ for some strictly increasing smooth function $\psi: \RR
\rightarrow \RR $ with $\psi(t + 1) =\psi(t ) + 1$, and hence ${\mc
P}(0,1]=\wt{\mc P}(\psi(0), \psi(0)+1]$. Conversely, for any $t \in
\RR $ and any  $({\mc P}, Z )$ we can act on it with element in
$\widetilde{GL}^+(2,\RR)$, so that the resulting $(\wt{{\mc P}},
\wt{Z} )$ satisfies ${\mc P}(t,t+1]=\wt{\mc P}(0, 1]$. \epr
Since the exceptional collection $\mc E$ in Definition \ref{sigma exceptional collection} has finite length, we have:
\begin{remark}  \label{equiv of open interval} The third condition in  Definition \ref{sigma exceptional collection} is equivalent to each of the following three conditions:
$\quad \{\phi(E_i)\}_{i=1}^n \subset (t,t+1)$ for some $t \in \RR$;  $\quad$  $\{\phi(E_i)\}_{i=1}^n \subset [t,t+1)$ for some $t \in \RR$; \\ $\max \left( \{\phi(E_i)\}_{i=0}^n \right ) - \min\left( \{\phi(E_i)\}_{i=0}^n \right )<1$.

 Furthermore, by Corollary \ref{coro of Macri} we have  $\Theta_{\mc E}'=\left \{\sigma  : \max \{\phi^\sigma_+(E_i)\}_{i=0}^n  - \min\{ \phi^\sigma_-(E_i)\}_{i=0}^n <1  \right \}=\left \{\sigma  : {\mc E}\subset \sigma^{ss} \ \mbox{and} \  \abs{\phi^\sigma(E_i)-\phi^\sigma(E_j)}<1   \ \mbox{for} \  i < j  \right \}$,
 therefore\footnote{For a a fixed nonzero object $X\in \mc T$ the functions $\sigma \mapsto \phi_\pm^\sigma(X)$ on
 the manifold $\mc T$ are continuous}  $\Theta_{\mc E}'$ is an open subset of $\st(\mc T)$.

 One can now easily  show  that the assignment:
 $$\Theta_{\mc E}'\ni \sigma = (\mc P,Z)\ \ \mapsto \ \ \left (\abs{Z(E_0)},\dots,\abs{Z(E_n)},\phi^\sigma(E_0),\dots,\phi^\sigma(E_n)\right )$$
 is well defined, and gives a homeomorphism between $\Theta_{\mc E}'$ and  the following  simply connected set:
 $$\left \{(x_0,\dots,x_n,y_0,\dots,y_n) \in \RR^{2(n+1)}\ \ : \ \ x_i>0,\ \abs{y_i-y_j}<1 \right \}.$$
\end{remark}
From the first part of this  remark and Corollary \ref{coro of Macri}  we see that for each $\sigma  \in \Theta_{\mc E}'$ we have an open interval, in which  $\mc P (x)$ is trivial (take $t\in \RR$
 and $\epsilon>0$ so that   $ \{\phi(E_i)\}_{i=0}^n  \subset (t,t+1]\cap  (t+\epsilon,t+\epsilon+1]$, then $(t,t+\epsilon)$ is such an interval). In particular(recall also that  ${\mathcal P}(x)[1]={\mathcal P}(x+1)$), we have:

\begin{remark} Let $\mc E$ be as in Corollary  \ref{coro of Macri}.  For each $\sigma  \in \Theta_{\mc E}'$ the set\footnote{see Lemma \ref{P_sigma} for the notation  $P_\sigma$}  $P_\sigma$ is not dense in $\SS^1$.
 \end{remark}

\section{Non-semistable exceptional objects in  hereditary abelian categories} \label{non-stable exc obj in...}

In this section   is written   an algorithm, denoted by $\mk{alg}$. In subsection \ref{pesumptions A is HER} we  define the
input data of the algorithm,  in subsection \ref{subs the cases} -
the data at the output. The rest sections of the text
refer   mainly to subsections \ref{pesumptions A is HER} and
\ref{subs the cases}.
\subsection{Presumptions} \label{pesumptions A is HER}  For  the rest of the paper ${\mathcal A}$ is    an abelian hereditary   $\hom$-finite  category, linear over an algebraically closed field $k$.\footnote{In all the sections  \ref{non-stable exc obj in...}, \ref{some terminilogy}, \ref{no bad after good}, \ref{sequence}, \ref{final}, \ref{constructing} the symbol  $\mc A$ denotes such a category.}  It can be shown\footnote{using some facts for modules over  unital associative  ring shown around page 302 of \cite{Lam}, see also \cite{Lenzing}}  that such a category   has  Krull-Schmidt property(Definition \ref{def of kr-schm}).
  Hence, by Lemma \ref{from Hered to derived with KSA}, the derived category  $D^b({\mathcal A})$ also satisfies the Krull-Schmidt property.
  For  brevity, we set ${\mathcal T}=D^b({\mathcal A})$.  Let $\sigma =({\mathcal P}, Z) \in \st({\mathcal T})$ be  a stability condition.
   In this setting  by Definition \ref{theta_sogma} we  obtain the function
   $\theta_\sigma : Ob({\mathcal T}) \rightarrow \NN^{(\sigma^{ss}_{ind}/\cong)}$.

 \textit{
 The input data of the algorithm $\mk{alg}$ is a non-semistable w. r. to $\sigma$
exceptional object $E\in \mc T$.  The output data is a
triangle, denoted by
 $ \mk{alg}(E)$.   We distinguish five  cases at the output, depending on the  features of the triangle  $ \mk{alg}(E)$, and  denote  them by \textbf{C1, C2, C3, B1, B2}. Only one of the five possible cases can occur at the output, i. e.  $ \mk{alg}(E)$ has all  the features of exactly  one case, say $X \in \{\textbf{C1, C2, C3, B1, B2}\}$, and then   $ \mk{alg}(E)$ is said to be of type $X$. }

 We note two facts, which we keep in mind further.
   \begin{remark}\label{from pre-exceptional to exceptional} It can be shown\footnote{by adapting the proof of this fact for quivers, given on \cite[p. 9,10]{WCB2}, to  $\mc A$}  that, under the given assumptions on $\mc A$,   if  $X\in \mc A_{ind}$ satisfies  ${\rm Ext}^1(X,X)=0$, then   $\Hom(X,X)=k$,  and hence  $X$ is an exceptional object.
\end{remark}
   \begin{remark} \label{remark for phi} Since $\mc A$ is a hereditary category, for any two indecomposable $A,B \in D^b(\mc A)$ with $\deg(A)=\deg(B)$ from $\phi_-(A)>\phi_+(B)+1$ it follows that $\Hom^*(A,B)=0$.
   \end{remark}
\noindent Another simple observation due to hereditariness, which
we will apply throughout, is:
  \begin{lemma} \label{lemmaHER} Let $\mathcal A$ be a hereditary abelian category and let
  $ 0 \rightarrow X \rightarrow Y \rightarrow Z \rightarrow 0 $ be a short exact sequence in $\mathcal A$.
  For each $W \in \mc A$  hold the following implications:
\begin{itemize}
    \item[\textbf{(a)}] If $\hom^1(Y,W)=0$, then $\hom^1(X,W)=0$
    \item[\textbf{(b)}] If $\hom^1(W,Y)=0$, then $\hom^1(W,Z)=0$.
\end{itemize}
  \end{lemma}
  \bpr To prove  \textbf{(a)} we apply $\Hom( \_,W[1])$ to the triangle  $X \rightarrow Y \rightarrow Z \rightarrow X[1]$, corresponding to the given exact sequence. It follows
$ 0=\Hom(Y,W[1]) \rightarrow \Hom(X, W[1])\rightarrow
\Hom(Z[-1],W[1])=0 $, where   the right vanishing is because
$\mathcal A$ is hereditary. In \textbf{(b)}  we apply
$\Hom(W,\_)$.
  \epr

We could work here with weaker assumptions on $\mc A$. More precisely:
\begin{remark} \label{weaker} Given that $\mc A$ is a hereditary $k$-linear abeilan category  with Krull Schmidt property as defined in
Definition \ref{def of kr-schm}, without assuming
$\hom$-finiteness and that $k$ is algebraically closed, then
everything in Sections \ref{non-stable exc obj in...}, \ref{some
terminilogy}, \ref{no bad after good}, \ref{sequence}, \ref{final}, \ref{constructing}
remains valid, if   we replace ``exceptional''   by
``pre-exceptional''.\footnote{By \textit{Pre-exceptional object}  we mean an indecomposable object $X\in
\mc T$ with $\Hom^i(X,X)=0$ for $i\neq 0$. \textit{Pre-exceptional
collection} is a sequence of pre-exceptional objects $(E_1, E_2,
\dots, E_n)$ with $\Hom^*(E_i,E_j)=0$ for $i>j$.}
 Under such seemingly weaker
 assumptions on $\mc A$,  we do not have the statement in Remark \ref{from pre-exceptional to exceptional}.
\end{remark}

\subsection{The cases} \label{subs the cases}
 Here we explain the features of each of the five  cases \textbf{C1, C2, C3, B1, B2} occurring at the output of $\mk{alg}$.  The other subsections  of \ref{non-stable exc obj in...}  contain  the  algorithm.

Let $E\in \mc T$ be a  non-semistable w. r. to $\sigma$
exceptional object. We recall that the meaning of the notation $\deg(E)$, used here, is explained in the paragraph \textbf{\textit{Some notations}} after the introduction.  The properties
(a),(b),(c) below are common features
of $\mk{alg}(E)$  for all the cases, property  (d) is
common for \textbf{C1, C2, C3}: \vspace{-5mm}
 \begin{itemize}\item[] \be \label{the diagram for all cases} \mk{alg}(E)= \begin{diagram}[size=1em]
U & \rTo      &     &       &   E \\
  & \luDashto &     & \ldTo &       \\
  &           & V &       &
\end{diagram} \qquad U\in \mc T,  V \in \sigma^{ss}, U \neq 0, V \neq 0, \ \mbox{where:}\ee
\item[(a)] $V$ is  the    degree $j$ component of\footnote{$\sigma_-(E)$ is defined in  Definition \ref{def of sigma_-}} $\sigma_-(E)$, where   $j\in \{\deg(E), \deg(E)+1\}$.
\item[(b)]
\ \ $ \theta_\sigma(U)<\theta_\sigma(E) \quad \Rightarrow \quad
\phi_-(U) \geq \phi(V)=\phi_-(E)  $.  \footnote{We write $f<g$ for
two functions $f,g \in \NN^{(\sigma^{ss}_{ind}/\cong) }$, if
$f(u)<g(u)$ for some $u \in \sigma^{ss}_{ind}/\cong$.}
\footnote{Note below that in cases \textbf{C3, B2} we have proper
inequality $\phi_-(U) > \phi(V)$.}
\item[(c)]  Any  $\Gamma \in Ind(V)$ satisfies $\hom(E,\Gamma)\neq 0$ (see Definition \ref{def of kr-schm} for the notation $Ind(V)$).
\item[(d)] In  the  cases \textbf{C1, C2, C3}  hold the vanishings  $\hom^*(U,V)=\hom^1(U,U)=\hom^1(V,V)=0$, in particular  for any $S \in Ind(V)$, $E' \in Ind(U)$ the pair $(S,E')$ is exceptional with $S\in\sigma^{ss}$.
\end{itemize}

We give now the complete lists of properties. For simplicity we assume that $E\in \mc A$, i. e. $\deg(E)=0$, for other degrees everything is shifted with the corresponding number.\\
{  \it {\rm \textbf{ C1}.}  The triangle is of the form $
\mk{alg}(E)=\begin{diagram}[size=1em]
A & \rTo      &     &       &   E \\
  & \luDashto &     & \ldTo &       \\
  &           & B &       &
\end{diagram}  $ with the properties:
\begin{itemize}
    \item[\textbf{C1.1}] $\{A,B\} \subset \mc A$, $A\neq 0$, $B\neq 0$, $\hom^1(A,A)=\hom^1(B,B)=\hom^*(A,B)=0$,
    \item[\textbf{C1.2}] $B$ is the zero degree component of $\sigma_-(E)$, in particular $B$ is semistable of phase $\phi_-(E)$,
    \item[\textbf{C1.3}] $\theta_\sigma(A) < \theta_\sigma(E)$ $\Rightarrow$ $\phi_-(A) \geq \phi_-(E)$,
    \item[\textbf{C1.4}]  any $\Gamma \in Ind(A)$ satisfies $\hom^1(B,\Gamma)\neq 0$.
 \end{itemize}
{\rm \textbf{ C2}.}  The triangle is of the form  \be \label{the C2 diagram}
\mk{alg}(E)=\begin{diagram}[size=1em]
A_1 \oplus A_2[-1] & \rTo      &     &       &   E \\
  & \luDashto &     & \ldTo &       \\
  &           & B &       &
\end{diagram}   \ee with the properties:
\begin{itemize}
    \item[\textbf{C2.1}] $\{ A_1,A_2,B \} \subset \mc A$, $A_2\neq 0$, $B\neq 0$, $A_1$ is a proper sub-object(in $\mc A$) of $E$,
    $\hom^1(A_2,A_2)=\hom^1(A_1,A_1)=\hom^*(A_1,B)=\hom^*(A_2,B)=\hom^*(A_1,A_2)=0$,
    \item[\textbf{C2.2}] $B$ is the zero degree component of $\sigma_-(E)$, in particular $B$ is semistable of phase $\phi_-(E)$,
    \item[\textbf{C2.3}] $\theta_\sigma(A_1)+ \theta_\sigma(A_2[-1])< \theta_\sigma(E)$, in particular $\phi_-(A_1)\geq \phi_-(E)$ and $ \phi_-(A_2[-1]) \geq \phi_-(E)$,
    \item[\textbf{C2.4}]  any  $\Gamma \in Ind(A_1)$ satisfies $\hom(B,\Gamma[1])\neq 0$,
     any  $\Gamma \in Ind(A_2)$ satisfies the three  conditions: $\hom(B,\Gamma) \neq 0$, $\hom(\Gamma,E[1]) \neq 0$, $\hom(E,\Gamma[1]) = 0$.
\end{itemize}\vspace{5mm}
{\rm \textbf{ C3}.}  The triangle is of the form \ $
\mk{alg}(E)=\begin{diagram}[size=1em]
A & \rTo      &     &       &   E \\
  & \luDashto &     & \ldTo &       \\
  &           & B[1] &       &
\end{diagram}   $ with the properties:
\begin{itemize}
    \item[\textbf{C3.1}] $\{A,B\} \subset \mc A$, $A\neq 0$, $B\neq 0$, $\hom^1(A,A)=\hom^1(B,B)=\hom^*(A,B)=0$,
    \item[\textbf{C3.2}]  $\mk{alg}(E)\cong {\rm H\! N}_-(E)$, hence
                                       $\theta_\sigma(A)< \theta_\sigma(E)$ and $\phi_-(A)> \phi_-(E)=\phi(B)+1$,
    \item[\textbf{C3.3}] any  $\Gamma\in  Ind(B)$ satisfies  $\hom^1(E,\Gamma)\neq 0$ and $\hom^1(\Gamma,E)= 0$,  any  $\Gamma\in Ind(A)$ satisfies   $\hom(B,\Gamma)\neq 0$ and $\hom(\Gamma,E)\neq 0$.
\end{itemize}\vspace{5mm}
{\rm \textbf{ B1}.} The  triangle is of the form $
\mk{alg}(E)=\begin{diagram}[size=1em]
A_1 \oplus A_2[-1] & \rTo      &     &       &   E \\
  & \luDashto &     & \ldTo &       \\
  &           & B &       &
\end{diagram}   $with the properties:
\begin{itemize}
    \item[\textbf{B1.1}] $\{A_1,A_2,B\} \subset \mc A$, $A_2\neq 0$, $B\neq 0$, $\hom^1(A_2,A_2)=\hom^1(A_1,A_1)=\hom^*(A_2,B)=0$, $A_1$ is a proper subobject(in $\mc A$)
    of $E$,
    \item[\textbf{B1.2}] $B$ is the zero degree component of $\sigma_-(E)$, in particular $B$ is semistable of phase $\phi_-(E)$,
    \item[\textbf{B1.3}] $\theta_\sigma(A_1)+ \theta_\sigma(A_2[-1])< \theta_\sigma(E)$, in particular $\phi_-(A_1)\geq \phi_-(E)$ and  $\phi_-(A_2[-1]) \geq \phi_-(E)$,
    \item[\textbf{B1.4}] there exists  $\Gamma \in Ind(A_2)$ with  $\hom^1(\Gamma,E)\neq 0$, $\hom^1(E,\Gamma)\neq 0$.\footnote{A comparison  with  \textbf{C2.4} shows  that \textbf{B1} and \textbf{C2} cannot appear together.}
\end{itemize}
{\rm \textbf{ B2}.} The  triangle is of the form $ \label{triangle B2}
\mk{alg}(E)=\begin{diagram}[1em]
A & \rTo      &     &       &   E \\
  & \luDashto &     & \ldTo &       \\
  &           & B[1] &       &
\end{diagram}   $ with the properties:
\begin{itemize}
    \item[\textbf{B2.1}] $\{A,B\} \subset \mc A$, $A\neq 0$, $B\neq 0$, $\hom^1(B,B)=\hom^*(A,B)=0$,
    \item[\textbf{B2.2}] $\mk{alg}(E)\cong {\rm H\! N}_-(E)$, hence
                                       $\theta_\sigma(A)< \theta_\sigma(E)$ and $\phi_-(A)> \phi_-(E)=\phi(B)+1$,
    \item[\textbf{B2.3}] there exists  $\Gamma\in B$ with $\hom^1(\Gamma,E)\neq 0$, $\hom^1(E,\Gamma)\neq 0$.\footnote{A comparison  with  \textbf{C3.3} shows  that \textbf{B2} and \textbf{C3} cannot appear together.}
\end{itemize}
   }
\subsection{The last HN triangle } Now we start explaining  $\mk{alg}$.

 Let   $E \in \mc A_{exc}$, $E \not \in \sigma^{ss}$.
  Macr\`i initiated in \cite[p. 10]{Macri} an analysis  of the last HN triangle of $E$, when $E \in  Rep_k(K(l))$.  The arguments  on \cite[p. 10]{Macri} are used here in formulas  \eqref{the last HN tr1}, \eqref{hom^leq(X,-)=0}, and in the derivation of the vanishings  \textbf{C3.1}(Subsection \ref{if B_0=0}).

  Consider the last HN triangle ${\rm H\! N}_-(E)$(see Definition \ref{def of sigma_-}):
 \begin{gather} \label{the last HN tr} {\rm H\! N}_-(E) \ = \ \begin{diagram} X & \rTo & E  & \rTo^f &  \sigma_-(E) & \rTo  & X[1] \end{diagram}, \qquad \qquad \phi_-(X) > \phi(\sigma_-(E))=\phi_-(E). \end{gather}

 \begin{lemma} \label{the last HN triangle} The triangle ${\rm H\! N}_-(E)$ is of the form (with $B_0,B_1 \in {\mathcal A}$): 
  \begin{gather} \label{the last HN tr1} \  \begin{diagram} X & \rTo & E  & \rTo^f &  B_0 \oplus B_1[1] & \rTo  & X[1] \end{diagram}, \quad \phi_-(X) > \phi(B_0)=\phi(B_1)+1=\phi_-(E), \\
 \label{hom^leq(X,-)=0} \hom^{\leq 0}(X,B_0)=\hom^{\leq 0}(X,B_1[1])=0 \\
   \label{theta_sigma(E)=theta_sigma(X)+...} \theta_\sigma(E)=\theta_\sigma(X)+\theta_\sigma(B_0)+\theta_\sigma(B_1[1]).\end{gather}

For any $i\in \{0,1 \}$,  $\Gamma \in Ind(B_i)$ the component of $f$
to $\Gamma[i]$ is non-zero and  $\hom(E,\Gamma[i])\neq 0$.

Any  $\Gamma \in Ind(X)$  satisfies  $\hom(\Gamma,E)\neq 0$ and
$\hom(B_0\oplus B_1[1], \Gamma[1])\neq 0$.
\end{lemma}
\bpr We show first that for each $\Gamma \in Ind(\sigma_-(E))$ the component of $f$ from $E$ to $\Gamma$ is non-zero.

Indeed,
 suppose that for some $\Gamma \in Ind(\sigma_-(E))$ this  component vanishes, then by the Krull-Schmidt property we can write  $\sigma_-(E)=U\oplus\Gamma$,
  and   $f$ is of the form: $f=(f':E\rightarrow U) \oplus (0
\rightarrow\Gamma).$ After summing the triangles
$\begin{diagram} X' & \rTo & E  & \rTo^{f'} &  U & \rTo  & X'[1]
\end{diagram}$  and $\begin{diagram} \Gamma [-1] & \rTo & 0  & \rTo &
\Gamma & \rTo  & \Gamma \end{diagram}$ we obtain a
triangle $ \begin{diagram} X'\oplus \Gamma[-1] & \rTo &
E  & \rTo^{f} &  \sigma_-(E) & \rTo  & X'[1]\oplus \Gamma
\end{diagram} $ (recall that $E\neq \sigma^{ss}$, hence $X'\neq
0$). From \eqref{the last HN tr} it follows that $X'\oplus
\Gamma[-1] \cong X$. From Corollary \ref{ineq for phi_-} we
see that
$\phi_-(X') \geq \phi_-(X) > \phi_-(E)=\phi(U).$ By this inequality and  the uniqueness of the HN filtration of $E$  we deduce that
$\sigma_-(E)\cong U$, i. e. $U \oplus \Gamma \cong U$,
which contradicts the Krull-Schmidt property.

Thus, for each $\Gamma \in Ind(\sigma_-(E))$ the component of $f$ to $\Gamma$ is non-zero and $\hom(E, \Gamma)\neq 0$. Now the triangle
\eqref{the last HN tr} reduces to  \eqref{the last HN tr1},  since $\mathcal A$ is hereditary. From $\phi_-(X) > \phi(B_i[i])$ ($i=0,1$) it
follows \eqref{hom^leq(X,-)=0}. Applying Lemmas \ref{additive property of theta_sigma 1}, \ref{additive
property of theta_sigma 2} to \eqref{the last HN tr1} we obtain
\eqref{theta_sigma(E)=theta_sigma(X)+...}.
 It remains to prove the
last  property.

Suppose that $\hom(\Gamma,E)=0$ for some  $\Gamma \in Ind(X)$. Then we can represent $X
\rightarrow E$ as a direct sum  $(U\rightarrow E) \oplus (\Gamma
\rightarrow 0).$  By the triangle \eqref{the last HN tr}  we get
 $Y'\oplus \Gamma[1] \cong \sigma_-(E)$, where $Y'$ is the cone of $U\rightarrow E$. From Corollary \ref{ineq for phi_-} we see
$\phi_-(U) \geq \phi_-(X) > \phi_-(E)=\phi(Y').$  Since
$\phi_-(U)>\phi_-(E)$, we have $U \not \cong E$ and $Y'
\neq 0$. It follows that ${\rm HN}_-(E)$ $=$ $\bd[1em] U & \rTo & E & \rTo & Y'& \rTo & U[1]  \ed$. Therefore   $X\cong
U$, i. e. $U \oplus \Gamma \cong U$, which contradicts the
Krull-Schmidt property.

Suppose that for some $\Gamma\in Ind(X)$ we have $\hom(B_0\oplus
B_1[1], \Gamma[1])= 0$, then by  similar arguments  we get $E\cong
E'\oplus \Gamma $, and hence $\Gamma \cong E$ (since $E$ is
indecomposable),  which contradicts $\phi_-(\Gamma) \geq
\phi_-(X)>\phi_-(E)$. The lemma is proved. \epr

 By  $ f_i $ will be denoted the component of $f$ to $B_i[i]$(see  \eqref{the last HN tr1}),  i. e.  we have
commutative diagrams (the right arrow is the projection)
\be \label{f_i=proj circ f} \begin{diagram}[1.5em]  E      & \rTo^f     &  B_0 \oplus B_1[1]  \\
                                              \dTo^{Id} &            &    \dTo           \\
                                              E      & \rTo^{f_i} &  B_i[i]   \\
                                \end{diagram} \qquad i\in \{0,1\}. \ee
  The algorithm $\mk{alg}$  tests now the condition $B_0=0$.
\subsection{If \texorpdfstring{$B_0=0$}{\space}}\label{if B_0=0} \mbox{}

 This condition  leads to one of the cases \textbf{C3, B2} depending on the outcome of one test.
 Since $B_0=0$,  the  triangle \eqref{the last HN tr1} is reduced to a short exact sequence  $
  \begin{diagram}[1em] 0& \rTo & B_1 & \rTo & X & \rTo & E & \rTo & 0 \end{diagram},  $  and  $X \in Ob({\mc A})$.
 Hence \eqref{hom^leq(X,-)=0} is now the same as
 $ \hom^*(X,B_1)=0,$ which by Lemma \ref{lemmaHER} \textbf{(a)} and the given exact sequence implies
 $ \hom^1(B_1,B_1)=0. $
 By Lemma \ref{the last HN triangle}  any  $\Gamma\in Ind(B_1)$ satisfies  $\hom(E,\Gamma[1])\neq 0$.
 Therefore,  if $\hom(\Gamma,E[1])\neq 0$ for some  $\Gamma\in Ind(B_1)$, then the
 triangle:
 \be \label{candidate for B2,C3}{\rm HN}_-(E)= \begin{diagram}[size=1em]
X & \rTo      &     &       &   E \\
  & \luDashto &     & \ldTo &       \\
  &           & B_1[1] &       &
\end{diagram}   \ee
 satisfies \textbf{B2.1, B2.2, B2.3} (with $A=X, B=B_1$). By setting $\mk{alg}(E)$ to  \eqref{candidate for B2,C3}  we get \textbf{B2}.

 It  remains to consider the case when  $\hom(\Gamma,E[1])= 0$ for each  $\Gamma \in Ind(B_1)$, i. e.
 \be \label{proof1 C3} \hom(B_1,E[1])= 0.\ee
Setting again $\mk{alg}(E)$ to   \eqref{candidate for B2,C3}(with $X$ replaced by $A$, $B_1$ replaced by $B$) we obtain  the property \textbf{C3.2} immediately.
 The property \textbf{C3.3} follows from Lemma \ref{the last HN triangle}. We have already all the features of \textbf{C3.1} except   the vanishing $\hom^1(X,X)=0$.

 The vanishing  $\hom^1(X,X)=0$ follows from \eqref{proof1 C3}, since  the triangle  \eqref{candidate for B2,C3} and  $\Hom(X,\_)$ give  an exact sequence
  $ \Hom^1(X,B_1)\rightarrow \Hom^1(X,X)\rightarrow \Hom^1(X,E) $,
  where the left and the right terms vanish.  The  vanishing $\Hom^1(X,B_1)=0$ is already shown (before \eqref{candidate for B2,C3}). The other vanishing
  $\hom^1(X,E) =0$  follows from  \eqref{proof1 C3}, $\hom^1(E,E)=0$, and    $\Hom(\_,E[1])$ applied to the same triangle.

 Thus, $\mk{alg}(E)$ is of type \textbf{C3}.

\subsection{If \texorpdfstring{$B_0 \neq 0$}{\space}} \label{if B_0 neq 0}  \mbox{}

 Under this condition we obtain  one of the cases \textbf{C1, C2, B1} at the output depending on the outcomes of additional tests.

 By Lemma \ref{the last HN triangle} we have $f_0 \neq 0$. Let us take kernel and cokernel of $f_0$ in $\mathcal A$: \be
\label{A_1 and A_2} \begin{diagram} A_1 & \rTo^{ker(f_0)} & E &
\rTo^{f_0} & B_0 & \rTo^{coker(f_0)} & A_2.  \end{diagram} \ee
Since $f_0 \neq 0$,   $ker(f_0)$ is a proper
subobject of $E$. Let  $  \begin{diagram}[1em] E & \rTo^{e_0} & B_0' & \rTo^{im(f_0)} & B_0 \end{diagram}  $ be a   decomposition of $f_0$ in $\mc A$,
where  $e_0$ is $\mc A$-epic and $im(f_0)$ is $\mc A$-monic. In particular, we have an exact sequence in $\mc A$
 \be \label{A_1 rto E rto B_0'}  \begin{diagram} 0& \rTo & A_1 & \rTo^{ker(f_0)} & E & \rTo^{e_0} & B_0' & \rTo & 0 \end{diagram}. \ee
  The next step of the algorithm $\mk{alg}$ is to  test the condition $A_2=0$. We show first  some preliminary facts, which do not depend on the vanishing of $A_2$.
 \subsubsection{Preliminary facts} These  facts are   \eqref{van1},\eqref{van2},\eqref{van3},  \eqref{f_0 neq 0 theta_sigma<..}, and Lemma \ref{for pro C1.4 C2.4}.

 The  equalities below will help us later to obtain \textbf{C1.1}, when $A_2= 0$, and \textbf{C2.1}, when $A_2\neq 0$:
\be \label{van1} \hom^1(A_1,A_1)= \hom^1(A_2,A_2)=0 \\
     \label{van2}  \hom(A_1,B_0)= \hom^*(A_2,B_0)=0 \\
    \label{van3} \hom^1(A_1,E)= \hom(A_1,A_2)=0 \ee
    The  inequality \eqref{f_0 neq 0 theta_sigma<..}  ensures \textbf{C1.3} and \textbf{C2.3}, and Lemma \ref{for pro C1.4 C2.4}   ensures      \textbf{C1.4} and
    half of  \textbf{C2.4}.
    \be \label{f_0 neq 0 theta_sigma<..} \theta_{\sigma}(A_1) +\theta_{\sigma}(A_2[-1]) < \theta_\sigma(E).\ee

To show these facts we start by recalling that  the  triangle in $\mathcal T$ containing $f_0$ is
\be \label{triangle woth f_0} \begin{diagram} E& \rTo^{f_0} & B_0
& \rTo & C(f_0) & \rTo E[1] \end{diagram} \ee where $C(f_0)$ is
the cochain complex ($B_0$ is in degree $0$) \be \label{C(f_0)}
\begin{diagram}  & \dots & \rTo & 0 & \rTo & E & \rTo^{f_0} & B_0
& \rTo & 0 & \rTo & \dots \end{diagram} \ee and the non-trivial part of the cochain maps $ B_0
\rightarrow  C(f_0) \rightarrow E[1]$ is $
 \begin{diagram}[1em]
                 0   & \rTo       & E         & \rTo &  E \\
                 \dTo  &            & \dTo  &      & \dTo \\
                 B_0 & \rTo  & B_0       & \rTo & 0  \\
                \end{diagram}.
$\\ Since ${\mathcal A}$ is  hereditary, we have $C(f_0) \cong
\bigoplus_{i} H^i(C(f_0))[-i]$, which we can reduce  by \eqref{A_1
and A_2} and  \eqref{C(f_0)} to \be \label{C(f_0) cong ...} C(f_0)
\cong  A_1[1] \oplus A_2. \ee

Since we have the commutative diagram \eqref{f_i=proj circ f} with $i=0$, by the $3\times 3$ lemma in triangulated categories \cite[Proposition 1.1.11]{BBD}
 we can  put the
 triangles \eqref{the last HN tr1}, \eqref{triangle woth f_0} in a diagram
\begin{diagram}[size=1.5em] E      & \rTo^f     &  B_0 \oplus B_1[1] & \rTo & X[1]   & \rTo  & E[1]      \\
             \dTo^{Id} &            &    \dTo            &      & \dTo   &       & \dTo^{Id} \\
                E      & \rTo^{f_0} &  B_0               & \rTo & C(f_0) & \rTo  & E[1]      \\
             \dTo      &            &    \dTo^0          &      & \dTo   &       & \dTo      \\
                0      & \rTo       &  B_1[2]            & \rTo & Y      & \rTo  & 0         \\
             \dTo      &            &    \dTo            &      & \dTo   &       & \dTo      \\
                E[1]   & \rTo^{f_0} &  B_0[1]\oplus B_1[2]& \rTo & X[2]  & \rTo  & E[1]      \\
\end{diagram}
where $ \begin{diagram}[1em]  X[1]     & \rTo     & C(f_0) & \rTo &  Y
& \rTo  & X[2] \end{diagram} $;   $ \begin{diagram}[1em] 0      & \rTo
&  B_1[2]            & \rTo & Y      & \rTo  & 0 \end{diagram} $
 are distinguished triangles. Hence $Y \cong B_1[2]$ and we obtain a distinguished triangle
\be \label{help dist triangle}  \begin{diagram} X     & \rTo     &
C(f_0)[-1] & \rTo &  B_1[1]   & \rTo  & X[1].  \end{diagram} \ee
The vanishings \eqref{van1},\eqref{van2},\eqref{van3} will be obtained from   triangles \eqref{help dist triangle},\eqref{triangle woth f_0}, and the exact
sequence \eqref{A_1 rto E rto B_0'}.

We apply $\Hom(-,B_0)$ and $\Hom(-,B_0[-1])$ to \eqref{help dist
triangle} and by \eqref{hom^leq(X,-)=0} the result is:
$\Hom(C(f_0),B_0[1])=\Hom(C(f_0),B_0)=0$. These vanishings and
\eqref{C(f_0) cong ...} imply \eqref{van2}. The vanishing
$\hom^1(A_1,E)=0$ (the first part of \eqref{van3}) follows from
$\hom^1(E,E)=0$, the exact sequence \eqref{A_1 rto E rto B_0'},
and Lemma \ref{lemmaHER} \textbf{(a)}. Now we can write
 $\hom(C(f_0),E[2])=\hom(A_1[1] \oplus A_2,E[2])=\hom^1(A_1,E)=0$. Having  $0=\hom(C(f_0),B_0[1])=\hom(C(f_0),E[2])$, we  apply $\Hom(C(f_0),\_)$
  to \eqref{triangle woth f_0} and obtain
$$ 0=\Hom(C(f_0),B_0[1]) \rightarrow \Hom(C(f_0),C(f_0)[1]) \rightarrow \Hom(C(f_0),E[2])=0.  $$
Hence $\hom(A_1[1]\oplus A_2, A_1[2]\oplus A_2[1] )=0$, which
contains \eqref{van1} and the second vanishing in \eqref{van3}.

The next step is to show \eqref{f_0 neq 0 theta_sigma<..}. From Lemma \ref{additive property of theta_sigma 1} and
the triangle \eqref{help dist triangle} we get
$\theta_\sigma(C(f_0)[-1])= \theta_\sigma(X) +
\theta_\sigma(B_1[1])$. From  $B_0 \neq 0$ it follows
$\theta_\sigma(B_0) > 0$, and hence:
$$ \theta_\sigma(C(f_0)[-1])= \theta_\sigma(X) + \theta_\sigma(B_1[1]) < \theta_\sigma(X) + \theta_\sigma(B_1[1])+\theta_\sigma(B_0) =\theta_\sigma(E),$$
where the last equality is taken from \eqref{theta_sigma(E)=theta_sigma(X)+...}. Now \eqref{f_0 neq 0 theta_sigma<..} follows from \eqref{C(f_0)
cong ...}.

Since $\mk{alg}(E)$ in both the cases $A_2 =0$ and $A_2 \neq 0$ will be set to  \eqref{triangle woth f_0},   the following corollary ensures  \textbf{C1.4},
and part of  \textbf{C2.4}.
\begin{lemma} \label{for pro C1.4 C2.4} Each  $\Gamma \in Ind(C(f_0)) =Ind( A_1[1]\oplus A_2)$ satisfies $\hom(B_0,\Gamma)\neq 0$, and each $\Gamma \in Ind(A_2)$
satisfies $\hom(\Gamma,E[1])\neq
0$.
\end{lemma}
\bpr Suppose that $\hom(B_0,\Gamma)=0 $ for some  $\Gamma \in Ind(C(f_0))$ and split $C(f_0)=U\oplus
\Gamma$, then the arrow $B_0 \rightarrow C(f_0)$ in
\eqref{triangle woth f_0} can be represented as $(B_0 \rightarrow U)\oplus(0 \rightarrow \Gamma)$. The
 sum of the triangle $ \begin{diagram}[1em] E'& \rTo & B_0 & \rTo & U & \rTo E'[1] \end{diagram} $ extending $B_0 \rightarrow U$
 and the triangle
$ \begin{diagram}[1em] \Gamma[-1] & \rTo & 0 & \rTo & \Gamma & \rTo \Gamma \end{diagram}   $
is  isomorphic to \eqref{triangle woth f_0}, hence
$E\cong E'\oplus \Gamma[-1]$. Since $E$ is exceptional and
$\Gamma\neq 0$, it follows $E'=0$ and $E\cong \Gamma[-1]$, hence
$\theta_\sigma(E)= \theta_\sigma(\Gamma[-1])\leq \theta_\sigma(C(f_0)[-1])<\theta_\sigma(E), $
where we used $C(f_0)=U\oplus
\Gamma$ and the inequality derived before this corollary. Thus,
we get a contradiction.

If $\hom(\Gamma, E[1]) =0$ for
some $\Gamma \in Ind(A_2)$, then  we can split $C(f_0)\cong A_1[1] \oplus A_2
\cong V \oplus \Gamma$, and  the last arrow in
\eqref{triangle woth f_0} is of the form
 $(V \rightarrow E[1])\oplus (\Gamma \rightarrow 0) $. It follows by similar arguments as above that  $B_0  \cong  U \oplus
\Gamma$ for some $U$. Therefore $\hom(A_2,
B_0)\neq 0$, which contradicts \eqref{van2}
\epr

\subsubsection{If $A_2=0$} Under this  condition  we  get here a triangle or  type \textbf{C1}.

Now  $f_0$ is epic(see \eqref{A_1 and A_2}) and \eqref{triangle woth f_0} becomes a  short exact sequence \be \label{sesC1} \bd  0
& \rTo &  A_1 &  \rTo^{ker(f_0)} & E & \rTo^{f_0} & B_0 & \rTo & 0.
\ed \ee
The triangle $\mk{alg}(E)$ is set  to   \eqref{sesC1},  so $A=A_1$, and $B=B_0$.  From
\eqref{f_0 neq 0 theta_sigma<..} we get $\theta_\sigma(A_1) <
\theta_\sigma(E)$, which is the same as \textbf{C1.3}. In Lemma \ref{for pro C1.4 C2.4} we have \textbf{C1.4},   and in \eqref{the last HN triangle} -  \textbf{C1.2}. It remains to show \textbf{C1.1}. We have $A_1\neq 0$, for otherwise $E$ would be
semistable.  We have also \eqref{van1} and \eqref{van2}, therefore we have to show only $\hom^1(A_1, B_0)=0=\Hom^1(B_0,B_0)$.

 By $\hom^1(A_1, E)=0$ (see \eqref{van3}), the sequence \eqref{sesC1}, and Lemma
\ref{lemmaHER} \textbf{(b)} we obtain $\hom^1(A_1, B_0)=0$.  The
same lemma and $\hom^1(E,E)=0$ imply $\hom^1(E,B_0)=0$, hence
$\Hom(\_,B_0[1])$ applied to \eqref{sesC1} gives: $
0=\Hom(A_1[1],B_0[1]) \rightarrow \Hom(B_0,B_0[1]) \rightarrow
\Hom(E,B_0[1])=0 $, i. e. $\hom^1(B_0,B_0)=0 .$

\subsubsection{If $A_2 \neq 0$.} \mbox{} \\
Under this condition we will obtain either the   case \textbf{C2}
or the  case \textbf{B1} depending on the outcome of one
additional test. The triangle $\mk{alg}(E)$  is set to
\eqref{triangle woth f_0}, which by $C(f_0) \cong A_1[1]\oplus
A_2$  can be rewritten as: \be \label{candidate for B1,C2}
\mk{alg}(E)=\begin{diagram}[1em]
A_1 \oplus A_2[-1] & \rTo      &     &       &   E \\
  & \luDashto &     & \ldTo &       \\
  &           & B_0 &       &
\end{diagram}   \ee

From Lemma \ref{for pro C1.4 C2.4}   we have $\hom^1(\Gamma,E)\neq
0$  for each $\Gamma \in Ind(A_2)$.
If $\hom^1(E,\Gamma)\neq 0$ for some  $\Gamma \in Ind(A_2)$,  then
the triangle \eqref{candidate for B1,C2} has all the
features of the case \textbf{B1} due to \eqref{A_1 rto E rto B_0'},
\eqref{van1}, \eqref{van2}, \eqref{f_0 neq 0 theta_sigma<..}.

Thus, it remains to show that if each  $\Gamma \in Ind(A_2)$ satisfies
$\hom(E,\Gamma[1])= 0$, in particular  \be \label{Hom(E,A_2[1])= 0}
\Hom(E,A_2[1])= 0 \ \ \ \Rightarrow \ \ \
\Hom(E,C(f_0)[1])=\Hom(E,A_1[2]\oplus A_2[1])=0, \ee then
  the triangle \eqref{candidate for B1,C2}  satisfies \textbf{C2.1},   \textbf{C2.2}, \textbf{C2.3},   \textbf{C2.4} (with $B=B_0$).

\textbf{C2.2} is in \eqref{the last HN tr1},  \textbf{C2.3} is \eqref{f_0 neq 0
theta_sigma<..}, and \textbf{C2.4} is contained in
\eqref{Hom(E,A_2[1])= 0}, Lemma \ref{for pro C1.4 C2.4}.  It remains to obtain  the vanishings
in \textbf{C2.1}, that are not claimed  in \eqref{van1}, \eqref{van2}, \eqref{van3}. These vanishings are $\hom^1(A_1,B_0)=\hom^1(A_1,A_2)=\hom^1(B_0,B_0)=0$.
We obtain them in this order below.

The  equality \eqref{Hom(E,A_2[1])= 0} together with
$\hom^1(E,E)=0$ and the triangle \eqref{triangle woth f_0} imply
$$ \hom(E,B_0[1])=0, $$ hence by the sequence \eqref{A_1 rto E rto
B_0'} and Lemma \ref{lemmaHER} we get $\hom^1(A_1,B_0)=0$.  From this vanishing it follows $\hom(C(f_0),B_0[2])=0$ and applying $\Hom(C(f_0),\_)$ to \eqref{triangle woth f_0} we obtain
 \begin{gather} 0=\Hom(C(f_0), B_0[2]) \rightarrow \Hom(C(f_0), C(f_0)[2]) \rightarrow \Hom(C(f_0), E[3])=0 \nonumber  \\
  \Rightarrow 0=\hom(C(f_0), C(f_0)[2])= \hom(A_1[1]\oplus A_2, A_1[3]\oplus A_2[2])  \ \ \Rightarrow \ \ \hom(A_1,A_2[1])=0. \nonumber \end{gather}

Finally, we apply $\Hom(\_,B_0[1])$ to
\eqref{triangle woth f_0}: $0= \Hom(C(f_0),B_0[1]) \rightarrow \Hom(B_0,B_0[1]) \rightarrow \Hom(E,B_0[1])=0$,
where the left vanishing is contained in \eqref{van2}, and the
right vanishing is above.

Now we have already the complete list \textbf{C2} for $\mk{alg}(E)$.

\section{Some terminology. The relation \texorpdfstring{$ R  \DashedArrow[densely dotted] (S,E) $}{\space} } \label{some terminilogy}
The terminology introduced here is important for the
 rest of the paper. All  definitions  in this section assume a  given
stability condition on\footnote{Here and in all sections that
follow  $\mc A$ is as in subsection \eqref{pesumptions A is HER}.}
$D^b(\mc A)$, which we  denote by $\sigma $. We divide the  non-semistable exceptional objects into
two types: $\sigma$-regular and $\sigma$-irregular (Definition
\ref{def regular irregular}). In turn the $\sigma$-regular objects
are divided into final and non-final (Definition \ref{def final
good}).

  We refer to \textbf{C1, C2, C3} as regular  cases   and to \textbf{B1, B2} as irregular cases.
More precisely:
\begin{df}\label{def regular irregular} Let $E \in D^b(\mc A)_{exc}$ and $E\not \in \sigma^{ss} $. If
 the triangle $\mk{alg}(E)$ given by  section \ref{non-stable exc obj in...} is of type  $X$, where $X$ is one  of  \textbf{C1,C2,C3,B1,B2}
, then  $E$ is said to be an  $X$ object
w. r. to $\sigma$.

The \textbf{Ci} objects(for $i=1,2,3$)  will be  called
\uline{$\sigma$-regular exceptional objects}  and the \textbf{Bi}
objects( for $i=1,2$) will be called \uline{$\sigma$-irregular
exceptional objects}.\footnote{In this text the adjectives
``$\sigma$-regular'', ``$\sigma$-irregular'' regard either
exceptional objects or the cases at the output of $\mk{alg}$. We
often omit ``exceptional object'' after these adjectives, when
this is by default. We sometimes omit ``$\sigma$-'', which is akin
to writing semistable instead of $\sigma$-semistable.}
\end{df}
We  introduce now the relation $\begin{diagram} R & \rDotsto^{X} &
(S,E) \end{diagram}$. It facilitates the next steps of the exposition.
\begin{df} \label{def of Dotsto}  Let $R,S,E \in D^b(\mc A)$ and let $X$ be one of the symbols   \textbf{C1,C2a,C2b,C3}.
By the notation $\begin{diagram} R & \rDotsto^{X} & (S,E) \end{diagram}$ we mean   the following data:
\begin{itemize}
    \item $R$ is a  $\sigma$-regular exceptional object, in particular  $\mk{alg}(R)$  is of type \textbf{Ci}($i\in \{ 1,2,3\}$),
    \item $S \in Ind(V)$, $E \in Ind(U)$, where $(V,U)$ are the  lower and the left  vertices  of $\mk{alg}(R)$   in \eqref{the diagram for all cases},
        \item if  $i\in \{1,3\}$ and  $R$ is  a \textbf{Ci} object, then we set $X= \textbf{Ci}$,
    \item if $R$ is a \textbf{C2} object and $E$ is a component of $A_2[-1]$ in diagram \eqref{the C2 diagram}, then we set $X= \textbf{C2a}$,
  \item   if $R$ is a \textbf{C2} object and $E$ is a component of $A_1$ in \eqref{the C2 diagram}, then we set $X= \textbf{C2b}$.
\end{itemize}
\end{df}
In the next sections we refer mainly  to the following
features(explained  below) of the pair  $(S,E)$:
 \begin{gather}\nonumber  \begin{diagram} R & \rDotsto^{X} & (S,E) \end{diagram} \ \ \  X \in \{ \textbf{C1, C2a, C2b, C3} \}\\
 \label{common properties for Dotsto 1} \{S,E \} \subset D^b({\mc A})_{exc}, \  \ \hom^*(E,S)=0, \  \ \deg(E)+1 \geq \deg(S)\geq \deg(R)\geq\deg(E) \\
\label{common properties for Dotsto
2}\theta_\sigma(E)<\theta_\sigma(R), \quad S \in \sigma^{ss},
\quad \theta_\sigma(R)(S)>0, \quad \phi_-(E)\geq
\phi(S)=\phi_-(R).\end{gather}

 The first two statements in \eqref{common properties for Dotsto 1} amount to   saying that $(S,E)$ is an exceptional pair,\footnote{In general, this pair is not uniquely determined by $R$, because we make  choices among  $Ind(U)$ and $Ind(V)$.} which is  the same as: $S$, $E$ are indecomposable and $\hom^*(E,S)=\hom^1(S,S)=\hom^1(E,E)=0$.
 This follows from  (d)  right after \eqref{the diagram for all cases} and $S \in Ind(V)$, $E \in Ind(U)$.
  In (a) right after \eqref{the diagram for all cases} is specified  that  $V$  is a direct summand of  $\sigma_-(R)$,    hence by $S\in Ind(V)$  and the
   definition of $\theta_\sigma$(Definition \ref{theta_sogma})  it follows that $\theta_\sigma(R)(S)>0$ and $S\in \sigma^{ss}$.  In (b) right after \eqref{the diagram for all cases} we have specified  $\theta_\sigma(U)<\theta_\sigma(R)$, $\phi_-(U) \geq \phi(V)=\phi_-(R)$, which by   $E\in Ind(U)$, $S \in Ind(V)$ implies  $\theta_\sigma(E)<\theta_\sigma(R)$, $\phi_-(E)\geq \phi(S)=\phi_-(R)$. Thus we obtain \eqref{common properties for Dotsto 2}.
  The degrees of $R,S,E$ are interrelated as shown  in the following table,\footnote{Recall that for $X\in \mc A$ and $j\in \ZZ$ we write $\deg(X[j])=j$.} which
   follows from the very definition of \textbf{C1,C2a,C2b,C3}:\footnote{the description of   \textbf{C1, C2, C3} is in subsection \ref{subs the cases}}
\begin{gather} \label{table with degrees} \begin{array}{ c | c | c | c}
      X     & \deg(S)- \deg(R) & \deg(R)- \deg(E)  &                          \\ \hline
    \textbf{C1, C2}b & 0              & 0               & \phi_-(E) \geq \phi(S) \\ \hline
    \textbf{C2a}     & 0              & +1              & \phi_-(E) \geq \phi(S) \\ \hline
   \textbf{ C3}      & +1             & 0               &  \phi_-(E) >\phi(S)
  \end{array}
  \end{gather}
The inequalities $\deg(E)+1 \geq \deg(S)\geq
\deg(R)\geq\deg(E)$ follow,  so  \eqref{common properties for Dotsto
1} is   shown completely.

We divide the  $\sigma$-regular objects into final and non-final as follows:
  \begin{df} \label{def final good} If $R$ is a $\sigma$-regular object  and all the  indecomposable components of $U$ (in  diagram
  \eqref{the diagram for all cases}) are  semistable, then $R$ is said to be  \uline{final}, otherwise - \uline{non-final}.
\end{df}
 If  $R$ is a non-final regular  object then  some indecomposable component of $U$ is not  semistable. By regularity this component is also an exceptional object  and then we can apply to it $\mk{alg}$.  Now we cannot exclude the occurrence of the irregular cases \textbf{B1, B2}, i. e. we cannot exclude the occurrence of an irregular  component of  $U$.
 \section{Regularity-preserving categories. RP prpoerties 1,2 } \label{no bad after good}
Recall that $\mk{alg}$ can be applied to any non-semistable
exceptional object.  Using the terminology  from
Section \ref{some terminilogy}, we can say that if   $R$ is   $\sigma$-regular and
non-final,  then from the output data  $\mk{alg}(R)$ we can extract
some  number of non-semistable exceptional objects (the non-semistable
components of $U$ in  diagram
  \eqref{the diagram for all cases}). The algorithm  $\mk{alg}$   can be applied    to any  of them  again.  If the category $\mc A$ has the property that  the cases \textbf{ B1, B2} cannot occur after this second iteration of $\mk{alg}$ we say that $\mc A$ is regularity-preserving. More precisely:
 \begin{df} \label{RP category} A hereditary abelian category $\mc A$ will be said to be \ul{regularity-preserving}, if for each $\sigma \in \st(D^b(\mc A))$ from the
 the following  data:

 $R\in D^b(\mc A)$ is a $\sigma$-regular object; $ \begin{diagram} R & \rDotsto^{X} & (S,E) \end{diagram}$, where $ X \in \{ \textbf{C1, C2a, C2b, C3} \}$; $E \not \in \sigma^{ss}$ \vspace{2mm}

it follows that  $E$   is  a  $\sigma$-regular object as well.
 \end{df}

In this  section \ref{no bad after good}  we show two restrictions  on the exceptional objects,  called RP property 1 and RP property 2, which ensure that   $\mc A$ is regularity-preserving.

\subsection{Ext-nontrivial couples}
Looking at the description of \textbf{B1, B2} (see \textbf{B1.4,
B2.3}) we see that in any of these cases occur couples $\{
L,\Gamma \} \subset \mc A$ of exceptional objects with
$\hom^1(L,\Gamma)\neq 0, \hom^1(\Gamma,L)\neq 0$. It is useful to
give a name to such a couple:
\begin{df} \label{Ext-nontrivial couple}An \uline{Ext-nontrivial couple} is a couple of exceptional objects $\{ L,\Gamma \} \subset {\mc A}_{exc}$, s. t.
$\hom^1(L,\Gamma)\neq 0$ and $\hom^1(\Gamma,L)\neq 0$.

\uline{Trivially coupling object} is an exceptional object $E\in
{\mc A}_{exc}$, s. t.  for each  $\Gamma\in {\mc A}_{exc}$
 we have $\hom^1(E,\Gamma)= 0$ or  $\hom^1(\Gamma,E) = 0$, i. e. for
each $\Gamma\in {\mc A}_{exc}$ the couple $\{E,\Gamma\}$ is not
Ext-nontrivial.
\end{df}
From \textbf{B1.4, B2.3} it follows
\begin{lemma} \label{from good to good case} If $E\in \mc A_{exc}$ is a trivially coupling object, then for each stability condition $\sigma \in \st(D^b({\mc A}))$ it is either
 $\sigma$-semistable or
$\sigma$-regular.
\end{lemma}
Thus, an object can be $\sigma$-irregular only if it is an element
of an Ext-nontrivial couple. The following lemma gives some
information about the other element of the  couple.
\begin{lemma}\label{no two bad cases} Let   each $X\in {\mc A}_{exc}$ satisfy the dichotomy  that it is
either trivially coupling or  there exists unique up to
isomorphism  another   object $Y\in {\mc A}_{exc}$ such that $\{
X,Y \}$ is an Ext-nontrivial couple.  Then for each Ext-nontrivial
couple $\{E,\Gamma\}\subset  {\mc A}_{exc}$ and each $\sigma \in
\st(D^b(\mc A))$ we have:
\begin{itemize}
 \item[\textbf{(a)}] If $E$ is  a \textbf{B2} object, then $\Gamma$ is semistable of phase $\phi_-(E)-1$.
 \item[\textbf{(b)}] If $E$ is  a \textbf{B1} object, then  $\phi_-(\Gamma)\geq \phi_-(E)+1$.
 \item[\textbf{(c)}] At most one of the objects $\{E,\Gamma\}$ can be  $\sigma$-irregular.
 \end{itemize}
\end{lemma}
\bpr

\textbf{(a)} By \textbf{B2.3} there exists a  semistable $X\in
{\mc A}_{exc}$ of phase $\phi_-(E)-1$, s. t. $\{ E,X\}$ is an
Ext-nontrivial couple. From the assumption of the lemma it follows
$X\cong \Gamma$.

\textbf{(b)} By \textbf{B1.3} and \textbf{B1.4} there exists $X\in
{\mc A}_{exc}$ with $\phi_-(X)\geq \phi_-(E)+1$, s. t. $\{ E,X\}$
is an Ext-nontrivial couple.
 From the assumption of the lemma we have $X\cong \Gamma$, hence $\phi_-(\Gamma)\geq \phi_-(E)+1$.

\textbf{(c)} It is enough to prove that if $E$ is
$\sigma$-irregular then $\Gamma$ is not $\sigma$-irregular. If $E$
is \textbf{B2}, then by \textbf{(a)}  $\Gamma$ is semistable, i.
e. it is not $\sigma$-irregular. By \textbf{(a)} applied to
$\Gamma$ it follows also that if $E$ is \textbf{B1} then $\Gamma$
is not \textbf{B2}. Whence, it remains to show that $E$ and
$\Gamma$ cannot both be \textbf{B1}. By \textbf{(b)} we see that
if both are \textbf{B1} then $\phi_-(\Gamma)\geq \phi_-(E)+1$ and
$\phi_-(E)\geq \phi_-(\Gamma)+1$ which is impossible. \epr The
next step is to show that even  with the presence of
Ext-nontrivial couples $\mc A$ could be regularity-preserving.
\subsection{RP property 1 and RP property 2} \label{RP properties} Our key to regularity-preserving of $\mc A$ are the following patterns of the Ext-nontrivial couples of $\mc A$.

\begin{df}
Let $\mc A$ be a hereditary category. We say that $\mc A$ has

   \textbf{RP Property 1}: if  for    each    Ext-nontrivial couple  $\{ \Gamma, \Gamma' \}\subset \mc A$  and    for each   $X \in {\mc A}_{exc}$ \\   from $\hom^*(\Gamma,X)=0 $ it follows $\hom^*(X,\Gamma')= 0$;

  \textbf{RP Property 2}: if for    each    Ext-nontrivial couple  $\{ \Gamma, \Gamma' \}\subset \mc A$  and for any two  $X,Y\in {\mc A}_{exc}$ \\
from  $ \hom(\Gamma,X)\neq 0, \hom(X,Y)\neq 0, \hom^*(\Gamma,Y) =
0$ it follows $ \hom(\Gamma',Y)\neq 0$.\footnote{note that $
\hom(\Gamma,X)\neq 0, \hom(X,Y)\neq 0, \hom^*(\Gamma,Y) = 0$ imply
$X \neq \Gamma$, $X \neq Y$}

\end{df}
The main result of Section \ref{no bad after good} is:
\begin{prop} \label{prop no bad after good} If $\mc A$ has RP Property 1 and  RP Property 2,\footnote{$\mc A$ is as in Subsection \ref{pesumptions A is HER}. } then  $\mc A$ is regularity-preserving.
\end{prop}
\subsection{Proof of Proposition \ref{prop no bad after good}}   We can assume that $R\in \mc A$. We split the proof in two lemmas.
The first lemma uses RP property 1, but does not use RP property 2.
\begin{lemma} Let $R$ be a \textbf{C3} object with
$ \mk{alg}(R)= \begin{diagram}[size=1em]
A & \rTo      &     &       &   R \\
  & \luDashto &     & \ldTo &       \\
  &           & B[1] &       &
\end{diagram} $. Then each non-semistable  $E \in Ind(A)$ is  $\sigma$-regular.
\end{lemma}
\begin{proof} Recall that in \textbf{C3.1, C3.2} we have $A,B\neq 0$, $\hom^*(A, B)=\hom^1(A,A)=\hom^1(B,B)=0$, and $\phi_-(A) > \phi(B)+1$. The last inequality,  together with   Corollary
 \ref{ineq for phi_-} and Remark \ref{remark for phi}, implies:
\be \label{proof of no bad after C3 1} \phi_-(E) > \phi(B)+1\ \
\Rightarrow \ \  \hom^*(E, B)=0. \ee If $E$ is a \textbf{B1}
object, then we get $ \mk{alg}(E)=
\begin{diagram}[1em]
A_1 \oplus A_2[-1] & \rTo      &     &       &   E, \\
  & \luDashto &     & \ldTo &       \\
  &           & B' &       &
\end{diagram} $
where $B'\in \mc A$ is a direct summand of $\sigma_-(E)$ (see \textbf{B1.2}). By \eqref{proof of no bad after C3 1}
we can apply Lemma \ref{lemma for hom leq 1(X,S)} to $E, B$ and
obtain $\hom^{\leq 1}(B',B)=0$, hence $\hom^{*}(B',B)=0$. From the
triangle
 $\mk{alg}(E)$ it follows $\hom^*(A_2,B)=0$. By \textbf{B1.4} there exists   $E'\in Ind(A_2)$ s. t. $\{E, E'\}$
  is an Ext-nontrivial couple. So, we obtained $\hom^*(E',B)=0$. Since $\hom^1(B,B)=0$,   RP property 1 in subsection \ref{RP properties}
  implies $\hom^*(B,E)=0$, which contradicts \textbf{C3.3}.

If $E$ is \textbf{B2} object, then we get $ \mk{alg}(E)=
\begin{diagram}[size=1em]
A' & \rTo      &     &       &   E \\
  & \luDashto &     & \ldTo &       \\
  &           & B'[1] &       &
\end{diagram}   $,
where $B'[1]=\sigma_-(E)$ (see \textbf{B2.2}). By \eqref{proof of
no bad after C3 1} we can apply Lemma \ref{lemma for hom leq
1(X,S)} to $E, B[1]$ and obtain $\hom^{\leq 1}(B'[1],B[1])=0$,
hence  $\hom^{*}(B',B)=0$. By \textbf{B2.3} there exists  $E' \in
Ind(B')$, s. t. $\{E, E'\}$ is an Ext-nontrivial couple. So, we
obtained $\hom^*(E',B)=0$ which by RP property 1 implies
$\hom^*(B,E)=0$. This contradicts \textbf{C3.3}. \epr The second
lemma uses both RP property 1 and RP property 2.
\begin{lemma}
Let $R$, $E \in {\mc A}_{exc}$, $R \not \in \sigma^{ss}$, $E \not
\in \sigma^{ss}$. If $R$, $E$ fit into  any of the following two
situations:
\begin{itemize} \item[\textbf{(a)}] $R$ is a \textbf{C1} object, $ \mk{alg}(R)= \begin{diagram}[size=1em]
A & \rTo      &     &       &   R \\
  & \luDashto &     & \ldTo &       \\
  &           & B &       &
\end{diagram}   $,  $E\in Ind(A)$;
\item[\textbf{(b)}]  $R$ is a  \textbf{C2} object, $ \mk{alg}(R)=\begin{diagram}[size=1em]
A_1 \oplus A_2[-1] & \rTo      &     &       &   R \\
  & \luDashto &     & \ldTo &       \\
  &           & B &       &
\end{diagram} $,   $E \in Ind(A_1)$ or $E \in Ind(A_2)$;
\end{itemize}
then $E$ is $\sigma$-regular.
\end{lemma}
\bpr  The arguments for  $E\in Ind(A)$, $R$ is \textbf{C1} and
$E\in Ind(A_1)$, $R$ is \textbf{C2} are similar. We give them
first. Recall that in \textbf{C1.3} and \textbf{C2.3} we have  $\phi_-(A)\geq \phi(B)$ and  $\phi_-(A_1)\geq \phi(B)$, respectively.

By Corollary \ref{ineq for phi_-}, in \textbf{C1} case we have
$\phi_-(E)\geq \phi_-(A)\geq \phi(B)$, and in \textbf{C2} case we have
$\phi_-(E)\geq \phi_-(A_1)\geq
\phi(B)$. In both the  cases (see \textbf{C2.1}, \textbf{C1.1}) we
have
  $\hom^*(E,B)=0$. In both the cases we  have also $\hom(E,R) \neq 0$ (recall that in \textbf{C2} case $A_2$ is a subobject
 of $R$), so  we can write
\be \label{proof of no bad after C1,C2 1} \phi_-(E) \geq \phi(B),
\ \  \hom^*(E,B)=0, \ \  \hom(E,R) \neq 0 \qquad E, B \in \mc
A.\ee If we take any  $X\in Ind(B)$, then $\hom(R,X) \neq 0$ (this
is valid in all the five cases\footnote{by the last part of Lemma \ref{the last HN triangle} and since $X$ is a direct summand of $\sigma_-(E)$}). Since $R$ is $\sigma$-regular,
we have $X,E \in {\mc A}_{exc}$ and  combining with \eqref{proof of no bad after C1,C2 1} we can write: \be \label{proof
of no bad after C1,C2 4} \hom(E,R) \neq 0, \hom(R,X) \neq 0,
\hom^*(E,X)=0, \qquad  X,E, R \in {\mc A}_{exc}.\ee

If $E$ is a \textbf{B2} object, then $\mk{alg}(E)$ is of the form \be \label{proof of
no bad after C1,C2 2}\mk{alg}(E)=
\begin{diagram}[1em]
A' & \rTo      &     &       &   E \\
  & \luDashto &     & \ldTo &       \\
  &           & B'[1] &       &
\end{diagram}.   \ee From \eqref{proof of no bad after C1,C2 1} we see that  Lemma \ref{lemma for hom leq 1(X,S)} can be applied, which
implies $\hom(B',B)=0$. By \textbf{B2.3}, there exists $E'\in
Ind(B')$, s. t.  $\{E, E'\}$ is an Ext-nontrivial couple. Then by
\eqref{proof of no bad after C1,C2 4} and RP property 2  we obtain
$\hom(E',X) \neq 0$, which contradicts $\hom(B',B)=0$.

If $E$ is \textbf{B1} object, then $\mk{alg}(E)$ is of the form \be
\label{proof of no bad after C1,C2 5}\mk{alg}(E)=
\begin{diagram}[1em]
A_1' \oplus A_2'[-1] & \rTo      &     &       &   E \\
  & \luDashto &     & \ldTo &       \\
  &           & B' &       &
\end{diagram}   \ee
with $B'\in \mc A$ and for some  $E' \in Ind(A_2')$ the couple
$\{E, E'\}$ is  Ext-nontrivial.
 From \eqref{proof of no bad after C1,C2 1} and Lemma \ref{lemma for hom leq 1(X,S)} it follows $\hom^{\leq 1}(B',B)=0$, hence   $\hom^{*}(B',B)=0$,
  which combined with $\hom^{*}(E,B)=0$ and the triangle \eqref{proof of no bad after C1,C2 5}, implies $\hom^{*}(A_2',B)=0$.
   Whence, we obtain $\hom^*(E',B)=0$, which by RP property 1 and $\hom^1(B,B)=0$ implies $\hom^*(B,E)=0$. The last contradicts \textbf{C1.4, C2.4}.

Suppose now that we are in the  situation \textbf{(b)} and  $E \in
Ind(A_2)$ is a \textbf{B2} object. Then we again have \eqref{proof
of no bad after C1,C2 2} and  some $E'\in Ind(B')$, s. t.  $\{E,
E'\}$ is an Ext-nontrivial couple. However, now in addition to
$\hom^*(E,B)=0$ we have $\phi_-(E) \geq
\phi_-(A_2)=\phi_-(A_2[-1])+1 \geq \phi(B)+1 = \phi(B[1])$. Now
Lemma \ref{lemma for hom leq 1(X,S)} gives $\hom^{\leq
1}(B'[1],B[1])=0$, i. e. $\hom^*(B',B)=0$. Thus, we obtain
$\hom^*(E',B)=0$, hence $\hom^*(B, E)=0$ by RP property 1, which
contradicts  \textbf{ C2.4}.

Finally, suppose that $E\in Ind(A_2)$ is a \textbf{B1} object.
Then we can use again \eqref{proof of no bad after C1,C2 5} and
take some $E'\in Ind(A_2')$, s. t.  $\{E, E'\}$ is an
Ext-nontrivial couple. As in the  preceding paragraph, in addition
to $\hom^*(E,B)=0$, we have again $\phi_-(E) \geq \phi(B[1])$. Now
Lemma \ref{lemma for hom leq 1(X,S)} gives $\hom^{\leq
1}(B',B[1])=0$, i. e. $\hom^*(B',B)=0$. Combining with
$\hom^{*}(E,B)=0$ and the triangle \eqref{proof of no bad after
C1,C2 5}  we obtain $\hom^*(E',B)=0$. As in the previous
paragraph, the last vanishing gives a contradiction. \epr

\section{Sequence of regular cases} \label{sequence}
In this section we assume that $\mc A$  is regularity-preserving.
If we are given  a  non-final $\sigma$-regular object $R$, then we can apply $\mk{alg}$ iteratively(Definition \ref{RP category}). As a result  we obtain a sequences
of exceptional pairs(between the subsequent  iterations we make a
choice, whence the resulting sequence is not uniquely determined
by $R$ in general):
\begin{gather} \label{sequence of cases}  \bd[height=1.5em] R & \rDotsto^{X_1} & (S_1,E_1) & \rMapsto^{proj_2}& E_1 & \rDotsto^{X_2} & (S_2,E_2)  & \rMapsto^{proj_2}& E_2 & \rDotsto^{X_3} & (S_3,E_3)& \rMapsto^{proj_2}&  \dots \\
  &   & \dMapsto^{proj_1} &    &  &    & \dMapsto^{proj_1} &  &  &  & \dMapsto^{proj_1} &  &  & \\
   &  &  S_1              &    &  &    &  S_2              &  &  &  & S_3               &  &  &  \ed \end{gather}
where $X_i \in \{ \textbf{C1, C2a, C2b,  C3} \} $.  Such a sequence will be called an \textit{$R$-sequence}.  The number of the objects $\{S_i\}$    will be called length of the
 $R$-sequence.\footnote{$R$ is the exceptional object, which is the origin of the sequence, so for example if the length is
 $\geq 2$, then after removing the first step $X_1$ we get an $E_1$-sequence.} We study here $R$-sequences.

The sequence \eqref{sequence of cases}
can be extended after $E_i$ iff  $E_i\not \in \sigma^{ss}$,
which is possible only if $E_{i-1}$ is not final (Definition
\ref{def final good}).  From \eqref{common properties for Dotsto 2}
it follows (recall that $\theta_\sigma(R)$ is an $\NN$-valued function with finite support)\be \label{theta dedreases}  \theta_\sigma(R) >
\theta_\sigma(E_1) >\theta_\sigma(E_2)>\dots. \ee
 Hence we see that \textit{after finitely many steps we   reach a final $\sigma$-regular object}.
 More precisely:
\begin{lemma} \label{R seq are finite} Let $R$ be  $\sigma$-regular. There does not exist an infinite $R$-sequence.
The lengths of all $R$-sequence are  bounded
above by $\sum_{u \in \sigma^{ss}_{ind}/\cong}
\theta_\sigma(R)(u)$.
\end{lemma}

 Some  features of the individual steps  in any $R$-sequence,
specified  in \eqref{common properties for Dotsto 1},
\eqref{common properties for Dotsto 2}, and   Lemma \ref{prop no bad
after good},     are readily integrated to the following basic
features of the whole $R$-sequence:

\begin{lemma} \label{lemma for R seq} Let $R$ be $\sigma$-regular. Let an $R$-sequence as \eqref{sequence of cases} have length $n$. Then
$\{ (S_i, E_i) \}_{i=1}^n $ is a sequence of exceptional pairs,
which, in addition to \eqref{theta dedreases}, satisfies the
following monotonicities:
\footnote{Recall that the notation $\deg(X)$ is explained in \textbf{Some notations} right after  the introduction.} \begin{gather}\label{phi(S_i) leq phi(S_i+1)}  \phi_-(R)=\phi(S_1) \leq \phi_-(E_1)=\phi(S_2)\leq \phi_-(E_2)=\phi(S_3)\leq\dots  \\
\label{deg(R) geq deg(E_1) geq dots}  \deg(R) \geq \deg(E_1) \geq
\deg(E_2) \geq \deg(E_3) \geq \dots \end{gather}
 where $\{S_i\}_{i=1}^n $ are semistable,  $\{ E_i \}_{i=1}^{n-1}$  are $\sigma$-regular, and
the last object $E_n$ is either semistable or again
$\sigma$-regular (and then the sequence can be extended).
\end{lemma}
In the rest of this section we  make various refinements of  Lemma
\ref{lemma for R seq}. Whence, \uline{ in the rest of this section
the objects $R$, $\{ (S_i, E_i) \}_{i=1}^n $, and  the integer
$n\in\NN$ will be as in Lemma \ref{lemma for R seq}, in particular
these objects fit in an $R$-sequence \eqref{sequence of cases},
which ends  at $E_n$.}  Assuming this data, we  will show that
under   additional conditions some of the inequalities in
\eqref{phi(S_i) leq phi(S_i+1)} are strict, and  vanishings, other
than  the already known $\{ \hom^*(E_i, S_i)=0 \}_{i=1}^n $,
appear.  The basic lemma is:
\begin{lemma} \label{lemma for R sequence and eqaul degrees}  Let $ 1\leq i <n$. Then  the following implications hold:
\begin{itemize}
    \item[\textbf{(a)}] If   $\deg(S_i)\geq \deg(S_{i+1})$, then $\hom^*(S_{i+1}, S_i)=0$.
    \item[\textbf{(b)}] If  $\deg(S_i)= \deg(S_{i+1})$, then  $\hom^*(S_{i+1}, S_i)=0$ and $\phi(S_{i+1})>\phi(S_{i})$.
        \item[\textbf{(c)}] If  $\deg(S_i)+1= \deg(S_{i+1})$, then $\hom^1(S_{i+1},S_{i})=0$.
\end{itemize}
 \end{lemma}
\bpr Since $E_i$ and $E_{i-1}$ are regular, all the four features
specified right after \eqref{the diagram for all cases} hold for
$\mk{alg}(E_{i-1})$ and  $\mk{alg}(E_{i})$.  Now we    unfold the
definitions and  use these features  to write:
\begin{gather} \mk{alg}(E_{i-1})= \begin{diagram}[size=1em]
U & \rTo      &     &       &   E_{i-1} \\
  & \luDashto &     & \ldTo &       \\
  &           & V   &       &
\end{diagram}   \ \  \mk{alg}(E_{i})=  \begin{diagram}[size=1em]
U' & \rTo      &     &       &   E_{i} \\
   & \luDashto &     & \ldTo &       \\
   &           & V' &       &
\end{diagram} \qquad \begin{array}{c} S_i \in Ind(V)\\E_i \in Ind(U)\\ S_{i+1}\in Ind(V') \end{array}  \begin{array}{c} \deg(S_i)=\deg(V)\\ \deg(S_{i+1})=\deg(V')\\ \phi(S_i)=\phi(V) \\ \phi(S_{i+1})=\phi(V') \end{array} \nonumber \end{gather}
 \be \label{proor of phi < phi 2} \hom^*(E_i,V)=0, \qquad  \phi(V')=\phi_-(E_i)\geq \phi(V), \qquad \theta_\sigma(E_{i})<\theta_\sigma(E_{i-1}).  \ee
The first two expressions in \eqref{proor of phi < phi 2} show
that we can apply Lemma  \ref{lemma for hom leq 1(X,S)} to $E_{i}$
and $V$. Since $V'$ is a direct summand of $\sigma_-(E_i)$ and  $\deg(S_{i+1})=\deg(V'), \deg(V)=\deg(S_i) $,  this
lemma  gives us:  $\hom^*(V',V)=0$, if $\deg(S_{i+1})\leq
\deg(S_i) $;  \ \ \ $\hom(V',V[1])=0$, if $\deg(S_{i+1})=\deg(S_i) +1$.

So far we proved \textbf{(a)}, \textbf{(c)}.  It remains to
 show  that the inequality  $\phi(S_{i+1}) \geq \phi(S_{i})$ given by \eqref{phi(S_i) leq phi(S_i+1)}  is  strict inequality $\phi(S_{i+1})>\phi(S_i)$
in \textbf{(b)}.  We first observe the following implication:
 \begin{gather} \label{proor of phi < phi 5} \phi(S_{i+1})=\phi(S_i) \qquad  \Rightarrow  \qquad  S_{i+1} \in Ind(\sigma_-(E_{i-1}))\cap Ind(\sigma_-(E_{i})). \end{gather}
  Indeed, by \eqref{common properties for Dotsto 2} we have $\theta_\sigma(E_{i})(S_{i+1})\neq 0$ . From \eqref{proor of phi < phi 2} it follows that $\theta_\sigma(E_{i-1})(S_{i+1})\neq 0$, hence $S_{i+1}$ is an indecomposable
  component of some HN factor of $E_{i-1}$. This must be $\sigma_-(E_{i-1})$, because the assumption $\phi(S_{i+1})=\phi(S_i)$ implies   $\phi_-(E_{i-1})=\phi(S_{i+1})$, so we obtain \eqref{proor of phi < phi 5}.

  Suppose that $\phi(S_i)=\phi(S_{i+1})$ and $\deg(S_i)=\deg(S_{i+1})$, then $\phi(V) = \phi(V')$ and $\deg(V)=\deg(V')=j$ for some $j \in \ZZ$.  Hence
       $V$ and $V'$ are the degree   $j$  terms of  $\sigma_-(E_{i-1})$ and $\sigma_-(E_i)$, respectively.
 Now   \eqref{proor of phi < phi 5} and Krull-Schmidt property imply
   $S_{i+1}\in Ind(V) \cap Ind(V')$, which contradicts the already proven  $\hom^*(V',V)=0$. Hence \textbf{(b)} and the lemma follow.
 \epr
\begin{coro} \label{coro for R seq initial C3,...} If for each $i\in \{ 1,2,\dots, n\}$ we have $\deg(S_1)\geq \deg(S_i)$, then:
\begin{itemize}
    \item[\textbf{(a)}] the vanishings $\hom^*(S_i,S_1)=\hom^*(E_i,S_1)=0$ hold for each integer $i$ with $ 2\leq i \leq n $,
\item[\textbf{(b)}] furthermore, if $\deg(S_i)=\deg(S_1)$ for some $i\geq 2$ then $\phi(S_1)<\phi(S_i)$.
\end{itemize}

The inequalities $\{ \deg(S_1)\geq \deg(S_i) \}_{i=1}^n$  hold in
any of the following cases:
\begin{itemize}
    \item $X_1=\textbf{C2a}$
    \item $X_1=\textbf{C3}$
    \item \textbf{C3} does not occur  in the sequence $\{X_1, X_2, X_3,\dots, X_n\}$.
\end{itemize}

\end{coro}
\bpr From Lemma \ref{lemma for R seq}    we have $\{\phi_-(E_i)
\geq \phi(S_1), \phi(S_i) \geq \phi(S_1)\}_{i=1}^n$ and
$\hom^*(E_1,S_1)=0$.

 Suppose that  for some $i$ with $1\leq i <n$  we are given  $\hom^*(E_{i},S_1)=0$ (here we make an induction assumption).
 We use the    triangle $\mk{alg}( E_i)$(it must be  of type  \textbf{C1, C2, C3}): \ben \mk{alg}( E_i)=\begin{diagram}[size=1em]
U & \rTo      &     &       &   E_i \\
  & \luDashto &     & \ldTo &       \\
  &           & V &       &
\end{diagram} \qquad \ \ \ U, V \in \mc T, U \neq 0, V \neq 0, \qquad \begin{array}{c}S_{i+1}\in Ind(V) \\ E_{i+1}\in Ind(U) \end{array} \een
where $V$ is a direct summand of  $\sigma_-(E_i)$
and $V$ is of pure degree.

By $\hom^*(E_{i},S_1)=0$,  $\phi_-(E_i) \geq \phi(S_1)$ we can
apply Lemma \ref{lemma for hom leq 1(X,S)} and we obtain \be
\hom^{\leq 1}(V,S_1)=0.\ee Therefore, if $\deg(S_{i+1})\leq
\deg(S_1)$, then $\hom^*(V,S_1)=0$, since   $
\deg(V)=\deg(S_{i+1})$. Now $\hom^{*}(V,S_1)=0$ together with the
induction assumption $\hom^*(E_{i},S_1)=0$ and the triangle
$\mk{alg}( E_i)$ give $\hom^*(U,S_1)=0$.  Hence
$\hom^*(E_{i+1},S_1)=0$ and $\hom^*(S_{i+1},S_1)=0$. Part
\textbf{(a)} follows.

We prove part \textbf{(b)} by contradiction.  Suppose that
 $\deg(S_i)=\deg(S_1)$ and $\phi(S_i)=\phi(S_1)$.   From \eqref{theta dedreases} and  \eqref{common properties for Dotsto 2} it follows
 $\theta_{\sigma}(R)(S_i) >\theta_{\sigma}(E_{i-1})(S_i) >0$, therefore $S_i$ is a direct summand of some HN factor of $R$. On the other
 hand by $\phi(S_1)=\phi_-(R)$,
   $\phi(S_i)=\phi(S_1)$, and  $\deg(S_i)=\deg(S_1)$  it follows $S_1, S_i \in Ind(V)$, where $V$ is the  degree $\deg(S_i)=\deg(S_1)$ term of $\sigma_-(R)$.    Therefore (recall also \textbf{C1.2, C2.2, C3.2}), we can write $\mk{alg}(R)=\begin{diagram}[size=1em]
U & \rTo      &     &       &   R \\
  & \luDashto &     & \ldTo &       \\
  &           & V &       &
\end{diagram}$ and $S_i \in Ind(V)$. The definition of $\bd &\rDotsto &\ed$ (Definition \ref{def of Dotsto})  implies that we can replace $S_1$ by $S_i$ in the
 $R$-sequence which we consider. However now   part \textbf{(a)} of the corollary says that $\hom^*(S_i,S_i)=0$, which  contradicts $S_i\neq 0$. Hence  $\phi(S_i)>\phi(S_1)$, if $\deg(S_i)=\deg(S_1)$ and part \textbf{(b)} is shown.

To prove the rest of the corollary, we use table \eqref{table with
degrees} for  comparing degrees.

If we are given $X_1=\textbf{C2a}$ or $X_1=\textbf{C3}$, then  $
\deg(E_1)=\deg(S_1)-1 $.   From \eqref{deg(R) geq deg(E_1) geq
dots} in Lemma \ref{lemma for R seq} we can write  that  $\deg(E_i) \leq \deg(E_1)=\deg(S_1)-1$ for
$i=1,2,\dots, n-1$, hence $\deg(E_i)+1\leq \deg(S_1)$. By  $\bd E_i & \rDotsto &(S_{i+1},E_{i+1})\ed$ and the last expression in \eqref{common properties for Dotsto 1} we have also $\deg(S_{i+1})\leq \deg(E_i)+1$. Hence, we obtain $\deg(S_{i+1}) \leq \deg(S_{1})$
 for $i=1,2,\dots, n-1$.

Finally, assume that  the sequence
$\{X_1, X_2, X_3,\dots, X_n\}$ does not contain  \textbf{C3}. By the already proven, we can
assume that $X_1=\textbf{C2b}$ or $X_1=\textbf{C1}$, which
implies
 $ \deg(E_1)=\deg(S_1)$.  Since  \textbf{C3} is forbidden, it follows   $\{ \deg(S_{i+1})=\deg(E_i) \}_{i=1}^{n-1}$, hence by \eqref{deg(R) geq deg(E_1) geq dots} we obtain
 $\{ \deg(S_{i+1}) \leq \deg(S_{1})\}_{i=1}^{n-1}$. The corollary is completely proved.
\epr
 Corollary \ref{coro for R seq initial C3,...} does not ensure the  vanishings $\{ \hom^*(S_i,S_1)=\hom^*(E_i,S_1)=0 \}_{i\geq 2}$ for $R$-sequences  with first step \textbf{C1} or \textbf{C2b} and containing a \textbf{C3} step.     The obstacle to obtain these vanishings for each $R$-sequence is that the data $\hom^*(X,S)=0$, $S\in\sigma^{ss}$, $\phi_-(X)\geq \phi(S)$ gives    $\hom^{\leq 1}(\sigma_-(X),S)=0$, but not  $\hom^*(\sigma_-(X),S)=0$ (see Lemma \ref{lemma for hom leq 1(X,S)}).

For certain $R$-sequences starting with a \textbf{C1} step and
ending with a \textbf{C3} step  we obtain these vanishings in the
next lemma,   but here we use the  property in Corollary \ref{RP
property 1,2 and.. for Q1} \textbf{(b)} for the first time.

\begin{lemma} \label{lemma for R-se starting with C1} Assume  that, besides being regularity-preserving,  the category $\mc A$ satisfies the following:
for any two $X, Y \in {\mc A}_{exc}$  at most one  degree in $\{
\hom^p(X, Y)\}_{p\in\ZZ}$ is nonzero.

If an $R$-sequence (as in Lemma \ref{lemma for R seq}) obeys the
following restrictions (all the three):\\
 $X_1= \textbf{C1}$; in the sequence $\{ X_2, X_3,\dots, X_{n-1}\}$ do
 not occur  \textbf{C2a} and \textbf{C3};  $X_n= \textbf{C3}$,\\
   then it satisfies  $ \hom^*(S_i,S_1)=\hom^*(E_i,S_1)=0 $ for $i=2,\dots,n$.
\end{lemma}
\bpr Applying  the previous lemma to the  sequence obtained  by
truncating  the last step $X_n$,  we obtain  the given vanishings
for $i<n$. We have to prove only  $
\hom^*(S_n,S_1)=\hom^*(E_n,S_1)=0 $.

We first observe that
from  $\begin{diagram} B & \rDotsto^{X} & (S,E) \end{diagram}$, $X\in \{ \textbf{C2b}, \textbf{C1} \}$ it
follows by  Definition \ref{def of Dotsto}  that $\deg(B)=\deg(E)$
and there exists a monic $E \rightarrow B$ in $\mc
A[\deg(B)]$. Therefore we can assume  that    $0=\deg(R)=\deg(E_1)=\dots =
\deg(E_{n-1})$ and  $E_1, E_2,\dots,E_{n-1}$ are $\mc
A$-subobjects of $R$.
 Since $X_n=\textbf{C3}$,  we have, by \textbf{C3.2}, that  $\mk{alg}(E_{n-1})\cong {\rm H\! N}_-(E_{n-1})$, and we can write:
 \begin{gather} \label{triangle C3 for proof}  \mk{alg}(E_{n-1})= \begin{diagram}[1em]
A & \rTo      &     &       &   E_{n-1} \\
  & \luDashto &     & \ldTo &       \\
  &           & B[1]&       &
\end{diagram}, \qquad\begin{array}{c} A,B \in {\mc {A}} \\ \phi_-(A) > \phi(B[1])=\phi_-(E_{n-1}) \end{array} \ \begin{array}{c} E_n\in Ind(A) \\ S_n \in Ind(B[1]). \end{array} \end{gather}
Let us take now any $\Gamma \in Ind(A)$.  From Lemma \ref{the
last HN triangle} we have  $\hom (\Gamma, E_{n-1}) \neq 0$.
Since $E_{n-1}$ is an  $\mc A$-subobject of $R$ and  $\Gamma \in
\mc A$,  it follows  that $\hom(\Gamma, R) \neq 0$.  By the given property of $\mc A$  it follows
 that $\hom^1(\Gamma,R)=0$ (any $\Gamma \in Ind(A)$ is exceptional object).  Therefore we obtain  $
\hom^1(A,R)=0. $  Since  $X_1=\textbf{C1}$, we have a diagram $
\mk{alg}(R)= \begin{diagram}[1em]
A'& \rTo      &     &       &   R, \\
  & \luDashto &     & \ldTo &       \\
  &           & B' &       &
\end{diagram} $
and $S_1 \in Ind(B')$, $E_1 \in Ind(A')$. By Lemma \ref{lemmaHER}
\textbf{(b)}
 it follows $\hom^1(A,B')=0$ . We have also
 $\phi_-(A) > \phi_-(E_{n-1})\geq \phi_-(R)=\phi(B')$,  therefore  $\hom(A,B')=0$. Thus,  we obtain $ \hom^*(A,B')=0$, and hence $ \hom^*(A,S_1)=0$. The triangle \eqref{triangle C3 for proof} and  $ \hom^*(E_{n-1},S_1)=0$ imply  $\hom^*(B,S_1)=0$.  The lemma follows.
\epr

\section{Final regular cases} \label{final}
Let $R$ be a final $\sigma$-regular object  and  $(S,E)$ be any
exceptional pair satisfying  $\bd R & \rDotsto^X & (S,E)\ed$,  $X
\in \{\textbf{C1, C2a, C2b, C3}\}$. We have that    $E \in
\sigma^{ss}$ from the very  definition of final(Definition
\ref{def final good}).  We show here that,  besides being
semistable, the exceptional pair $(S,E)$ satisfies
$\phi(S)<\phi(E)$(Corollary \ref{coro for final good case}).
Furthermore, if $R$ is the middle term  of an exceptional triple
$(S_{min},R, S_{max})$(see Corollary \ref{between min max 1}),
then  the quadruple $(S_{min},S,E, S_{max})$ is also exceptional.

All results here, except  the second part of Corollary \ref{sigma pairs}, hold without   regularity-preserving.

The first  lemma ensures some strict inequalities. In this respect
it is similar to Lemma \ref{lemma for R sequence and eqaul
degrees} \textbf{(b)} and Corollary \ref{coro for R seq initial
C3,...} \textbf{(b)}. As in their proofs, the function
$\theta_\sigma$ will be useful again here.

\begin{lemma} \label{lemma for final good case objects}
Let $R$ be a $\sigma$-regular object with $  \mk{alg}(R)=
\begin{diagram}[1em]
U & \rTo      &     &       &   R \\
  & \luDashto &     & \ldTo &       \\
  &           & V &       &
\end{diagram}$.  For each $\Gamma \in Ind(U)   $
 from  $\Gamma \in \sigma^{ss}$ it follows that $\phi(V)<\phi(\Gamma)$.     In particular, if $R$ is a final, then $\phi_-(U) > \phi(V)$.
\end{lemma}
\bpr For simplicity, let $R \in \mc A$. If $R$ is  a \textbf{C3}
object, then the lemma is true by  \textbf{C3.2}, so we can assume
that $R$ is  a \textbf{C1} or a \textbf{C2} object. Then the
triangle $\mk{alg}(R)$ is of the form (if $R$ is \textbf{C1}, then
$A_2=0$, otherwise $A_2 \neq 0$)
 \be \label{C1 and C2 united} \begin{diagram}[1em]
A_1 \oplus A_2[-1] & \rTo      &     &       &   R \\
  & \luDashto &     & \ldTo &       \\
  &           & B &       &
\end{diagram}   \qquad \qquad \begin{array}{c} \hom^*(A_1,B)= \hom^*(A_2,B)=0 \\  A_1, A_2, B \in \mc A \\ \theta_\sigma(A_1 \oplus A_2[-1]) < \theta_\sigma(R).\end{array} \ee
We consider first the case $\Gamma \in \sigma^{ss} \cap Ind(A_1)$.
Then $\theta_\sigma(\Gamma)\leq \theta_\sigma(A_1 \oplus A_2[-1])
< \theta_\sigma(R)$. Since $\Gamma$ is semistable, the last
inequality implies $\theta_\sigma(R)(\Gamma) \neq 0$, hence
$\Gamma$ is an indecomposable component of some HN factor of $R$.
If $\phi(\Gamma)=\phi(B)$ then this must be the minimal HN factor $\sigma_-(R)$.
On the other hand  $ \deg(\Gamma)=0$ and $B$ is the zero degree of
$\sigma_-(R)$. Therefore,  we see that if
$\phi(\Gamma)=\phi(B)$, then $\Gamma \in Ind(B)$, which contradicts
$\hom^*(A_1,B)=0$.

Now let $\Gamma \in \sigma^{ss}\cap Ind(A_2)$. Then
$\theta_\sigma(\Gamma[-1])\leq \theta_\sigma(A_1 \oplus A_2[-1]) <
\theta_\sigma(R)$ and as in the previous case we deduce that
$\Gamma[-1]$ is an indecomposable component of an HN factor of
$R$. If $\phi(\Gamma[-1])=\phi(B)$ then this must be
$\sigma_-(R)$, but $\deg(\Gamma[-1])=-1$, which contradicts Lemma
\ref{the last HN triangle} \textbf{(a)}.

If $R$ is final, then each $\Gamma \in  Ind(U)$ is semistable and
the lemma follows. \epr

By this lemma and Definition
\ref{def of Dotsto} we obtain:
\begin{coro} \label{coro for final good case} Let $R$ be  final $\sigma$-regular. Let $ \bd R & \rDotsto^{X} & (S,E) \ed $. Then $S,E \in \sigma^{ss} $  and $\phi(E)> \phi(S)$.
\end{coro}
Having    $\phi(S)< \phi(E)$, it follows that $(S,E[-i])$ is a   $\sigma$-pair(Definition \ref{sigma exceptional collection})  for some $i\geq 1$. Indeed, we have  $\phi(S)-1 <
\phi(E[-i])\leq \phi(S)$ for some $i\geq 1$. Since $\deg(S)\geq \deg(E)$(recall  \eqref{common properties
for Dotsto 1}), the pair
$(S,E[-i])$ has all the features of a $\sigma$-pair.  Thus, we obtain the first part of the following corollary:
\begin{coro} \label{sigma pairs}  Each  final $\sigma$-regular object  implies the existence of a $\sigma$-exceptional pair.

In particular, if $\mc A$ is regularity-preserving, then   each  $\sigma$-regular object induces such a  pair.
\end{coro}
\bpr
If there exists a $\sigma$-regular object, then by preserving of regularity and Lemma \ref{R
seq are finite} we get a final $\sigma$-regular
object. Hence, by the first part, we obtain a $\sigma$-exceptional pair. \epr
If $\mc A$ has not  Ext-nontrivial couples, then each non-semistable exceptional object is  $\sigma$-regular  for  each stability condition, hence:
\begin{remark} \label{if there are not Ext-nontrivila pairs} If there are not Ext-nontrivial couples in $\mc A$, as in $\mc A=Rep_k(K(l))$, then each non-semistable  exceptional object induces a $\sigma$-exceptional pair.
\end{remark}
The origin of  our main $\sigma$-triples criterion(Proposition \ref{between min and max}) is in the next corollary.
\begin{coro} \label{between min max 1} If we are given the following data:
\begin{itemize}
    \item  $S_{min}, S_{max} \in \sigma^{ss}\cap  {\mc A}_{exc}$ with   $\phi(S_{min})\leq \phi(A)\leq \phi(S_{max})$ for each $A \in \sigma^{ss}\cap {\mc A}_{exc}$
    \item $(S_{min}, R, S_{max})$ is an exceptional  triple, s. t. $R\in {\mc A}_{exc}$  is   final and $\sigma$-regular
    \item $ \bd R & \rDotsto^{X} & (S,E) \ed $, $X \in \textbf{C1,C2a,C2b,C3}$,
\end{itemize}
then $(S_{min}, S,E, S_{max})$  is a semistable exceptional
quadruple (and no two of ${R, S, E}$ are isomorphic).
\end{coro}
\bpr  We have $\hom^*(E,S)=0$ (in particular $S \not \cong  E$)
and we  must show that
$\hom^*(S_{max},S)=\hom^*(S_{max},E)=\hom^*(S,S_{min})=\hom^*(E,S_{min})=0$.
By assumption $R$ is final and then both  $S, E$ are semistable.
Since $R$ is not semistable, it cannot be isomorphic  to $S$
or to $E$.

 Let us assume first that $R$ is a  \textbf{C3} object. Then we have a triangle $ \mk{alg}(R)=
\begin{diagram}[1em]
A & \rTo      &     &       &   R \\
  & \luDashto &     & \ldTo &       \\
  &           & B[1] &       &
\end{diagram} $ with $ \hom^*(A,B)=0 $
and $E\in Ind(A) $, $S\in Ind(B[1])$. The assumptions on $S_{min}$,$S_{max}$ and \textbf{C3.2}  imply
$$ \phi(S_{max}) \geq \phi(E) > \phi(B)+1=\phi(S)+1 \geq \phi(S_{min})+1.$$
Hence $\hom^*(S_{max},B)=0$, which, combined with  $\hom^*(S_{max}, R)=0$ and the triangle $\mk{alg}(R)$,
implies $\hom^*(S_{max},A)=0$. Thus, we get $\hom^*(S_{max},S)=\hom^*(S_{max},E)=0$.  Since
 each  $\Gamma \in Ind(A)$ satisfies  $\phi(\Gamma) > \phi(B)+1 \geq \phi(S_{min})+1$, we have  $\hom^*(A,S_{min})=0$. Now
   $\hom^*(R, S_{min})=0$ and $\mk{alg}(R)$ imply $\hom^*(B,S_{min})=0$. Thus, we get $\hom^*(S,S_{min})=\hom^*(E,S_{min})=0$ as well.

  Let us assume now that $R$ is a \textbf{C1} or \textbf{C2} object.  Then the triangle $ \mk{alg}(R)$ is of the  form (if $R$ is \textbf{C1}, then $A_2=0$,
   otherwise $A_2 \neq 0$):
\be \label{C1 and C2 united 1} \mk{alg}(R)=\begin{diagram}[1em]
A_1 \oplus A_2[-1] & \rTo      &     &       &   R \\
  & \luDashto &     & \ldTo &       \\
  &           & B &       &
\end{diagram}   \ \ \ \qquad \begin{array}{c}A_1,A_2,B \in \mc A \\ \hom^*(A_1,B)= \hom^*(A_2,B)=0 \\ E\in Ind(A_1 \oplus A_2[-1]),\ S\in Ind(B).\end{array}\ee
  Since $ B\neq 0 $ is semistable and $\hom^1(B,B)=0$, it follows
$\phi(S_{max}) \geq \phi(B) \geq \phi(S_{min})$. On the other hand we have
$\hom^*(R,S_{min})=0$ and $\phi_-(R)= \phi(B)$.
From Lemma \ref{lemma for hom leq 1(X,S)} it follows
$\hom^*(B,S_{min})=0$, which, combined with $\hom^*(R,S_{min})=0$
and the triangle $\mk{alg}(R)$, implies $\hom^*(A_1 \oplus
A_2[-1],S_{min})=0$.  So,  we obtained  $\hom^*(S,S_{min})=\hom^*(E,S_{min})=0$ and it remains to show $\hom^*(S_{max},S)=\hom^*(S_{max},E)=0$.    From Lemma \ref{lemma for final good case
objects} it follows that  for each indecomposable component $\Gamma$ of $A_1$,
resp $A_2$, we have $\phi(\Gamma) > \phi(B)$, resp.
$\phi(\Gamma[-1]) > \phi(B)$, and combining with $\phi(S_{max})
\geq \phi(\Gamma)$ we see that $\phi(S_{max}) > \phi(B)$, hence
$ \hom(S_{max},B)=0.$

Furthermore, if $R$ is \textbf{C2}, then $A_2\neq 0$ and  $\phi(S_{max})
\geq \phi(\Gamma)$, $\phi(\Gamma[-1]) > \phi(B)$  for each $\Gamma \in Ind(A_2)$. Therefore $\phi(S_{max}) >
\phi(B)+1$ and $\hom^*(S_{max},B)=0$. The latter together with
$\hom^*(S_{max},R)=0$ imply $\hom^*(S_{max},A_1 \oplus A_2[-1])=0$,
and   the corollary follows.

Finally, if $R$ is \textbf{C1}, then
$A_2=0$ in the triangle \eqref{C1 and C2 united 1} and we have a short exact sequence
$0 \rightarrow  A_1 \rightarrow R \rightarrow B \rightarrow 0 $.
Hence, by Lemma \ref{lemmaHER} and $\hom(S_{max},R[1])=0$ we get
$\hom(S_{max},B[1])=0$. We showed already that
$\hom(S_{max},B)=0$, therefore $\hom^*(S_{max},B)=0$. Using again  the triangle \eqref{C1 and C2 united 1}  and
$\hom^*(S_{max},R)=0$ we obtain $\hom^*(S_{max},A_1 )=0$.
 The corollary follows.
\epr

\section{Constructing \texorpdfstring{$\sigma$}{\space}-exceptional triples}     \label{constructing}   \mbox{}
So far, the property of Corollary \ref{RP property 1,2 and.. for
Q1} \textbf{(b)} was used  only in Lemma  \ref{lemma for R-se
starting with C1}. In this section it is  used throughout.   We
start with a simple observation:
 \begin{lemma} \label{another prop of exc in quiver} Let $\mc A$ be as in Subsection \ref{pesumptions A is HER}. Let Corollary \ref{RP property 1,2 and.. for Q1} \textbf{(b)} hold for $\mc A$.  Then for any two    non-isomorphic exceptional
 objects $A,B \in  \mc A$  we have $\hom(A,B)=0$ or $\hom(B,A)=0$.

In particular, if $C \in \mc A$ satisfies $\hom^1(C,C)=0$, then
for any two non-isomorphic $A,B \in Ind(C)$ one of the pairs $(A,B)$,
$(B,A)$ is  exceptional.
 \end{lemma}
 \bpr  Let $\hom(A,B)\neq 0$. Take  a nonzero $u:A\rightarrow B$. By  Corollary \ref{RP property 1,2 and.. for Q1} \textbf{(b)} it follows  $\hom^1(A,B)= 0$. One can show that  \cite[Lemma 1, page 9]{WCB2} holds for $\mc A$, so  $\hom^1(A,B)= 0$ implies that  every nonzero $f\in \hom(B,A)$ is either monic or epic. Suppose that  $f\in \hom(B,A)$  is epic, then $u\circ f \in \hom(B,B)=k$ is nonzero,  hence $f$ is also monic. Therefore $f$ is invertible, which contradicts the assumptions.  If $f$ is monic, then we consider $f\circ u $ and again get a contradiction.

 The second part follows from Remark \ref{from pre-exceptional to exceptional}.
  \epr

\uline{Besides the restrictions   of Subsection \ref{pesumptions A
is HER},   we assume throughout   Section \ref{constructing} that
Corollary \ref{RP property 1,2 and.. for Q1} \textbf{(b)} holds
for $\mc A$ and that $\mc A$ is regularity-preserving.}
   In Subsection \ref{with the additional RP property}, besides  these features, we assume    that $\mc A$ has  the additional RP property (Corollary \ref{additional RP property}) and that Corollary \ref{coro for isom
triples} holds for it.   In particular, all results hold for $\mc
A = Rep_k(Q_1)$ (the preserving of regularity follows from
Corollary \ref{RP property 1,2 and.. for Q1} \textbf{(a)} and
Proposition \ref{prop no bad after good}).

 We denote  $D^b(\mc A)$ by $\mc T$,  and choose any $\sigma \in \st(\mc T)$.
  In Corollary \ref{sigma pairs} is shown that any $\sigma$-regular object  $R$ induces a $\sigma$-pair. If $R$ is   final, then  this pair is of the form $(S,E[-j])$ with $j\geq 0$, for any $\bd R&\rDotsto&(S,E)\ed$. Using a $\sigma$-regular object $R$,  we will obtain in this section  various  criteria  for existence of $\sigma$-exceptional triples in $\mc T$.
 To obtain a $\sigma$-triple we utilize  three   approaches:  using long $R$-sequences(of length greater than one); combining the $\sigma$-pairs induced by several single step $R$-sequences with a final $R$;  combining a $\sigma$-pair induced  from   $R$ with  a semistable $S\in \mc A_{exc}\cap \sigma^{ss}$ of phase close to the minimal/maximal phase.  The  minimal and maximal phases are defined by\footnote{For the notation $\sigma^{ss}$ see \eqref{sigma^{ss}} and recall that by ${\mc A}_{exc}$ we denote the set of exceptional
 objects of $\mc A$}
 \begin{gather}\label{phi_min,phi_max} \phi_{min}=\inf(\{\phi(S) : S \in \sigma^{ss} \cap {\mc A}_{exc}\})  \qquad \phi_{max}=\sup(\{\phi(S) : S \in \sigma^{ss} \cap {\mc A}_{exc}\}).\end{gather}

    Note  that if  Corollary \ref{coro for isom triples} holds for $\mc A$, which is assumed in  Subsection \ref{with the additional RP property}, then we have
$-\infty  < \phi_{min}\leq \phi_{max}<\infty$.
 Indeed, if some of the strict inequalities fails,
 then we can construct a sequence $S_1,S_2,S_3,S_4,\dots, S_n$ (as long as we want) of semistable exceptional objects in $\mc A$, s. t. $\{ \phi(S_i)+1 < \phi(S_{i+1}) \}_{i=1}^{n-1}$,  which contradicts Corollary \ref{coro for isom triples}.

 We denote by $S_{min}/S_{max}$ objects in $\mc A_{exc}\cap \sigma^{ss}$ satisfying   $\phi(S_{min})=\phi_{min}$/$\phi(S_{max})=\phi_{max}$,  this  can be expressed by writing  $S_{min/max}\in \mc P(\phi_{min/max})\cap \mc A_{exc}$.

We note in advance that by replacing ``\textbf{C3}'' with ``\textbf{C2}'' and ``$>\phi_{min}$'' with ``$<\phi_{max}$'' we obtain the  criteria in  which     $R$ is a  \textbf{C2} object  from those in which $R$ is  a \textbf{C3} object.
  However,     the proof of the \textbf{C2} versions  demands more efforts and  more assumptions on $\mc A$(the additional RP property and  Corollary \ref{coro for isom
triples}).

   The criteria using  long $R$-sequences with a   \textbf{C1} object $R$  are weaker  than those   with   \textbf{C2}/\textbf{C3}.

  The distinction between \textbf{C1, C2, C3} is not essential in Lemma \ref{lemma for final good case objects in main th for Q1 pr} (based on the second approach, where $R$ is final)  and in Proposition \ref{between min and max}.  Furthermore,  Proposition \ref{between min and max} asserts  that if $\phi_{min}-\phi_{max}>1$, then any  non-semistable $E \in \mc A_{exc}$, which is a middle term of an exceptional triple $(S_{min},E,S_{max})$ induces a $\sigma$-exceptional triple (the regularity of $E$ follows).

\subsection{Constructions without assuming the additional RP property} \label{without the additional RP} \mbox{}\\
 Recall(Definition \ref{def of sigma-exceptional collection}) that an exceptional triple $(S_0,S_1,S_2)$ is said to be $\sigma$-exceptional under  three conditions: it must be semistable, it must satisfy $\hom^{\leq 0}(S_0,S_1)=\hom^{\leq 0}(S_0,S_2)=\hom^{\leq 0}(S_1,S_2)=0$, and the phases of its elements must be in $ (t,t+1]$ for some $t\in \RR$.  If we are given only that $(S_0,S_1,S_2)$ is semistable, then we can always ensure  the second or the third condition by applying the shift functor to $S_1,$ $S_2$, but both  together - not always. For example
 if $\phi(S_i)=\phi(S_{i+1})$, $\hom(S_i,S_{i+1})\neq 0$
($i=0,1$),  then this cannot be achieved (similarly, if
$\phi(S_i)=\phi(S_{i+1})+1$, $\hom^1(S_i,S_{i+1})\neq 0$).
 In the following lemma are given  some cases  in which  this can be achieved.  We give the arguments for one of them. The rest are also easy. Keeping in mind Remark \ref{equiv of open interval} is useful, when checking these implications.
\begin{lemma} \label{inequalities} Let $(S_0,S_1, S_2)$ be a semistable exceptional triple,
where $S_0,S_1, S_2 \in \mc A$. If any of the following conditions
holds:
\begin{itemize}
    \item[\textbf{(a)}] $\phi(S_0)<\phi(S_1)<\phi(S_2)$, $1+\phi(S_0)<\phi(S_2)$
    \item[\textbf{(b)}] $\phi(S_0)\leq \phi(S_1) < \phi(S_2)$, $\hom(S_0,S_1)=0$
    \item[\textbf{(c)}] $\phi(S_0)< \phi(S_1) \leq \phi(S_2)$, $\hom(S_1,S_2)=0$
    \item[\textbf{(d)}] $\phi(S_0)<\phi(S_2)\leq \phi(S_1) < \phi(S_2)+1$, $\hom(S_1,S_2)=0$
    \item[\textbf{(e)}] $\phi(S_0)< \phi(S_1)+1$, $\phi(S_1) < \phi(S_2)$, $\phi(S_0)< \phi(S_2)$, $\hom(S_0,S_1)=0$
    \item[\textbf{(f)}] $\phi(S_0)< \phi(S_1)+1$, $\phi(S_0)< \phi(S_2)+1$, $\phi(S_1)< \phi(S_2)+1$, $\hom(S_0,S_1)=\hom(S_0,S_2)=\hom(S_1,S_2)=0$
    \item[\textbf{(g)}] $\phi(S_0)< \phi(S_2)$,$\phi(S_0)+1 < \phi(S_1)$,$\phi(S_2)\neq \phi(S_1[-1])$, $\hom(S_0,S_2)=\hom(S_1,S_2)=0$,
 \end{itemize}
then for some integers $0\leq i$, $0\leq  j$ the triple $(S_0,
S_1[-i], S_2[-j])$ is $\sigma$-exceptional.
\end{lemma}
\bpr \textbf{(d)} From $\phi(S_0)<\phi(S_2)$ it follows that
$\phi(S_2[-j])\leq \phi(S_0)<\phi(S_2[-j])+1$ for some $j\geq 1$.
From $\phi(S_2)\leq \phi(S_1) < \phi(S_2)+1$ it follows $
\phi(S_2[-j])\leq \phi(S_1[-j])<\phi(S_2[-j])+1$. Now
$\hom(S_1,S_2)=0$  implies  that $(S_0,S_1[-j],S_2[-j])$ is a
$\sigma$-exceptional triple.
 \epr

The next lemma is a  step in the proof of our basic long $R$-sequences criterion Proposition
\ref{non-existence of some sequeces}.
\begin{lemma} \label{lemma for C1, C2b} Let  $ \begin{diagram} R & \rDotsto^{\textbf{X}} & (S,E) \ed $, where $X \in \{ \textbf{C1, C2b}\}$. Then there exists $S'$, such that
$$ \begin{diagram} R & \rDotsto^{X} & (S',E) \end{diagram}, \  \hom(S',E) = 0, \ \hom(R,S')\neq 0, \ \hom(E,R)\neq 0. $$ \end{lemma}
\bpr By Definition \ref{def of Dotsto} with $X \in \{ \textbf{C1,
C2b}\}$, there is a triangle of the form\footnote{If
$X=\textbf{C1}$, then $A_2=0$.  If $X=\textbf{C2b}$, then $A_2\neq
0$.}  $ A_1\oplus A_2[-1] \rightarrow R \rightarrow B \rightarrow
A_1[1]\oplus A_2$   and $S \in Ind(B)$, $E \in Ind(A_1)$.
Furthermore,   any  $A'\in Ind(A_1)$, $B'\in Ind(B)$  satisfy
$\hom^1(B,A')\neq 0$, $\hom(R,B')\neq 0$, $\hom(A',R)\neq 0$(see
\textbf{C1, C2} and  Lemma \ref{the last HN triangle}).  In
particular, there exists  $S' \in Ind(B)$, with $\hom^1(S',E)\neq
0$. By Corollary \ref{RP property 1,2 and.. for Q1} (b) it follows
$\hom(S',E) = 0$. The lemma follows. \epr

Now we obtain $\sigma$-triples from    certain, but not all, long\footnote{by ``long'' we mean of length greater than one} $R$-sequences.

  \begin{prop} \label{non-existence of some sequeces} If there exists an  $R$-sequence
 \begin{gather} \label{sequence of cases 1}
\bd[size=1.5em]
R & \rDotsto^{X_1}& (S_1,E_1)         & \rMapsto^{proj_2}& E_1& \rDotsto^{X_2}& (S_2,E_2)         & \rMapsto^{proj_2}& E_2 & \rDotsto^{X_3}& \dots & \rMapsto^{proj_2}&E_{n-1}& \rDotsto^{X_n}&(S_n,E_n) & \rMapsto^{proj_2}& E_n \\
  &               & \dMapsto^{proj_1} &                  &    &               & \dMapsto^{proj_1}&                  &     &               &       &  & & & \dMapsto^{proj_1} &  &   \\
  &               &  S_1              &                  &    &               &  S_2             &                  &     &               & \dots & &  & & S_n      &  &    \ed  \end{gather}
 with $n\geq 2$, $E_{n-1}$ is final, and $\{ \deg(S_1)\geq \deg(S_i) \}_{i=1}^n$, then   there exists a $\sigma$-exceptional triple.

\end{prop}
\bpr Assume that such a sequence exists. Since $E_{n-1}$ is final,
Corollary \ref{coro for final good case} implies that $S_n$ and
$E_n$
  are both semistable and $\phi(E_n) > \phi(S_n)$. Since  $ \deg(S_1)\geq \deg(S_i)$ for each $i=\{1,2,\dots,n\}$, by Corollary \ref{coro for R seq initial C3,...}  and table
   \eqref{table with degrees} we obtain
 \begin{gather} \hom^*(S_n,S_1)=\hom^*(E_n,S_1)=0  \nonumber \\
 \deg(S_1) \geq \deg(S_n) \geq \deg(E_n), \qquad \phi(S_1)\leq  \phi(S_n) < \phi(E_n). \nonumber \end{gather}
   In particular, the exceptional triple $(S_1,S_n,E_n)$ is  semistable   and after shifting we obtain a triple of the form $(A,B[-i],C[-i-j])$ with
   $0\leq i, 0 \leq j$, $\phi(A)\leq \phi(B[-i])<\phi(C[-j-i])$, $A,B,C \in \mc A$. If $i\neq 0$, then
    Lemma \ref{inequalities}, \textbf{(a)} can be applied to the triple $(A,B,C)$ and the proposition follows.

   If  $i=0$, then  $\deg(S_1)=\deg(S_n)$. By  Corollary \ref{coro for R seq initial C3,...} (b)
   it follows $\phi(S_1)<\phi(S_n)$. Whence,  we obtain a semistable triple $(A,B,C[-j])$ with
   $ 0 \leq j$, $\phi(A)< \phi(B)<\phi(C[-j])$. If $j\neq 0$, then the triple
   $(A,B,C)$ satisfies the conditions in Lemma \ref{inequalities}, \textbf{(a)}.
   If   $j =  0$, then  $X_n\in \{\textbf{C2b,C1} \}$ and due to Lemma \ref{lemma for C1, C2b}  we can assume that $\hom(S_n,E_n) =\hom(B,C)=0$.
    Now the  triple $(A,B,C)$ satisfies the conditions in Lemma \ref{inequalities}, \textbf{(c)}. The proposition follows.
\epr It follows that any long $R$-sequence starting with a
\textbf{C3} or  a \textbf{C2a} step induces a $\sigma$-triple:
 \begin{coro} \label{all C3 are final Q1}  From the data: $ \bd R & \rDotsto^{X}& (S,E) \ed $, $X\in\{\textbf{C3, C2a}\}$, $E \not \in \sigma^{ss}$ it follows  that there  exists a $\sigma$-exceptional triple. In particular  each non-final \textbf{C3} object implies such a triple.
\end{coro}
\bpr  
Since $E \not \in \sigma^{ss}$,  by Lemma \ref{R seq are finite}
we obtain an $R$-sequence with maximal length $n\geq 2$ and with
first step  the given $ \bd R & \rDotsto^{X}& (S,E) \ed $. This
sequences is of the form \eqref{sequence of cases 1}  with $X_1=X$, $n\geq
2$.   As far as the sequence is  of maximal length, the object
$E_{n-1}$ must be  final and $\sigma$-regular.   Since $X_1=X \in
\{\textbf{C3, C2a}\}$, Corollary \ref{coro for R seq initial
C3,...} gives  $\{\deg(S_1)\geq \deg(S_i)\}_{i=1}^n$. Thus, we
constructed an $R$-sequence \eqref{sequence of cases 1} with   the
three properties used in  Proposition \ref{non-existence of some
sequeces}.  The corollary follows. \epr

The next lemma uses a final regular object $R$, so we do not have long $R$-sequences here.

 \begin{lemma} \label{lemma for final good case objects in main th for Q1 pr}
Let $R$ be a final $\sigma$-regular object
with $\mk{alg}(R)= \begin{diagram}[1em]
U & \rTo      &     &       &   R \\
  & \luDashto &     & \ldTo &       \\
  &           & V &       &
\end{diagram}.   $  Then we have:
\begin{itemize}
    \item[\textbf{(a)}] If $\mk{alg}(R)$ is not the HN filtration of $R$, then $U$ is not semistable.
    \item[\textbf{(b)}] If  $U$ is not semistable,   then there exists a $\sigma$-exceptional triple.
\end{itemize}
\end{lemma}
\bpr   Without loss of generality we can assume that  $R \in \mc
A$.
  Since $R$ is  a final $\sigma$-regular object,  any $\Gamma\in
Ind(U)$ is a semistable exceptional object, and hence  by Lemma
\ref{lemma for final good case objects} it satisfies $\phi(\Gamma)
> \phi(V)$. Now part  \textbf{(a)}  is clear and it remains to prove \textbf{(b)}.

If $U$ is not semistable, then there exists  a   pair of
non-isomorphic   $\Gamma_1,\Gamma_2 \in Ind(U)$ with different
phases.  We can assume  $\phi(\Gamma_2)>\phi(\Gamma_1)$. In
particular, for the rest of the proof we can use \be \label{ineq}
\hom(\Gamma_2,\Gamma_1)=0 \qquad
\phi(\Gamma_2)>\phi(\Gamma_1)>\phi(V). \ee

First, assume that $R$ is a \textbf{C1} object. Then the triangle
$\mk{alg}(R)$ and some of its properties are
\begin{gather} \mk{alg}(R)=\begin{diagram}[1em]
A & \rTo      &     &       &   R \\
  & \luDashto &     & \ldTo &       \\
  &           & B &       &
\end{diagram} \qquad  A,B \in {\mc A}, \hom^1(A,A)= \hom^1(B,B)=\hom^*(A,B)=0. \nonumber \end{gather}
By $\hom^1(A,A)=0$ we have $\hom^1(\Gamma_2,\Gamma_1)=0$, which,
combined with $\hom(\Gamma_2,\Gamma_1)=0$, implies
$\hom^*(\Gamma_2, \Gamma_1)=0$. By $\hom^*(A,B)=0$ it follows that
for each  $\Gamma \in Ind(B)$ we have $\hom^*(\Gamma_i,
\Gamma)=0$, $i=1,2$.  Hence for each  $\Gamma \in Ind(B)$ the
triple $(\Gamma, \Gamma_1, \Gamma_2)$ is  exceptional   and
$\phi(V)=\phi(\Gamma)<\phi(\Gamma_1)<\phi(\Gamma_2)$. By
\textbf{C1.4} we have $\hom^1(B,\Gamma_1)\neq 0$, and hence we can
choose $\Gamma$ so that $\hom^1(\Gamma,\Gamma_1)\neq 0$, which by
Corollary \ref{RP property 1,2 and.. for Q1} \textbf{(b)} implies
$\hom(\Gamma,\Gamma_1) = 0$. Thus, we constructed an exceptional
triple $(\Gamma, \Gamma_1, \Gamma_2)$ with $\hom(\Gamma,\Gamma_1)
= 0$, $\phi(\Gamma)<\phi(\Gamma_1)<\phi(\Gamma_2)$. By Lemma
\ref{inequalities} \textbf{(b)}, after shifting  this triple
becomes $\sigma$-exceptional.

In \textbf{C3} case:
\begin{gather} \mk{alg}(R)=\begin{diagram}[1em]
A & \rTo      &     &       &   R \\
  & \luDashto &     & \ldTo &       \\
  &           & B[1] &       &
\end{diagram} \qquad  A,B \in {\mc A}\setminus \{0\}, \hom^1(A,A)= \hom^1(B,B)=\hom^*(A,B)=0. \nonumber \end{gather}
As in the previous case we obtain that for each  $\Gamma \in
Ind(B)$ the triple $(\Gamma, \Gamma_1, \Gamma_2)$ is  exceptional.
Now \eqref{ineq} becomes
$\phi(V)=\phi(\Gamma)+1<\phi(\Gamma_1)<\phi(\Gamma_2)$ and Lemma
\ref{inequalities}, \textbf{(a)} gives a $\sigma$-triple.

In \textbf{C2} case the triangle $\mk{alg}(R)$ and some of its
properties are:
\begin{gather} \label{one C2 case triangle} \begin{diagram}[1em]
A_1 \oplus A_2[-1] & \rTo      &     &       &   R \\
  & \luDashto &     & \ldTo &       \\
  &           & B &       &
\end{diagram} \ \ \ \qquad  \begin{array}{c}A_2,B \in {\mc A}\setminus \{0\} \\
                  \hom^1(A_1,A_1)= \hom^1(A_2,A_2)=\hom^1(B,B)=0\\
                 \hom^*(A_1,A_2)= \hom^*(A_1,B)=\hom^*(A_2,B)=0 \end{array}.  \end{gather}

If both $\Gamma_1, \Gamma_2 \in Ind(A_1)$, then the arguments are
the same as in \textbf{C1} case.

If both $\Gamma_1, \Gamma_2\in Ind(A_2[-1])$, then
$\hom^1(B,B)=\hom^*(A_2,B)=0$ imply
  that for each  $\Gamma \in Ind(B)$ the triple $(\Gamma, \Gamma_1,
\Gamma_2)$ is  exceptional  and now $\Gamma_i[1] \in \mc A$,
$\phi(\Gamma_i[1]) > \phi(B)+1=\phi(\Gamma)+1$, i. e.
$\phi(\Gamma)+1 < \phi(\Gamma_1[1]) < \phi(\Gamma_2[1])$. From
this data Lemma \ref{inequalities} \textbf{(a)} produces a
$\sigma$-exceptional triple.

Before we continue with the other  possibility, we note that \be
\hom(A_2,A_1)=0.\ee Indeed, by \textbf{C2.4} for each  $\Gamma \in
Ind(A_2)$  we have $\hom(\Gamma, R[1])\neq 0$, then by Corollary
\ref{RP property 1,2 and.. for Q1} \textbf{(b)} it follows  $\hom(\Gamma,
R) = 0$, i. e. $\hom(A_2, R) = 0$. Now $\hom(A_2,A_1)=0$ follows
from the fact that $A_1$ is a proper subobject of $R$ in $\mc A$.

If $\Gamma_1\in Ind(A_1)$, $\Gamma_2 \in Ind(A_2[-1])$, then(see
\eqref{one C2 case triangle})  for each $\Gamma \in Ind(B)$  the
triple  $(\Gamma, \Gamma_2, \Gamma_1)$ is exceptional. We will
show that $\Gamma \in Ind(B)$ can be chosen so that the conditions
of Lemma \ref{inequalities} \textbf{(g)} hold with the triple
$(\Gamma, \Gamma_2[1], \Gamma_1)$.  These conditions are:
$\phi(\Gamma)< \phi(\Gamma_1)$, $\phi(\Gamma)+1 <
\phi(\Gamma_2[1])$, $\phi(\Gamma_2)\neq \phi(\Gamma_1)$,
$\hom(\Gamma,\Gamma_1)=\hom(\Gamma_2[1],\Gamma_1)=0$.

 By \textbf{C2.4} we see that
$\Gamma$ can be chosen so that $\hom^1(\Gamma, \Gamma_1)\neq 0$
and then by Corollary \ref{RP property 1,2 and.. for Q1} \textbf{(b)}
$\hom(\Gamma,\Gamma_1)=0$.  We have the vanishing
$\hom(\Gamma_2[1],\Gamma_1)=0$ by $\hom(A_2, A_1)=0$. The
inequalities $\phi(\Gamma_1)> \phi(\Gamma)$,   $\phi(\Gamma_2[1])>
\phi(\Gamma)+1$ hold because $\Gamma_1, \Gamma_2$ are components
of $U=A_1 \oplus A_2[-1]$. Finally, we have   $\phi(\Gamma_2)\neq
\phi(\Gamma_1)$ by assumption and the conditions of Lemma
\ref{inequalities} (g) are verified. The lemma
follows.\footnote{We do not need to consider separately the case:
$\Gamma_1\in Ind(A_2)$, $\Gamma_2 \in Ind(A_1)$, for the relation
$\phi(\Gamma_2)\neq \phi(\Gamma_1)$ is symmetric.} \epr

\begin{coro} \label{from C3 to min} Let  $R\in {\mc A}_{exc}$ be a \textbf{C3} object with $\mk{alg}(R)= \begin{diagram}[1em]
A & \rTo      &     &       &   R \\
  & \luDashto &     & \ldTo &       \\
  &           & B[1] &       &
\end{diagram} $.
If $\mk{alg}(R)$ differs from the HN filtration of $R$ or they coincide and  $\phi_{min}<\phi(B)$,  then there exists a
$\sigma$-exceptional triple.
 \end{coro}
 \bpr By the previous lemma  and Corollary \ref{all C3 are final Q1} we can assume that
 $\mk{alg}(R)$ is the HN filtration, hence  $A$ is semistable and $\phi(A)> \phi(B)+1$. If  $\phi_{min}<\phi(B)$, then   $\phi(B)>\phi(S)$ for some $S\in {\mc A}_{exc}\cap \sigma^{ss}$, and by $\phi(A)> \phi(B)+1$ we obtain $\hom^*(A,S)=0$. Since we have $\phi_-(R)=\phi(B)+1 > \phi(S)+1$,
 it follows $\hom^*(R,S)=0$, which due to $\mk{alg}(R)$  gives $\hom^*(B,S)=0$.  Thus, we see that for any  $A' \in Ind(A)$, $B'\in Ind(B)$ the triple $(S,B',A')$ is  semistable and exceptional  with $\phi(S) < \phi(B') < \phi(A')$, $\phi(S)+1 < \phi(A')$. Now the corollary follows from Lemma \ref{inequalities}, \textbf{(a)}.
 \epr
We  obtain now $\sigma$-triples from some $R$-sequences starting  with a \textbf{C1} object $R$.
  \begin{lemma}  \label{coro for C3 after C1} Let $R \in \mc A$ be a \textbf{C1} object. Let
   $\bd R & \rDotsto^{\textbf{C1}} & (S_1,E_1) & \rMapsto^{proj_2} & E_1 & \rDotsto^{\textbf{C3}} & (S_2[1],E_2)\ed $ be an $R$-sequence.
   Then $(S_1,S_2,E_2)$ is an  exceptional triple with $\phi(S_2)+1<\phi_-(E_2)$ and $\hom(S_1, S_2)=0. $

   Furthermore, any of the three  conditions   $ E_2\not \in \sigma^{ss};  \quad \phi(S_2) > \phi_{min};  \quad \phi(S_1)\neq \phi(S_2)+1 $  implies  an existence of a  $\sigma$-exceptional triple.

   \end{lemma}
   \bpr  By  Lemma \ref{lemma for R-se starting with C1}  we see that $(S_1,S_2,E_2)$ is an exceptional triple.
    Since $E_1$  is a \textbf{C3} object, we can write  $\mk{alg}(E_1)= \begin{diagram}[1em]
A' & \rTo      &     &       &   E_1 \\
  & \luDashto &     & \ldTo &       \\
  &           & B'[1] &       &
\end{diagram} $ and(see \textbf{C3.2})  $\phi_-(A')>\phi(B')+1$. From   $E_2 \in Ind(A'), S_2 \in Ind(B')$ we obtain the first property
 $\phi(S_2)+1 < \phi_-(E_2)$.

Next, we consider the vanishing $\hom(S_1, S_2)=0.$   From
\textbf{C3.3} it follows $\hom(E_2,E_1)\neq 0$. As far as $R$ is a
\textbf{C1} object, we  can write $ \mk{alg}(R)=
\begin{diagram}[1em]
A & \rTo      &     &       &   R \\
  & \luDashto &     & \ldTo &       \\
  &           & B &       &
\end{diagram}  $ and $E_1\in Ind(A)$, $S_1\in Ind(B)$. In particular  $E_1$ is a subobject of $R$ in $\mc A$.
Now by $E_1, E_2, R \in \mc A$ and $\hom(E_2,E_1)\neq 0 $ it
follows that  $\hom(E_2,R)\neq 0$, and hence
Corollary \ref{another prop of exc in quiver} implies $\hom(R,E_2)
= 0$. These arguments hold  for each element in $Ind(A')$, hence
$\hom(R, A')=0$. By the exact sequence $\mk{alg}(E_1)$ we get
$\hom(R, B')=0$, and  by the exact sequence $\mk{alg}(R)$ we get
$\hom(B,B')=0$, hence  $ \hom(S_1, S_2)=0. $

 If $E_2 \not \in \sigma^{ss}$,  then we get a $\sigma$-triple from Corollary \ref{all C3 are final Q1}, so let $E_2  \in \sigma^{ss}$.
If $\phi(S_2) > \phi_{min}$, then by Corollary \ref{from C3 to
min} the lemma follows.

Finally, consider the condition $\phi(S_1)\neq \phi(S_2)+1$. Since  we have also $\phi(B'[1])=\phi_-(E_1) \geq
\phi(S_1)$,  we
can write
   $ \phi(S_1) < \phi(S_2)+1. $ We already obtained $\phi(S_2)+1 < \phi(E_2)$ in the beginning of the proof.
 Thus, the triple  $(S_1,S_2,E_2)$ satisfies $ \phi(S_1) < \phi(S_2)+1$,   $\phi(S_2) < \phi(E_2)$, $\phi(S_1) < \phi(E_2)$, $\hom(S_1,S_2)=0$ and by
   Lemma \ref{inequalities} \textbf{(e)} it produces a $\sigma$-exceptional triple.
   \epr

\subsection{Constructions assuming the additional RP property} \label{with the additional RP property}

In this subsection we  restrict $\mc A$ further by assuming that
the properties  in Corollaries \ref{additional RP property},  \ref{coro for isom
triples}
hold.\footnote{to which we refer  as the additional RP property}

In the previous subsection we obtained  a $\sigma$-triple (without
using the additional RP property) from any  long $R$-sequence with
a \textbf{C3} object $R$. One difficulty to  obtain  analogous
criterion when $R$ is  a \textbf{C2} or a \textbf{C1} object  is
mentioned before Lemma \ref{lemma for R-se starting with C1}. It
makes it difficult  to obtain the vanishings
$\{\hom^*(S_1,S_i)=\hom^*(S_1,E_i)\}_{i\geq 2}$  and so to obtain
an exceptional triple. Nevertheless, when $R$ is \textbf{C2}, with
some extra efforts  and utilizing the additional RP property and the property in Corollary \ref{coro for isom
triples} we
obtain exceptional triples in Proposition \ref{all C2 are final}. Furthermore, we show
 that these exceptional triples
can be shifted to    $\sigma$-triples.       We have not an
analogous criterion with a  \textbf{C1} object. \footnote{Lemma \ref{coro
for C3 after C1} and Corollary \ref{after C1} cover all
$R$-sequences  with a \textbf{C1} object $R$ and of length
greater than two.}

\begin{prop} \label{all C2 are final} Each non-final \textbf{C2} object produces   a $\sigma$-exceptional triple.
 \end{prop}
 \bpr  Let  $R\in \mc A$ be a  non-final \textbf{C2} object. Consider the triangle $\mk{alg}(R)$:
   \begin{gather} \label{all C2 are final 15} \begin{diagram}[1em]
A_1 \oplus A_2[-1] & \rTo      &     &       &   R \\
  & \luDashto &     & \ldTo &       \\
  &           & B &       &
\end{diagram}  \ \ \ \ \  \ \begin{array}{c}A_2,B \in {\mc A}\setminus \{0\} \\
                  \hom^1(A_1,A_1)= \hom^1(A_2,A_2)=\hom^1(B,B)=0\\
                 \hom^*(A_1,A_2)= \hom^*(A_1,B)=\hom^*(A_2,B)=0 \end{array}.  \end{gather}
                 For any $\Gamma_0 \in Ind(B)$, $\Gamma \in Ind(A_2[-1])$  we have $\bd R & \rDotsto^{\textbf{C2a}} & (\Gamma_0, \Gamma) \ed$, hence
                  by Corollary  \ref{all C3 are final Q1} if  $\Gamma\not \in \sigma^{ss}$, the proposition follows. Thus, we can assume that   all components of $A_2$ are semistable and $A_1 \neq 0$.

                  For any
   $\Gamma_0 \in Ind(B),\Gamma_1 \in Ind(A_2), \Gamma_2 \in Ind(A_1)$ the triple $(\Gamma_0,\Gamma_1, \Gamma_2)$ is  exceptional, hence by Corollary
   \ref{coro for isom triples} we see that each of $Ind(A_1), Ind(A_2), Ind(B)$ has up to isomorphism unique element. Whence  we can write
   \be \label{A_1,A_2,B} A_1=\Gamma_2^p, \ A_2=\Gamma_1^q, \ B=\Gamma_0^r \qquad (\Gamma_0,\Gamma_1, \Gamma_2) \ \mbox{is exceptional triple}.\ee
  We explained that $\Gamma_1 \in \sigma^{ss}$, furthermore by Lemma \ref{lemma for final good case objects} it follows
   $\phi(\Gamma_1[-1]) > \phi(\Gamma_0)$:
   \be \label{all C2 are final 0}\Gamma_0, \Gamma_1 \in \sigma^{ss},  \qquad  \phi(\Gamma_1) > \phi(\Gamma_0) +1. \ee

      By Proposition  \ref{prop no bad after good}, we know that $\Gamma_2$ is $\sigma$-regular, so   $\mk{alg}(\Gamma_2)$ is  of type  $X \in $ $\{$\textbf{C1}, \textbf{C2}, \textbf{C3}$\}$. We will construct a $\sigma$-exceptional triple in each case.

   \uline{If $\Gamma_2$ is a \textbf{C3} object}, then by Corollary  \ref{all C3 are final Q1}  we can assume that $\Gamma_2$ is final.  For the triangle
  \begin{gather} \label{all C_2 are final 2} \mk{alg}(\Gamma_2)= \begin{diagram}[1em]
A' & \rTo      &     &       &   \Gamma_2 \\
  & \luDashto &     & \ldTo &       \\
  &           & B'[1] &       &
\end{diagram} \ \ \ \qquad \begin{array}{c}A',B' \in {\mc A} \setminus \{0\} \\
                  \hom^1(A',A')= \hom^1(B',B')=0\\
               \hom^*(A',B')=0 \end{array}
 \end{gather} due to  Lemma \ref{lemma for final good case objects in main th for Q1 pr} \textbf{(b)} and  Corollary \ref{from C3 to min}, we can assume also that $A'$ is semistable with  $\phi(A') > \phi(B')+1$ and $\phi(B')=\phi_{min}$. We have also
 $\phi(B')+1=\phi_-(\Gamma_2) \geq \phi_-(A_1) \geq \phi(B)=\phi(\Gamma_0) \geq \phi_{min}=\phi(B')$. Therefore  we can write
\be \label{all C2 are final 1} \phi(A') >\phi(B')+1 \geq
\phi(\Gamma_0)=\phi(B)\geq \phi(B').\ee

For any  $A'' \in Ind(A'), B'' \in Ind(B')$ we have $ \bd R &
\rDotsto^{\textbf{C2b}}& (\Gamma_0,\Gamma_2) & \rMapsto^{proj_2}&
\Gamma_2& \rDotsto^{\textbf{C3}}& (B''[1],A'') \ed $, hence  by
$\deg(\Gamma_0)+1=\deg(B''[1])$ and  Lemma \ref{lemma for R
sequence and eqaul degrees} \textbf{(c)} we get
$\hom(B'',\Gamma_0)=0$, hence \be\label{all C_2 are final 4}
\hom(B',B)=0.\ee

We show now an implication, which  will be used twice later: \be \label{all C2 are
final 7} \mbox{If} \  \ \ \ \hom(B,B')=0 \ \ \mbox{and} \ \ A''
\in Ind(A'), \ \ \  \ \mbox{then} \ \ \ \ A'' \not \cong \Gamma_1.
\ee Indeed, if $A''  \cong \Gamma_1$, then  by  \textbf{C2.4}
applied to \eqref{all C2 are final 15}  and recalling
\eqref{A_1,A_2,B} we obtain  $\hom(B,A')\neq 0$, and then by the
 short exact sequence \eqref{all C_2 are final 2} and $\hom(B,B')=0$ we get $\hom(B,\Gamma_2) \neq 0$. Now from  Corollary \ref{RP property 1,2 and.. for Q1} \textbf{(b)} it follows
 $\hom^1(\Gamma_0,\Gamma_2) = \hom^1(B,A_1)= 0$, which contradicts \textbf{C2.4}.

Keeping \eqref{all C2 are final 1} in mind, we consider two
options $\phi(A') >\phi(B)+1$ and $\phi(A') \leq \phi(B)+1$.

If $\phi(A') >\phi(B)+1$, then $\hom^*(A',B)=0$, which, together
with $\hom^*(\Gamma_2,B)=0$, implies $\hom^*(B',B)=0$. Therefore
(see \eqref{A_1,A_2,B})
$\hom^*(A',\Gamma_0)=\hom^*(B',\Gamma_0)=\hom^*(A',B')=0$, which
by Corollary \ref{coro for isom triples} imply that
$Ind(A')/\cong$, $Ind(B')/\cong$ have unique elements,say
$A'',B''$, and $(\Gamma_0,B'',A'')$  is a semistable exceptional
triple with $\phi(B'')=\phi(B')$, $\phi(A'')=\phi(A')$.  \\
Next, we show that the inequality $\phi(\Gamma_0)\leq \phi(B'')+1$
in \eqref{all C2 are final 1} must be an equality. Indeed, if
$\phi(\Gamma_0)<\phi(B'')+1$, then we have
$\phi(\Gamma_0)<\phi(B'')+1$, $\phi(B'')<\phi(A'')$,
$\phi(\Gamma_0)<\phi(A'')$ and by
  Lemma \ref{inequalities} \textbf{(e)} we can assume $\hom(\Gamma_0,B'')\neq 0$, so $\hom(\Gamma_0,B')\neq 0$. Hence, the triangle $\mk{alg}(\Gamma_2)$ implies
  $\hom(\Gamma_0,A')\neq 0, \hom(\Gamma_0,A'')\neq 0$.   Now Corollary \ref{RP property 1,2 and.. for Q1} \textbf{(b)} implies $\hom^1(\Gamma_0,A'')=\hom^1(\Gamma_0,A')= 0$.
   From the exact sequence $0\rightarrow B'\rightarrow A' \rightarrow \Gamma_2 \rightarrow 0$ and Lemma \ref{lemmaHER} it follows $\hom^1(\Gamma_0,\Gamma_2)= 0$.
   The latter is the same as  $\hom^1(B,A_1)=0$, which contradicts \textbf{C2.4}. So, we obtained $\phi(\Gamma_0)= \phi(B'')+1$ and \eqref{all C2 are final 1} becomes:
    \be \label{all C2 are final 3} \phi(B)=\phi(\Gamma_0)=\phi(B'')+1=\phi(B')+1  \ \ \Rightarrow \ \  \hom(B,B')=0.\ee
Now we utilize the   semistable $\Gamma_1$ in \eqref{all C2 are
final 0}.  If $\phi(\Gamma_1)>\phi(B'')+1$, then
$\hom^*(\Gamma_1,B'')=0$ as well as $\hom^*(\Gamma_1,\Gamma_0)=0$,
hence the  triple $(\Gamma_0,B'',\Gamma_1)$ is exceptional. From
Corollary \ref{coro for isom triples} and the triple
$(\Gamma_0,B'',A'')$ it follows $\Gamma_1 \cong A''$,
   which contradicts \eqref{all C2 are final 7}.
  Therefore  $\phi(\Gamma_1)\leq \phi(B'')+1$. Now \eqref{all C2 are final 3} implies $\phi(\Gamma_1)\leq \phi(B)$.
   Since we consider the subcase $\phi(A') >\phi(B)+1$, therefore $\phi(A')=\phi(A'') >\phi(\Gamma_1)+1$.
   Hence, in addition to $\hom^*(A'',\Gamma_0)=\hom^*(\Gamma_1,\Gamma_0)=0$, we get $\hom^*(A'',\Gamma_1)=0$. Whence, the assumption $\phi(A') >\phi(B)+1$
    leads us to an exceptional triple $(\Gamma_0,\Gamma_1, A'')$. However, the triple $(\Gamma_0,\Gamma_1, \Gamma_2)$ implies  $\Gamma_2 \cong A''$, which contradicts $\Gamma_2 \not \in \sigma^{ss}$, $A''  \in \sigma^{ss}$.

  Therefore, it remains to consider the subcase $\phi(A') \leq \phi(B)+1$. The latter together with  $\phi(B')+1<\phi(A')$, taken from \eqref{all C2 are final 1}, imply $\phi(B')<\phi(B)$.  Combining  with \eqref{all C2 are final 0} and \eqref{all C2 are final 1} we get
   \be \phi(B')<\phi(B)\leq\phi(B')+1 <\phi(A')\leq \phi(B)+1 <\phi(\Gamma_1). \ee
 These inequalities show that, in addition to $\hom(B',B)=0$ (equality \eqref{all C_2 are final 4}) and $\hom^*(\Gamma_1,\Gamma_0)=0$, we get $\hom(B,B')=0$ and $\hom^*(\Gamma_1,B')=0$.
 For clarity, we put together these vanishings:
 \be  \hom(B',\Gamma_0)= \hom(\Gamma_0,B')= 0,  \qquad \ \hom^*(\Gamma_1,\Gamma_0)= \hom^*(\Gamma_1,B')=0. \ee
 The vanishings $\hom^*(\Gamma_1,\Gamma_0)= \hom^*(\Gamma_1,B')=0$ and the additional RP property (Corollary \ref{additional RP property})  show that for each $B'' \in Ind(B')$  the couple $\{\Gamma_0,B''\}$ is not
   Ext-nontrivial, i. e. we have $\hom^1(\Gamma_0,B'')=0$ or $\hom^1(B'',\Gamma_0)=0$. Therefore,  for each $B'' \in Ind(B')$ we have
   $\hom^*(\Gamma_0,B'')=0$ or $\hom^*(B'',\Gamma_0)=0$.  If $\hom^*(\Gamma_0,B'')=0$ for some  $B'' \in Ind(B')$,
  then $(B'',\Gamma_0,\Gamma_1)$ is a semistable exceptional triple with $\phi(B'') < \phi(\Gamma_0)<\phi(\Gamma_1)$, $\hom(B'',\Gamma_0)=0$ and we can apply Lemma \ref{inequalities} \textbf{(b)}. Hence, we can assume that for each  $B'' \in Ind(B')$ we have  $\hom^*(B'',\Gamma_0)=0$  and $(\Gamma_0,B'',\Gamma_1)$ is an exceptional triple.  Therefore the set $Ind(B')/\cong$ has unique element,
 say $B''$. Thus, we arrive at an exceptional triple \be (\Gamma_0,B'',\Gamma_1), \ \ \hom(\Gamma_0,B'')=0, \ \ \qquad B'\cong (B'')^s. \ee
 On the other hand, the vanishings  $\hom^*(B',\Gamma_0)=\hom^*(\Gamma_2,\Gamma_0)=0$ and the triangle \eqref{all C_2 are final 2} imply $\hom^*(A',\Gamma_0)=0$. The last vanishing and $\hom^*(A',B')=0$ give rise  to a  triple $(\Gamma_0,B'', A'')$ with  $(A'')^u \cong A'$. Both the triples $(\Gamma_0,B'', A'')$, $(\Gamma_0,B'',\Gamma_1)$ imply $A'' \cong \Gamma_1$, which
 contradicts \eqref{all C2 are final 7}.  Thus, the proposition follows,   when $\Gamma_2$ is a \textbf{C3} object.

  \uline{If $\Gamma_2$ is a \textbf{C2} object},  then $\mk{alg}(\Gamma_2)$ and some of its features are
 \begin{gather} \begin{diagram}[1em]
A_1' \oplus A_2'[-1] & \rTo      &     &       &   \Gamma_2 \\
  & \luDashto &     & \ldTo &       \\
  &           & B' &       &
\end{diagram}  \qquad \begin{array}{c}A_2',B' \in {\mc A}\setminus \{0\} \\
                  \hom^1(A_1',A_1')= \hom^1(A_2',A_2')=\hom^1(B',B')=0\\
                 \hom^*(A_1',A_2')= \hom^*(A_1',B')=\hom^*(A_2',B')=0. \end{array}\ \nonumber \end{gather}
For any $A'' \in Ind(A_1' \oplus A_2'[-1]), B'' \in Ind(B')$ we have an $R$-sequence \\
$\bd R & \rDotsto^{\textbf{C2b}}& (\Gamma_0,\Gamma_2) &
\rMapsto^{proj_2}& \Gamma_2& \rDotsto^{\textbf{C2a/b}}& (B'',A'')
\ed $ without a \textbf{C3}-step in it. From  Corollary \ref{coro
for R seq initial C3,...} (the last case) it follows
$\hom^*(B',\Gamma_0)=
\hom^*(A_1',\Gamma_0)=\hom^*(A_2',\Gamma_0)=0 $. Combining these
vanishings  with $\hom^*(A_1', B')=\hom^*(A_2', B')=0$,
$A_2'\neq 0$ we  conclude  by Corollary \eqref{coro for isom
triples} that
\begin{gather}\label{all C2 are final 11}  A_2'\cong (A'')^s; \  B'\cong (B'')^t;  \  (\Gamma_0,B'', A'') \ \mbox{is exceptional;} \
\ \mbox{if}  \ \ A_1'\neq 0  \ \ \mbox{then} \ \ A_1'\cong
(A'')^u
\end{gather}  for some  $A'', B'' \in {\mc A}_{exc}$.  By Corollary \ref{all C3 are final Q1} and $\bd \Gamma_2& \rDotsto^{\textbf{C2a}}& (B'',A''[-1]) \ed $ we reduce to the case  $A''\in \sigma^{ss}$. Thus,   $\Gamma_2$ becomes  final.
Furthermore, by $\deg(B)=\deg(B')$  we have $\deg(\Gamma_0)=\deg(B'')$ and we see that the  $R$-sequence   $\bd R & \rDotsto^{\textbf{C2b}}& (\Gamma_0,\Gamma_2) & \rMapsto^{proj_2}& \Gamma_2& \rDotsto^{\textbf{C2a}}& (B'',A''[-1]) \ed $ satisfies  the three conditions of Proposition \ref{non-existence of some sequeces}. This proposition ensures a $\sigma$-exceptional triple.  It remains to consider:

  \uline{$\Gamma_2$ is a \textbf{C1} object}. Denote the corresponding triangle as follows:
  \begin{gather} \label{all C_2 are final 8} \mk{alg}(\Gamma_2)= \begin{diagram}[1em]
A' & \rTo      &     &       &   \Gamma_2 \\
  & \luDashto &     & \ldTo &       \\
  &           & B' &       &
\end{diagram} \ \ \ \qquad  \begin{array}{c}A',B' \in {\mc A} \setminus \{0\} \\
                  \hom^1(A',A')= \hom^1(B',B')=0\\
               \hom^*(A',B')=0. \end{array} \end{gather}
             Now we have again $\deg(B')=\deg(\Gamma_0)$. It follows from  Corollaries \ref{coro for R seq initial C3,...}, \ref{coro for isom triples} that
             \begin{gather}\label{all C2 are final 13}  A'\cong (A'')^s, \ \ B'\cong (B'')^t,  \ \ \ (\Gamma_0,B'', A'') \ \mbox{is exceptional,}\
\\ \label{all C2 are final 14} \phi(B'')>\phi(\Gamma_0).
\end{gather} for some  $A'', B'' \in {\mc A}_{exc}$.  The arguments which give \eqref{all C2 are final 13} are as those giving  \eqref{all C2 are final 11}, and  \eqref{all C2 are final 14} follows from Corollary \ref{coro for R seq initial C3,...} \textbf{(b)}.
  If $A''\in \sigma^{ss}$, then $\Gamma_2$ is final, and  Proposition \ref{non-existence of some sequeces} produces  a  $\sigma$-sequence  from the $R$-sequence   $\bd R & \rDotsto^{\textbf{C2b}}& (\Gamma_0,\Gamma_2) & \rMapsto^{proj_2}& \Gamma_2& \rDotsto^{\textbf{C1}}& (B'',A'') \ed $.     Therefore, we can  assume that $A'' \not \in \sigma^{ss}$.

  If $A''$ is  \textbf{C1} or \textbf{C2}, then we get an $R$-sequence, in which a \textbf{C3} step does not appear as follows:
  \begin{gather}
\bd[height=1.5em]
R & \rDotsto^{\textbf{C2b}}& (\Gamma_0,\Gamma_2)         & \rMapsto^{proj_2}& \Gamma_2 & \rDotsto^{\textbf{C1}}& (B'',A'')  & \rMapsto^{proj_2}& A'' & \rDotsto^{X_3}& (S,E)& \rMapsto^{proj_2}& E\\
  &               & \dMapsto^{proj_1} &                  &    &               & \dMapsto^{proj_1}&                  &     &               &  \dMapsto^{proj_1} &  &   \\
  &               &  \Gamma_0             &                  &    &               &  B''             &                  &     &               &  S      &  &    \ed   \ \ X_3 \in \{ \textbf{C1} , \textbf{C2a}, \textbf{C2b} \}. \nonumber \end{gather}
 From Corollary
 \ref{coro for R seq initial C3,...} it follows that the sequence $(\Gamma_0, B'',S,E)$ is  exceptional, which contradicts Corollary \ref{coro for isom triples}.

 Therefore $A''$ must be a \textbf{C3} object, which ensures a  $\Gamma_2$-sequence  of the form \\
 $\bd \Gamma_2 & \rDotsto^{\textbf{C1}} & (B'',A'') & \rMapsto^{proj_2} & A'' & \rDotsto^{\textbf{C3}} & (S[1],E)\ed $. In
 Lemma \ref{coro for C3 after C1} is shown that   the triple   $(B'', S, E)$   is  exceptional. The  criteria given there show that $E \in \sigma^{ss}$ and  reduce the phases of $(B'', S, E)$ to
 \be \label{all A_2 stable 3}  \phi(B'') = \phi(S)+1=\phi_{min}+1 < \phi(E); \qquad (B'', S, E) \ \mbox{is semistable and exceptional}. \ee From Corollary  \ref{coro for isom triples} it follows that
 $  \mk{alg}(A'')=  \begin{diagram}[1em]
E^i & \rTo      &     &       &   A'' \\
  & \luDashto &     & \ldTo &       \\
  &           & S[1]^j &       &
\end{diagram} $ for some integers  $i,j \in \NN$.

If $\phi(E)>\phi(\Gamma_0)+1$, then $\hom^*(E, \Gamma_0)=0$, which,
combined  with $\hom^*(A'',\Gamma_0)=0$ (see \eqref{all C2 are
final 13}),  implies $\hom^*(S,\Gamma_0)=0$. These vanishings and the exceptional triples
$(\Gamma_0, B'', A'')$,  $(B'', S, E)$ imply that $(\Gamma_0,
B'',S, E)$ is an exceptional sequence, which is impossible.

Thus, $\phi(E)\leq \phi(\Gamma_0)+1$ and we can write (see also
\eqref{all C2 are final 0})\be \label{all A_2 are stabel}
\phi(S)+1<\phi(E)\leq \phi(\Gamma_0)+1<\phi(\Gamma_1) \ \ \ \
\Rightarrow \ \ \ \   \hom^*(\Gamma_1,S)=0. \ee Since
$\hom^*(\Gamma_1,\Gamma_0)=0$ as well,  the additional RP
property(Corollary \ref{additional RP property}) ensures  that the
couple $\{S,\Gamma_0\}$ is not  Ext-nontrivial, therefore
$\hom^1(\Gamma_0,S)=0$ or $\hom^1(S, \Gamma_0)=0$. We show below
that $\hom(\Gamma_0,S)=\hom(S,\Gamma_0)=0$, hence
$\hom^*(\Gamma_0,S)=0$ or $\hom^*(S, \Gamma_0)=0$. It follows
that some of the triples $ \ (S,\Gamma_0,\Gamma_1), 
(\Gamma_0,S,\Gamma_1) $  is exceptional.

If  $(S,\Gamma_0,\Gamma_1)$ is exceptional, then  Lemma
\ref{inequalities}, \textbf{(a)}  produces $\sigma$-exceptional
triple, due to the inequalities $\phi(S)<\phi(\Gamma_0)$,
$\phi(\Gamma_0)+1<\phi(\Gamma_1)$ (see \eqref{all A_2 are
stabel}).

If $(\Gamma_0,S,\Gamma_1)$ is exceptional, then due to the inequalities
$\phi(S)<\phi(\Gamma_1)$, $\phi(\Gamma_0)<\phi(\Gamma_1)$,
  $\phi(\Gamma_0)<\phi(S)+1$  (the last  comes from  \eqref{all C2 are final 14},
\eqref{all A_2 stable 3}) and $\hom(\Gamma_0,S)=0$ we can apply
Lemma \ref{inequalities} \textbf{(e)}.

The used in advance $\hom(\Gamma_0,S)=0$ follows from
$\phi(S)<\phi(\Gamma_0)$ (see \eqref{all A_2 are stabel}). The other vanishing
$\hom(S, \Gamma_0)=0$ follows from $\phi_-(A'')\geq
\phi(\Gamma_0)$, $\hom^*(A'',\Gamma_0)=0$ (see \eqref{all C2 are
final 13}), and Lemma \ref{lemma for hom leq 1(X,S)}.

 Now the proposition is
completely proved.
      \epr
   It follows now   the  \textbf{C2}-analogue of    Corollary \ref{from C3 to min}. After a proper reformulation,\footnote{The part of Corollary \ref{from C3 to min} using $\phi_{min}$   can be reformulated as  saying that the data: a final \textbf{C3} object $R\in \mc A_{exc}$, $\bd R&\rDotsto& (S,F) \ed$, $X\in \{S,F\}$, $\deg(X)\neq 0$, and  $\phi(X)-\deg(X)>\phi_{min}$ implies a $\sigma$-triple.} Corollary \ref{from C3 to min} is transformed to Corollary \ref{from C2 to max} by replacing ``\textbf{C3}'' with ``\textbf{C2}'' and ``$>\phi_{min}$'' with ``$<\phi_{max}$''.
      \begin{coro} \label{from C2 to max} Let $R\in \mc A$ be a  \textbf{C2} object with $\mk{alg}(R)= \begin{diagram}[1em]
A_1 \oplus A_2[-1] & \rTo      &     &       &   R \\
  & \luDashto &     & \ldTo &       \\
  &           & B   &       &
\end{diagram} $.  If either $\mk{alg}(R)$ differs from the HN filtration of $R$  or they coincide and $\phi(A_2)<\phi_{max}$,  then there exists a $\sigma$-triple.
 \end{coro}
 \bpr  Due to the criteria given in  Proposition \ref{all C2 are final} and  Lemma \ref{lemma for final good case objects in main th for Q1 pr},   we reduce to the case:
  $R$ is final and $\mk{alg}(R)$ is the HN filtration of $R$. In particular   $A_1 \oplus A_2[-1]  \in \sigma^{ss}$.

 If $\phi(A_2)<\phi_{max}$, then $\phi(S)>\phi(A_2)$ for some $S \in {\mc A}_{exc}\cap \sigma^{ss}$. Since  $\mk{alg}(R)$ is the HN filtration of $R$, it follows that $\phi_+(R)=\phi(A_2)-1$.
 Therefore  $\phi(S)>\phi_+(R) +1>\phi(B)+1$, which implies $\hom^*(S,R)=\hom^*(S,B)=0$. From the triangle  $\mk{alg}(R)$ we obtain also $\hom^*(S,A_2)=0$.  Therefore, for any  $A'\in Ind(A_2)$, $B'\in Ind(B)$ the semistable triple
  $(B',A',S)$ is   exceptional and it satisfies $\phi(B') < \phi(A') < \phi(S)$, $\phi(B')+1 < \phi(S)$. Now
  Lemma \ref{inequalities} \textbf{(a)} produces a $\sigma$-triple.
 \epr

In the next corollary we obtain $\sigma$-triples from some, but not all, long $R$-sequences with a \textbf{C1} object $R$.

\begin{coro} \label{after C1} Let $ \bd R & \rDotsto^{\textbf{C1}}& (S_1,E_1) \ed $. If $E_1$ is either a  \textbf{C2} or  a \textbf{C1} object, then there exists  a $\sigma$-exceptional triple.
\end{coro}
  \bpr
   If $E_1$ is  \textbf{C2}, then we have an  $R$-sequence
   $\bd R & \rDotsto^{\textbf{C1}} & (S_1,E_1) & \rMapsto^{proj_2} & E_1 & \rDotsto^{\textbf{C2a}} & (S_2,E_2[-1])\ed$. By Proposition \ref{all C2 are final}, we can assume that $E_1$ is final, and then  Proposition \ref{non-existence of some sequeces} ensures a $\sigma$-triple.

   If $E_1$ is  \textbf{C1}, then we get a second step  $\bd E_1 & \rDotsto^{\textbf{C1}} & (S_2,E_2)  \ed $, for some $(S_2,E_2) $, and then we go on further   until  a final object occurs, which will certainly happen by Lemma \ref{R seq are finite}.    We can assume that in this process a \textbf{C2}  step does not occur (otherwise the corollary follows by the proven case). By Corollary \ref{all C3 are final Q1} we can assume that all
    \textbf{C3} objects are final. Hence, if a \textbf{C3} step occurs, then this is the last step. The other possibility   is to reach a final \textbf{C1} case and then   Proposition  \ref{non-existence of some sequeces}  gives a $\sigma$-triple.
    Whence, we reduce to  an $R$-sequence with $n\geq 3$ of the form:
  \begin{gather}
\bd[height=1.5em]
R & \rDotsto^{\textbf{C1}}& (S_1,E_1)         & \rMapsto^{proj_2}& E_1 & \rDotsto^{\textbf{C1}}& (S_2,E_2)         & \rMapsto^{proj_2}& E_2 & \rDotsto^{\textbf{C1}}& \dots & \rMapsto^{proj_2}&E_{n-1}& \rDotsto^{\textbf{C3}}&(S_n,E_n) & \rMapsto^{proj_2}& E_n \\
  &               & \dMapsto^{proj_1} &                  &    &               & \dMapsto^{proj_1}&                  &     &               &       &  & & & \dMapsto^{proj_1} &  &   \\
  &               &  S_1             &                  &    &               &  S_2             &                  &     &               & \dots & &  & & S_n      &  &    \ed. \nonumber \end{gather}
 We apply  Lemma \ref{lemma for R-se starting with C1}  to the $R$-sequence above and to the $E_1$-sequence in it, and obtain:\\
$
\hom^*(S_n,S_1)=\hom^*(E_n,S_1)=\hom^*(S_n,S_2)=\hom^*(E_n,S_2)=0
$.
  Furthermore, by Lemma \ref{lemma for R sequence and eqaul degrees} \textbf{(b)} and $\deg(S_2)=\deg(S_1)=0$ (see table \eqref{table with degrees}) it follows
 $ \hom^*(S_2,S_1)=0. $     These vanishings imply that $(S_1,S_2,S_n,E_n)$ is a semistable exceptional sequence, which is a contradiction.
  \epr

 We summarize now the  results concerning  $R$-sequences with a \textbf{C1} object $R$.
 \begin{coro} \label{lemma for final good case objects in main th for Q1} Let  there be no a $\sigma$-exceptional triple.
   If
   $\bd R & \rDotsto^{\textbf{C1}} & (S_1,E_1) \ed $, then the object $E_1$ is either semistable or a \textbf{C3} object. If $E_1$ is a  \textbf{C3} object, then   for each $R$-sequence\\
$\bd R & \rDotsto^{\textbf{C1}} & (S_1,E_1) & \rMapsto^{proj_2} &
E_1 & \rDotsto^{\textbf{C3}} & (S_2[1],E_2)\ed $  the triple
$(S_1,S_2,E_2)$ is   exceptional, semistable, and it satisfies:   $\phi(S_2)=\phi_{min}$, $
\phi(S_1)=\phi(S_2)+1<\phi(E_2),  \ \ \hom(S_1, S_2)=0, \ \
\hom^1(S_1, S_2) \neq 0. $
\end{coro}
\bpr  Follows from Corollary \ref{after C1} and Lemma \ref{coro for C3
after C1}. \epr

 A next step to the proof of Proposition
\ref{between min and max} is to show that, given a \textbf{C1}-object $R$, each long $R$-sequence induces a $\sigma$-triple, when    $R$ is  part of an exceptional pair $(R,S_{max})$ or $(S_{min},R)$.

 \begin{lemma} \label{when C1 is final} Let  $R \in \mc A$ be a  non-final \textbf{C1} object. If we are given one of the following:
\begin{itemize}
    \item[\textbf{(a)}] $S_{min} \in {\mc A}_{exc}$ with $\phi(S_{min})=\phi_{min}$ and $\hom^*(R, S_{min})=0$,
    \item[\textbf{(b)}]  $S_{max} \in {\mc A}_{exc}$ with $\phi(S_{max})=\phi_{max}$ and $\hom^*(S_{max},R)=0$,
\end{itemize}
then there  exists   a $\sigma$-exceptional triple.
  \end{lemma}
  \bpr By the criterion given in  Corollary \ref{after C1}  we can assume that there exists an  $R$-sequence  of the form
  $\bd R & \rDotsto^{\textbf{C1}} & (S_1,E_1) & \rMapsto^{proj_2} & E_1 & \rDotsto^{\textbf{C3}} & (S_2[1],E_2)\ed $.
  The triple   $(S_1,S_2,E_2)$ is exceptional  by Lemma \ref{coro for C3 after C1} and using the criteria given there we can assume that
  it is semistable and:
    \ben \phi(S_1)=\phi(S_2)+1=\phi_{min}+1<\phi(E_2),\ \quad  \phi(S_2)=\phi_{min}. \een
      In part (a) we are given that $\hom^*(R, S_{min})=0$. We claim that the triple  $(S_{min},S_1,E_2)$ is exceptional. Indeed, we have:  $\hom^*(E_2,S_{min})=0$ by
$\phi(E_2)>\phi(S_{min})+1$, and $\hom^*(E_2,S_1)=0$ by the
exceptional triple $(S_1,S_2,E_2)$. Finally
$\hom^*(S_1,S_{min})=0$ by $\hom^*(R,S_{min})=0$,
$\phi_-(R)=\phi(S_1) \geq \phi(S_{min})$ and Lemma \ref{lemma for
hom leq 1(X,S)}. Thus, we constructed a semistable exceptional
triple $(S_{min},S_1,E_2)$ with
$\phi(S_{min})<\phi(S_1)=\phi(S_{min})+1<\phi(E_2)$.  Now Lemma
\ref{inequalities} \textbf{(a)} produces a $\sigma$-triple.

Let $\hom^*(S_{max}, R)=0$ for some $S_{max}\in {\mc A}_{exc}$
with maximal phase.  Unfolding the definition of \textbf{C1} we
get a short exact sequence
 $0 \rightarrow E \rightarrow R \rightarrow S \rightarrow 0$ with
$E_1 \in Ind(E)$, $S_1 \in Ind(S)$, $\phi(S)=\phi(S_1)$.  Since
$S_{max}$ is of maximal phase, we have $\phi(S_{max})\geq
\phi(E_2)
>\phi(S_2)+1=\phi(S_1)=\phi(S)$, which implies
$\hom^*(S_{max},S_2)=0$, $\hom(S_{max},S)=0$.   By Lemma
\ref{lemmaHER} and $\hom^*(S_{max},R)=0$ we get also
$\hom(S_{max},S[1])=0$, hence $\hom^*(S_{max},S)=0$, which in turn
implies $\hom^*(S_{max},E)=0$. So far, using the conditions of
\textbf{(b)}, we obtained \be
\hom^*(S_{max},S_1)=\hom^*(S_{max},E_1)=\hom^*(S_{max},S_2)=0.\ee
We  show below that $\hom^*(S_{max},E_2)$ also vanishes, and then
the sequence $(S_1,S_2,E_2, S_{max})$ becomes exceptional, which
is a contradiction.  Then the corollary follows.

Since  any relation of the form $\bd E_1 & \rDotsto^{\textbf{C3}} & (X[1],Y)\ed $ gives by  Lemma \ref{coro for C3 after C1} an  exceptional
triple  $(S_1,X,Y)$, it follows from  Corollary \ref{coro for isom triples}
 that  $\mk{alg}(E_1)=  \begin{diagram}[1em]
E_2^i & \rTo      &     &       &  E_1 \\
  & \luDashto &     & \ldTo &       \\
  &           & S_2[1]^j &       &
\end{diagram}$. This triangle and the already shown $\hom^*(S_{max},E_1)=\hom^*(S_{max},S_2)=0$ give the desired  $\hom^*(S_{max},E_2)=0$.
  \epr

The additional RP property   gives us another situation,
 where  the irregular  cases
\textbf{B1} and \textbf{B2} cannot occur. This is shown in Lemmas \ref{no B2 0}, \ref{no bad between min and max} below.  In this respect these lemmas  are similar to
Proposition \ref{prop no bad after good}, but the latter uses RP properties 1,2.
\begin{lemma}  \label{no B2 0}  If $(S_{min}, E)$ is an exceptional  pair in $\mc A$ with  $S_{min}  \in\mc P(\phi_{min})$, then $E$  is not  \textbf{B2}. \end{lemma}
 \bpr
If $E$ is a \textbf{B2} object, then  $ \mk{alg}(E) =
\begin{diagram}[1em]
A& \rTo      &     &       &   E \\
  & \luDashto &     & \ldTo &       \\
  &           & B[1] &       &
\end{diagram}$ with  $B\in \sigma^{ss}$,   $\phi(B)+1=\phi_-(E)$, $ \phi_-(A) >\phi(B)+1  $, and for some
$\Gamma \in Ind(B)$ the couple $\{E,\Gamma\}$ is  Ext-nontrivial.
From $\Gamma \in {\mc A}_{exc}\cap \sigma^{ss}$ it follows  that
$\phi(\Gamma)=\phi(B)\geq \phi_{min}$, therefore
$\phi_-(A)>\phi_{min}+1$ and $\hom^*(A,S_{min})=0$. The vanishings
$\hom^*(A,S_{min})=0$, $\hom^*(E,S_{min})=0$ imply
$\hom^*(B,S_{min})=0$. Thus, we obtain an Ext-nontrivial couple
$\{\Gamma, E\}$ and $S_{min} \in {\mc A}_{exc}$ with
$\hom^*(E,S_{min})= \hom^*(\Gamma,S_{min})=0$, which contradicts
the additional RP property (Corollary \ref{additional RP
property}). \epr

\begin{lemma} \label{no bad between min and max}  Let $\phi_{max} > \phi_{min}+1$.  If $(S_{min}, E, S_{max})$ is an exceptional  triple in $\mc A$  with $S_{min} \in \mc P(\phi_{min})$,  $S_{max} \in \mc P(\phi_{max})$, then $E$
 is not  $\sigma$-irregular. \end{lemma}
 \bpr In the previous lemma we showed that $E$ is not a \textbf{B2} object. Suppose that  $E$ is a \textbf{B1} object. Then  $ \mk{alg}(E) = \begin{diagram}[1em]
A_1 \oplus A_2[-1] & \rTo      &     &       &   E \\
  & \luDashto &     & \ldTo &       \\
  &           & B &       &
\end{diagram}$ with $B\in \sigma^{ss}$, $\phi_-(A_1 \oplus A_2[-1])\geq \phi(B)$, $\phi(B)=\phi_-(E)$, and for some  $\Gamma \in Ind(A_2)$ the couple $\{E,\Gamma\}$ is  Ext-nontrivial.

If $\phi(B) > \phi(S_{min})$, then we have $\phi_-(\Gamma[-1])
\geq \phi_-(A_1 \oplus A_2[-1])\geq \phi(B)  >\phi(S_{min}) $,
hence  $\phi_-(\Gamma)  >\phi(S_{min})+1 $. However, this implies
$\hom^*(\Gamma, S_{min}) = 0$ and we have also $\hom^*(E, S_{min})
= 0$, which contradicts the additional RP property(Corollary \ref{additional RP property}).

 If
$\phi(B) \leq \phi(S_{min})$, then by
$\phi_{max} > \phi_{min}+1$ we have $\hom^*(S_{max},B)=0$, which, combined with
$\hom^*(S_{max},E)=0$ and the triangle $\mk{alg}(E)$, implies
$\hom^*(S_{max},A_2)=0$. Thus,  we have  $\hom^*( S_{max},\Gamma)
=\hom^*( S_{max},E) = 0$, which contradicts  Corollary
\ref{additional RP property}. \epr We can prove now easily:
\begin{prop} \label{between min and max} Let $\phi_{max}-\phi_{min} > 1$.  Let $(S_{min}, E, S_{max})$ be an exceptional  triple in $\mc A$
 with  $S_{min} \in \mc P(\phi_{min})$,  $S_{max} \in \mc P(\phi_{max})$.  If $E \not \in \sigma^{ss}$,  then there exists a $\sigma$-exceptional triple.
\end{prop}
\bpr From Lemma \ref{no bad between min and max}  and $E\not \in
\sigma^{ss}$  it follows that $E$ is regular. From Corollary
\ref{between min max 1} it follows that  $E$ cannot be final
(due to Corollary \ref{coro for isom
triples} there are not exceptional sequences of length 4).  Now the
existence of a $\sigma$-exceptional triple follows from Corollary
\ref{all C3 are final Q1}, Proposition \ref{all C2 are final}, and
Lemma \ref{when C1 is final}. \epr

\section {Application to  \texorpdfstring{$\st(D^b(Q_1))$}{\space}} \label{main theo for Q_1}
 The criteria of Section \ref{constructing} hold for ${\mc A}=Rep_k(Q_1)$, due to Section \ref{two examples}.
In this section we apply these criteria  to $Rep_k(Q_1)$.   The result  is  the following theorem:
\begin{theorem}\label{main theorem for Q_1} Let $k$ be an algebraically closed field.
For each $\sigma \in \st(D^b(Rep_k(Q_1)))$ there exists a
$\sigma$-exceptional triple.
\end{theorem}
In Remark \ref{weaker} we pointed out a variant of  Sections \ref{non-stable exc obj in...}, \ref{some terminilogy}, \ref{no bad after good}, \ref{sequence}, \ref{final}, \ref{constructing} in which  $k$ is any field.
 We cannot point out a variant of  Theorem \ref{main theorem for Q_1} without the restriction  that    $k$ is  algebraically closed.

\begin{coro}\label{connectedness}  The manifold  $\st(D^b(Rep_k(Q_1)))$ is connected.
\end{coro}
\bpr  Let  $\mc E=(E_1^0,M,E_3^0)$.  Let $\Sigma_{\mc E}$ be as in
\eqref{theta_{mc E}}. From Corollary \ref{RP property 1,2 and..
for Q1} \textbf{(b)}   we see that all triples in $D^b(Q_1)$ are regular.
 Therefore
 $\Sigma_{\mc E}$ is connected \cite[Corollary 3.20]{Macri}.
     From \cite{WCB1} it follows that  all exceptional triples in $D^b(Q_1)$ are obtained by shifts and   mutations of  $\mc E$.
  Recalling  Corollary \ref{coro of Macri}   we see that  Theorem \ref{main theorem for Q_1 in intro} is the same as the equality $\st(D^b(Q_1))=\Sigma_{\mc E}$. The corollary follows.
\epr

Throughout the proof of Theorem \ref{main theorem for Q_1}(the entire Section \ref{main theo
for Q_1}) we fix  the  notations ${\mc A} = Rep_{k}(Q_1)$ and ${\mc T} = D^b(\mc A)$.
We  prove the theorem by
contradiction.

\uline{Let $\sigma
\in \st(D^b({\mc A}))$. In all subsections of Section  \ref{main theo for Q_1}, except  subsection \ref{basic lemmas}, we assume that there does not
 exist a $\sigma$-exceptional triple.}

Loosely speaking,  this assumption leads   to  certain
``non-generic'' situations (see \eqref{non-generic}). However,
using the  locally finiteness of $\sigma$, we show that these
situations cannot occur (Corollaries \ref{coro for E_2 and E_1},
\ref{coro for E_3 and E_4}) and so we get a contradiction.

The notations $M, M', E_1^m, E_2^m, E_3^m, E_4^m$ are explained in
Proposition \ref{exceptional objects in Q1}. We will refer  often
to table \eqref{Q1 table} and Corollary \ref{exceptional
colleections}. Whenever  we claim that a triple $(A_0,A_1,A_2)$ is
an exceptional triple(with $A_0,A_1,A_2$ one of the symbols $M,
M', E_1^m, E_2^m, E_3^m, E_4^m$),  then we refer implicitly to
Corollary \ref{exceptional colleections}, and whenever we discuss
$\hom^*(A,B)$ with $A,B$ varying in these symbols, we refer to
table \eqref{Q1 table}.

\begin{remark}  \label{inc dec seq} Recall that(see right after Definition \ref{def of sigma_-})  $\hom(A,B) \neq 0$ implies  $\phi_-(A)\leq \phi_+(B)$.
Using table  \eqref{Q1 table} we can write  for any $n \in \NN$
\begin{itemize}
    \item $\hom(E_1^{n+1}, E_1^n) \neq 0  $  hence  $\phi_-(E_1^{n+1})\leq \phi_+(E_1^{n})$
    \item $\hom(E_2^{n}, E_1^{n+1}) \neq 0$  hence $\phi_-(E_2^{n})\leq \phi_+(E_2^{n+1})$
    \item $\hom(E_3^{n}, E_3^{n+1}) \neq 0$   hence $\phi_-(E_3^{n})\leq \phi_+(E_3^{n+1})$
        \item $\hom(E_4^{n+1}, E_4^{n}) \neq 0$   hence $\phi_-(E_4^{n+1})\leq \phi_+(E_4^{n})$.
\end{itemize}

\end{remark}

\subsection{Basic lemmas} \label{basic lemmas} The facts explained  here  are   basic tools used in the following subsections. These facts are   individual for $Q_1$. The reader may skip this subsection on a first reading and return to it only  when we refer to these tools.

In this subsection we do not put any restrictions on $\sigma \in \st(\mc T)$. In all the rest subsections $\sigma$ is assumed not to admit a $\sigma$-exceptional triple.

\subsubsection{Useful short exact sequences in ${\mc A}$ and two  corollaries based on locally finiteness.} It is easy to check: \begin{lemma} There exist arrows in $\mc A$ as shown below, so  that  the resulting  sequences are  exact:
\begin{gather} \label{ses1} \bd 0 & \rTo & E_2^{m-1} & \rTo & E_1^m & \rTo & ( E_1^0 )^2 & \rTo & 0 \ed \\
\label{ses2} \bd 0 & \rTo & E_3^{m} & \rTo & E_2^m & \rTo & M & \rTo & 0 \ed \\
\label{ses3} \bd 0 & \rTo & E_3^{m-1} & \rTo & E_4^m & \rTo & ( E_4^0 )^2 & \rTo & 0 \ed  \\
\label{ses4} \bd 0 & \rTo & M & \rTo & E_4^m & \rTo & E_1^m & \rTo & 0 \ed
\end{gather}
\end{lemma}
These short exact sequences combined with the locally finiteness of $\sigma$ result in  Corollaries \ref{coro for E_2 and E_1}, \ref{coro for E_3 and E_4}. These corollaries  exclude the following  two  situations: \begin{gather} \label{non-generic} \{  E_2^m  \}_{m \in \NN}\subset \mc P(t),  \{  E_1^m  \}_{m \in \NN}\subset \mc P(t+1)  \ \  \mbox{or } \ \  \{  E_3^m  \}_{m \in \NN}\subset \mc P(t),  \{  E_4^m  \}_{m \in \NN}\subset \mc P(t+1). \end{gather}
We will sometimes  refer to these two cases  as non-locally finite cases.
\begin{coro} \label{coro for E_2 and E_1} Assume that  $ \{ E_1^m, E_2^m \}_{ m \in \NN }\subset \sigma^{ss}$  and $\{  E_2^m  \}_{m \in \NN}\subset \mc P(t) $ for some $t \in \RR$. Then for
each  $m\in \NN $ we have $t \leq \phi(E_1^m) \leq t+1$, and  there exists $n \in \NN $ with $t \leq \phi(E_1^n) < t+1$.
\end{coro}
\bpr By table \eqref{Q1 table} we have $\hom(E_2^m, E_1^n) \neq 0$ and $ \hom(E_1^n,E_2^m[1])\neq 0$ for $m\geq 1$,
 hence $t=\phi(E_2^m)\leq \phi(E_1^n) \leq \phi(E_2^m)+1=t+1$. It remains to show the last claim.

 The short exact sequence \eqref{ses1}  gives a distinguished triangle
$\bd[1em]  E_1^{m} & \rTo & (E_1^0)^2 & \rTo & E_2^{m-1}[1] & \rTo & E_1^{m}[1] \ed$. Suppose that $\phi(E_1^m)=t+1$ for  each  $m$. Then      $\{ E_1^m$, $(E_1^0)^2$, $ E_2^{m-1}[1]\}_{m\in\NN} \subset \mc P(t+1)$.  It follows  that $\bd[1em] 0& \rTo & E_1^{m} & \rTo & (E_1^0)^2 & \rTo & E_2^{m-1}[1] & \rTo & 0 \ed$ is a short exact sequence  in the abelian category $\mc P(t+1)$ for each $m\in \NN$ (see the  beginning of subsection \ref{comments on stab cond}). Hence  $ \bd[1em] E_1^{m} & \rTo & (E_1^0)^2  \ed $ is a monic arrow in $\mc P(t+1)$
 for each $m \in \NN$. It follows by  Lemma \ref{finiteness coro} that the set  $\{[E_1^{m}]\}_{m \in \NN}$ is a finite subset of $K(D^b(\mc A))$.
 On the other hand (see Lemma \ref{exceptional objects in Q1}) we can write
$ \{ \ \ [E_1^{m}]=(m+1)[E_1^0] + m [M] + m [E_3^0] \ \ \}_{m \in
\NN},$ which is infinite in $K(D^b(\mc A))$. Thus,  the assumption
that  $\phi(E_1^n)=t+1$ for each $n$ leads to a contradiction.
\epr

\begin{coro} \label{coro for E_3 and E_4} Assume that  $ \{ E_3^m, E_4^m  \}_{m \in \NN }\subset \sigma^{ss}$  and $\{  E_3^m  \}_{m \in \NN }\subset \mc P(t) $ for some $t \in \RR$.
Then for each $m\in \NN $ we have $t \leq \phi(E_4^m) \leq t+1$, and  there exists $l \in \NN $ with $t \leq \phi(E_4^l) < t+1$.
\end{coro}
\bpr By  table \eqref{Q1 table} we have $\hom(E_3^m,E_4^n) \neq 0$ and  $\hom(E_4^n,E_3^m[1]) \neq 0$ for $m\geq 1$, hence $t \leq \phi(E_4^n) \leq t+1$.
The rest of the proof is the same as the proof of Corollary \ref{coro for E_2 and E_1}, but  one must use   the short exact sequence \eqref{ses3} instead of  \eqref{ses1}.
\epr
The short exact sequences with middle terms  $E_2^0, E_4^0, M'$ are unique:
\begin{lemma} \label{ses for M' etc} If $0\rightarrow A \rightarrow C \rightarrow B \rightarrow 0$ is a short exact sequence in
$\mc A$ with $A\neq 0$ and $B\neq 0$,  then we have the following implications:
\begin{itemize}
    \item if $C \cong E_2^0$, then $A\cong E_3^0$ and $B\cong M$;
    \item if $C \cong E_4^0$, then $A\cong M$ and $B\cong  E_1^0$;
    \item if $C \cong M'$, then $A\cong E_3^0$ and $B\cong  E_1^0$.
\end{itemize}
\end{lemma}
\bpr See the representations $E_1^0, E_2^0, E_3^0,  E_4^0, M, M'$ in Proposition \ref{exceptional objects in Q1}.\epr

\subsubsection{Comments on \textbf{C1} objects}

 Recall(see Lemma \ref{lemma for C1, C2b}) that for  any \textbf{C1} object $R\in \mc A$  there exists  an exceptional pair $(X,Y)$  in $\mc A$  satisfying
  $\bd R & \rDotsto^{\textbf{C1}} & (X,Y) \ed$, $\hom(X,Y)=0$, $\hom(R,X)\neq 0$, $\hom(Y,R)\neq 0$.  A list of the exceptional pairs in $\mc A$ is given in
  Lemma \ref{coro Q1 exceptional pairs}. Using table \eqref{Q1 table} we see that the exceptional pairs $(X,Y)$ in $\mc A$ with  $\hom(X,Y)=0$ are
 \be \label{some pairs ass to C1}  (E_1^0,E_2^0), (E_1^0,E_3^0), (E_4^0,E_3^0), (E_1^m,M), (M,E_3^m),(M',E_2^m),(E_4^m,M')\  \qquad  m\in \NN. \ee

By setting $R$ to specific objects in $\mc A_{exc}$ we can shorten this list further as follows:

\begin{lemma} \label{comments on C1} Let $R \in \{E_i^m: m\in \NN, 1\leq i \leq 4\}$ and let  $R$ be a \textbf{C1} object.
Then there exists  a pair $(X,Y)\in P_R$ which satisfies $\bd R & \rDotsto^{\textbf{C1}} & (X,Y) \ed$, where $P_R$ is a
set of pairs depending on $R$ as shown in the table: \be\label{C1
table} \begin{array}{|c|c|} \hline R & P_R \\ \hline
  E_1^m, m \geq 1 & \{ (E_1^0,E_2^0) , \ (E_4^0,E_3^0), \ (E_1^0,E_3^0) \} \cup \{ (E_4^n,M'):   n<m \} \\   \hline
  E_2^m, m \geq 0 & \{ (E_1^0,E_2^0) , \ (E_4^0,E_3^0), \ (E_1^0,E_3^0) \} \cup\{ (M,E_3^n):  n\leq m \} \\ \hline
   E_3^m, m \geq 1&  \{ (E_1^0,E_2^0) , \ (E_4^0,E_3^0), \ (E_1^0,E_3^0) \}\cup \{ (M',E_2^n):  n <  m\} \\ \hline
     E_4^m, m \geq 0 & \{ (E_1^0,E_2^0) , \ (E_4^0,E_3^0), \ (E_1^0,E_3^0)\}\cup \{(E_1^n,M): n  \leq  m \} \\ \hline \end{array}\ee
\end{lemma} \bpr We shorten the list  \eqref{some pairs ass to C1} using $\hom(R,X)\neq 0, \hom(Y,R)\neq 0$ and table
\eqref{Q1 table}.  \epr
Recall that for each \textbf{C1} object $C\in \mc A$
we have a short exact sequence $0 \rightarrow A \rightarrow C
\rightarrow B \rightarrow 0$ with $A\neq 0$, $B \neq 0$. It follows  the first part of:
\begin{lemma} \label{C1 cases for M' etc} The simple objects $E_1^0$, $E_3^0$, $M$ cannot be \textbf{C1} objects. Furthermore: \\
If $\bd E_2^0 & \rDotsto^{\textbf{C1}} & (X,Y) \ed$, then $(X,Y)\cong (M,E_3^0)$. If $\bd E_4^0 & \rDotsto^{\textbf{C1}} & (X,Y) \ed$, then $(X,Y)\cong (E_1^0,M)$. \\
If $\bd M' & \rDotsto^{\textbf{C1}} & (X,Y) \ed$, then $(X,Y)\cong
(E_1^0,E_3^0)$.
\end{lemma}
\bpr   The rest of the lemma follows from
Lemma \ref{ses for M' etc}.\epr

 \subsubsection{\texorpdfstring{$\sigma$}{\space}-exceptional triples from the low dimensional exceptional objects \texorpdfstring{$\{ E_i^0\}_{i=1}^4$, $M$, $M'$}{\space}}\mbox{}\\
We have   the following  corollaries of Lemma \ref{inequalities}
 \begin{coro} \label{ineq M} Let $\{ E_1^0$, $E_2^0$, $E_3^0$, $M \} \subset \sigma^{ss}$. If $\phi(E_2^0)>\phi(E_1^0)$ or $\phi(E_3^0)>\phi(E_1^0)$,
then there exists a $\sigma$-exceptional triple.
\end{coro}
\bpr If $\phi(E_3^0)>\phi(E_1^0)$, then by $\phi(E_3^0)\leq
\phi(E_2^0)$ (since $\hom(E_3^0,E_2^0)\neq 0$) we have
$\phi(E_2^0)>\phi(E_1^0)$. Therefore, it is enough to  construct a
$\sigma$-exceptional triple  assuming
$\phi(E_2^0)>\phi(E_1^0)$.

By $\hom(E_2^0,M)\neq 0$ we have $\phi(E_2^0) \leq \phi(M)$. If
$\phi(E_2^0) < \phi(M)$, then we obtain a $\sigma$-exceptional
triple  from   the triple $(E_1^0,E_2^0,M)$ with
$\hom(E_1^0,E_2^0)=0$ and  Lemma \ref{inequalities} \textbf{(b)}. Hence, we
reduce   to the case $ \phi(E_2^0) =  \phi(M)>\phi(E_1^0).$

Next, we consider the triple $(E_1^0,M,E_3^0)$ with
$\hom(E_1^0,M)=\hom(E_1^0,E_3^0)=\hom(M, E_3^0)=0$. By
$\hom^1(M,E_3^0)\neq 0$ it follows $\phi(M)\leq \phi(E_3^0)+1$. If
$\phi(M)<\phi(E_3^0)+1$,   then we  obtain a  $\sigma$-triple from
Lemma \ref{inequalities} (f), due to the inequalities
$\phi(E_1^0)<\phi(M)<\phi(E_3^0)+1$. Thus, it remains  to consider
the case  $ \phi(E_1^0)<\phi(E_3^0)+1=\phi(E_2^0) =  \phi(M)$. In
this case we apply Lemma \ref{inequalities} (e) to the triple
$(E_1^0,E_3^0,E_2^0)$ with $\hom(E_1^0,E_3^0)=0$ and obtain a
$\sigma$-triple.
 \epr
 \begin{coro} \label{ineq M'} Let $\{ E_1^0$, $E_4^0$, $E_3^0$, $M' \} \subset \sigma^{ss}$.  If  $\phi(E_3^0)>\phi(E_4^0)$ or $\phi(E_3^0)>\phi(E_1^0)$,
 then there exists a $\sigma$-exceptional triple.
\end{coro}
\bpr By $\hom(E_4^0, E_1^0)\neq 0$, we see that  $\phi(E_3^0)>\phi(E_1^0)$ implies $\phi(E_3^0)>\phi(E_4^0)$. Hence,
it is enough to show that the inequality $\phi(E_3^0)>\phi(E_4^0)$ induces a $\sigma$-triple.

The triple $(E_4^0,E_3^0,M')$ has $\hom(E_4^0,E_3^0)=0$ and
$\hom(E_3^0,M')\neq 0$, therefore $\phi(E_4^0)<\phi(E_3^0)\leq
\phi(M')$. By Lemma \ref{inequalities} \textbf{(b)} we reduce to the case
$\phi(E_3^0)=\phi(M')>\phi(E_4^0)$.

Now, the triple $(E_4^0,M',E_1^0)$ has $\hom(E_4^0,M')=0$,
$\hom(M',E_1^0)\neq 0$ and  $\phi(E_4^0)<\phi(M')\leq
\phi(E_1^0)$.  Therefore, by Lemma \ref{inequalities} \textbf{(b)} we  can
reduce  the phases to
$\phi(E_4^0)<\phi(E_1^0)=\phi(E_3^0)=\phi(M')$.

Due to the obtained setting of the phases and
$\hom(E_1^0,E_3^0)=0$, Lemma \ref{inequalities} \textbf{(c)}
produces a $\sigma$-triple from the exceptional triple
$(E_4^0,E_1^0,E_3^0)$. The corollary follows.
 \epr

  \subsection{On the existence of \texorpdfstring{$S_{min}, S_{max}$}{\space}} \label{S_min S_max} For the rest  of section \ref{main theo for Q_1} we assume that $\sigma \in \st(D^b(Q_1))$ does
not admit a $\sigma$-exceptional triple. Hence, Corollaries \ref{from C3 to min}, \ref{from C2 to max} imply:
\begin{coro} \label{C2,C3 are final} If $R$   is  a \textbf{C2} or a \textbf{C3} object, then the  HN filtration of  $R$ is $\mk{alg}(R)$ and $R$  is final.
  \end{coro} Moreover, by   Corollary \ref{from C3 to min}/\ref{from C2 to max}, any \textbf{C3}/\textbf{C2} object
  induces  a semistable $S_{min/max}\in {\mc A}_{exc}$ with $\phi(S_{min/max})=\phi_{min/max}$, i. e. each \textbf{C3}/\textbf{C2} object ensures that $\mc P(\phi_{min/max})\cap \mc A_{exc} \neq \emptyset$. In this subsection we generalize these implications. The main proposition here is in terms of the numbers
   $\phi_{min}$, $\phi_{max}$ defined in \eqref{phi_min,phi_max}.  The following lemma gives  some information about these numbers.
   \begin{lemma} \label{from C2 C3 to phi_max - phi_min > 1} If there exists $R\in {\mc A}_{exc}$ which is  either  \textbf{C2} or \textbf{C3} object, then  $ \phi_{max} -\phi_{min} >1.  $
\end{lemma}
\bpr We use that $R$ is final and apply Corollary \ref{coro for
final good case}. Therefore,  we have either $\bd  R &
\rDotsto^{\textbf{C2}} & (S,E[-1]) \ed$ with $\phi(S)<\phi(E[-1])$
or $\bd  R & \rDotsto^{\textbf{C3}} & (S[1],E) \ed$ with
$\phi(S[1])<\phi(E)$, where $S,E \in \sigma^{ss}\cap {\mc
A}_{exc}$.
 Hence there exist $S,E \in \sigma^{ss}\cap {\mc A}_{exc}$ with $\phi(E)>\phi(S)+1$, therefore $\phi_{max} -\phi_{min} >1$.
\epr
The main proposition of this subsection  is:
    \begin{prop} \label{from min to max} If $\phi_{max}-\phi_{min}>1$,  then  $ \mc P(\phi_{min})\cap {\mc A}_{exc}\neq \emptyset$ and  $ \mc P(\phi_{max})\cap {\mc A}_{exc}\neq \emptyset $.
  \end{prop}
In the proof of Proposition \ref{from min to max} we use Corollaries \ref{M and E_3^m ess}, \ref{M' and
E_2^m ess}, proved later   independently.
\bpr[of Proposition \ref{from min to max}] Suppose first that $\mc P(\phi_{max})\cap \mc A_{exc}=\emptyset$. It follows that there exists a
sequence $\{S_i\}_{i \in \NN} \subset \sigma^{ss}\cap {\mc
A}_{exc}$ such that
\begin{gather} \label{from min to max 1} \phi_{min}+1 < \phi(S_0)<\phi(S_1)<\dots <\phi(S_i) < \phi(S_{i+1}) < \dots < \phi_{max} \\
 \label{from min to max 2} \lim_{i\rightarrow \infty} \phi(S_i)=\phi_{max}. \end{gather}
The objects $\{S_i\}_{i \in \NN}$ are  pairwise
non-isomorphic.   Since $\phi(S_0)-1 > \phi_{min}$,  there exists
$S \in \sigma^{ss}\cap {\mc A}_{exc}$ with $\phi(S_0)-1>\phi(S)
\geq \phi_{min}$. In particular, for each $i \in \NN$ holds
$\hom^*(S_i,S)=0$. From table \eqref{Q1 table} it follows that
either $S=M$ or $S=M'$,  i. e. there can be at most two elements
in $\sigma^{ss}\cap {\mc A}_{exc}$ with phase strictly smaller
than $\phi(S_0)-1$ and  such an element exists. Whence, there
exists $S_{min} \in \sigma^{ss}\cap {\mc A}_{exc}$
 of minimal phase, i. e. $\phi(S_{min})=\phi_{min}$. Furthermore $S_{min}\in \{M,M'\}$.

\uline{If $S_{min}=M$.}  Now, due to $\hom^*(S_i,M)=0$, table \eqref{Q1 table}
 shows that $\{ S_i \}_{i\in\NN}\subset \{ E_3^m, E_4^m\}_{m
\in \NN}$.  From Remark \ref{inc dec seq} and the monotone
behavior \eqref{from min to max 1} it follows that $S_i=E_3^{m_i}$
and $ m_i<m_{i+1}$ for big enough $i\in \NN. $ Later in Corollary
\ref{M and E_3^m ess} \textbf{(a)} we show  that such a   sequence
$\{S_i\}_{i \in \NN}$ with \eqref{from min to max 2} and the
equality $\phi(M)=\phi_{min}$ imply  that all elements of $\{
E_3^j \}_{j\in \NN}$ are semistable. Therefore, from
\begin{gather} \phi(M)+1 < \phi(E_3^{m_i}) \leq
\phi(E_3^{{m_i}+1})\leq \phi(E_3^{{m_i}+2})\leq..\leq
\phi(E_3^{{m_{i+1}}-1})\leq \phi(E_3^{m_{i+1}}); \nonumber \\
\phi(E_3^{m_i})< \phi(E_3^{m_{i+1}}) \nonumber \end{gather} it
follows that for some $j\in \{m_i, m_i+1,\dots,m_{i+1}\}$  we have
$\phi(M)+1<\phi(E_3^j)<\phi(E_3^{j+1})$, hence we can apply Lemma
\ref{inequalities} \textbf{(a)} to the triple
$(M,E_3^{j},E_3^{j+1})$, which  contradicts our assumption on
$\sigma$.

\uline{If $S_{min}=M'$.} Now table \eqref{Q1 table} shows that
$\{S_i\}_{i \in \NN} \subset \{ E_1^m, E_2^m\}_{m \in \NN}$ and
Remark \ref{inc dec seq} shows
 that for big enough $i\in \NN$ we have $  S_i=E_2^{m_i}, m_i<m_{i+1} $.  By Corollary \ref{M' and E_2^m ess} (a) we obtain
  $\{ E_2^j \}_{j\in \NN} \subset \sigma^{ss}$. Now similar  arguments as in the previous case (with an exceptional  triple   $(M',E_2^j,E_2^{j+1})$ for some $j\in \NN$) lead us to a contradiction.

So far, we derived  that there exists $S_{max} \in \mc P(\phi_{max})\cap \mc A_{exc} $.  Next, suppose that
 $\mc P(\phi_{min})\cap \mc A_{exc}=\emptyset$. Then we have a sequence
 $\{S_i\}_{i \in \NN} \subset \sigma^{ss}\cap {\mc A}_{exc}$  with
 \begin{gather} \label{from min to max 11} \phi_{max}-1 >\phi(S_i) > \phi(S_{i+1}) > \phi_{min} \qquad  \lim_{i\rightarrow \infty} \phi(S_i)=\phi_{min}. \end{gather}
It is clear that  $\hom^*(S_{max},S_i)=0$ for each $i\in \NN $, hence (by table \eqref{Q1 table})
  we see that $S_{max}\in \{M,M'\}$.

\uline{If $S_{max}=M'$.} In this case from table \eqref{Q1 table}
it follows that $\{S_i\}_{i \in \NN}\subset \{ E_3^m, E_4^m\}_{m
\in \NN}$.  By Remark \ref{inc dec seq} and the monotone behavior
\eqref{from min to max 11} we can construct the sequence so that $
S_i=E_4^{m_i},   m_i<m_{i+1}$ for $i\in \NN$. Now Corollary \ref{M
and E_3^m ess} \textbf{(b)} shows that  $\{ E_4^j \}_{j\in \NN}\subset
\sigma^{ss}$. Hence,  for some $j\in \NN$ we can apply Lemma
\ref{inequalities} \textbf{(a)} to the triple $(E_4^{j+1},E_4^{j}, M')$,
which is a contradiction.

\uline{If $S_{max}=M$.} Since we have $\{\hom^*(M,S_i)=0\}_{i\in
\NN}$,  table \eqref{Q1 table} shows  that $\{S_i\}_{i\in \NN}
\subset \{ E_1^m, E_2^m\}_{m \in \NN}$.  From Remark \ref{inc dec
seq}   we get  $  S_i=E_1^{m_i}$, $m_i<m_{i+1}$ for $i\in \NN$.
Corollary \ref{M' and E_2^m ess}
 \textbf{(b)} shows that  $\{ E_1^j \}_{j\in \NN}\subset
\sigma^{ss}$, hence for some $j\in \NN$ we can use Lemma
\ref{inequalities} \textbf{(a)} with the triple
$(E_1^{j+1},E_1^{j}, M)$, which gives us a  contradiction. The
proposition is proved.
 \epr

 We divide the proof of Corollaries  \ref{M and E_3^m ess}, \ref{M' and E_2^m ess} in  several lemmas.
  \begin{lemma} \label{lemma for R,E_min} Let $S_{min} \in \mc P(\phi_{min})\cap{\mc A}_{exc}$. Let $R\in {\mc A}_{exc}$ be either a \textbf{C2} object or a \textbf{C3} object. If $\hom^*(R,S_{min})=0$, then
  there exists  $S \in \sigma^{ss}\cap {\mc A}_{exc}$ with  $\hom^*(S,S_{min})=0$ and $\phi(S)+1 < \phi_{max}$.
  \end{lemma}
  \bpr Presenting  the arguments below we keep in mind Corollary \ref{C2,C3 are final}.

   If $R$ is \textbf{C2}, then we have $\mk{alg}(R)=\begin{diagram}[1em]
A_1 \oplus A_2[-1] & \rTo      &     &       &   R \\
  & \luDashto &     & \ldTo &       \\
  &           & B &       &
\end{diagram} $,  $A_2, B \in {\mc A}\setminus \{0\} $, $\phi_{max}\geq \phi(A_2)>\phi(B)+1$. From $\phi_-(R)=\phi(B)\geq \phi(S_{min})$, $\hom^*(R,S_{min})=0$,
and Lemma \ref{lemma for hom leq 1(X,S)} it follows, that
 $\hom^*(B,S_{min})$. Any $S \in Ind(B)$ satisfies  the desired
properties and the lemma follows.

If $R$ is \textbf{C3}, then  $\mk{alg}(R)=\begin{diagram}[1em]
A & \rTo      &     &       &   R \\
  & \luDashto &     & \ldTo &       \\
  &           & B[1] &       &
\end{diagram} $,  $A, B \in {\mc A}\setminus \{0\} $, $\phi_{max}\geq \phi(A)>\phi(B)+1\geq \phi(S_{min})+1$, hence $\hom^*(A,S_{min})=0$, which, together
with $\hom^*(R,S_{min})=0$, implies $\hom^*(B,S_{min})=0$. Now the
lemma follows with any $S \in Ind(B)$.
  \epr

  \begin{lemma}\label{lemma S_max,R}  Let $S_{max} \in {\mc A}_{exc}$ satisfy $\phi(S_{max})= \phi_{max}$, and let $R\in {\mc A}_{exc}$ be either
  a \textbf{C2} or a \textbf{C3} object. If $\hom^*(S_{max},R)=0$, then there exists  $S \in \sigma^{ss}\cap {\mc A}_{exc}$ with
  $\hom^*(S_{max},S)=0$ and $\phi(S) > \phi_{min}+1$.
  \end{lemma}
  \bpr If $R$ is \textbf{C2}, then we can write $\mk{alg}(R)=\begin{diagram}[1em]
A_1 \oplus A_2[-1] & \rTo      &     &       &   R \\
  & \luDashto &     & \ldTo &       \\
  &           & B &       &
\end{diagram} $,  $A_2, B \in {\mc A}\setminus \{0\} $, $\phi_{max}=\phi(S_{max})\geq \phi(A_2)>\phi(B)+1 \geq \phi_{min}+1$. Hence
$\hom^*(S_{max},B)=0$, which, together with $\hom^*(S_{max},R)=0$,
implies $\hom^*(S_{max},A_2)=0$. Now the lemma follows with $S \in
Ind(A_2)$.

If $R$ is \textbf{C3}, then
$\mk{alg}(R)=\begin{diagram}[1em]
A & \rTo      &     &       &   R \\
  & \luDashto &     & \ldTo &       \\
  &           & B[1] &       &
\end{diagram} $,  $A, B \in {\mc A}\setminus \{0\} $, $\phi(S_{max})\geq \phi(A)>\phi(B)+1\geq \phi_{min}+1$, hence
 $\hom^*(S_{max},B)=\hom^*(S_{max},A)=0$. Now any $S \in Ind(A)$ has  the desired properties.
  \epr

   \begin{lemma} \label{M and E_3^m} Let  $M\in\mc P(\phi_{min})$ or $M'\in\mc P(\phi_{max})$. If  for
   some $m>0$ we have  $E_3^m \in \sigma^{ss}$ or $E_4^m \in \sigma^{ss}$, then  there is not a  \textbf{C1} object in the set $\{ E_3^j$, $E_4^j  \}_{j \in \NN}$.
  \end{lemma}
  \bpr
 Suppose that some $R\in \{ E_3^j$, $E_4^j \}_{ j \in \NN}$ is a \textbf{C1} object. From Lemma \ref{when C1 is final} and $\hom^*(E_{3/4}^j,M)=\hom^*(M',E_{3/4}^j)=0$
  for each $j\in \NN$ we see that $R$ must be final,\footnote{Recall that we have Corollary \ref{C2,C3 are final} at our disposal, due to our assumption on $\sigma$.} hence $\mk{alg}(R)$ is the HN filtration of $R$. In particular, from
  $\bd R & \rDotsto^{\textbf{C1}} & (X,Y)\ed $ it follows that $X, Y$ are semistable and $\phi(Y)>\phi(X)$. Now  Lemma \ref{comments on C1}
  (look at the last two rows in the  table) contradicts the following negations:\footnote{For a statement $p$, when we write $\neg p$ we mean: ``p is not true''. }

$\neg$ $\left ( E_3^0,E_1^0 \in \sigma^{ss}\ \  \mbox{and} \ \ \phi(E_3^0)>\phi(E_1^0) \right )$. \  \textit{Proof:} If $E_3^0,E_1^0 \in \sigma^{ss}$, then
from   $E_3^m \in \sigma^{ss}$ or $E_4^m \in \sigma^{ss}$,$m>0$ and $\hom(E_3^0,E_{3/4}^m)\neq 0$, $\hom(E_{3/4}^m,E_1^0)\neq 0$ it follows
$\phi(E_3^0)\leq \phi(E_1^0)$.

$\neg$ $( E_3^0,E_4^0 \in \sigma^{ss}\ \  \mbox{and} \ \ \phi(E_3^0)>\phi(E_4^0))$.  \  \textit{Proof:} We are given $m>0$ with  $E_3^m \in \sigma^{ss}$ or $E_4^m \in \sigma^{ss}$,
 hence    $\hom(E_3^0,E_{3/4}^m)\neq 0$,
 $\hom(E_{3/4}^m,E_4^0)\neq 0$ imply $\phi(E_3^0)\leq \phi(E_4^0))$.

$\neg$ $\left ( E_2^0,E_1^0 \in \sigma^{ss} \ \  \mbox{and} \ \ \phi(E_2^0)>\phi(E_1^0) \right ) $.   \  \textit{Proof:} If $\phi(M)=\phi_{min}$,
then from $\hom(E_2^0,M)\neq 0$
it follows $\phi(E_2^0)=\phi(M)=\phi_{min}\leq \phi(E_1^0)$. If $\phi(M')=\phi_{max}$, then  $\hom(M',E_1^0)\neq 0$ implies
$\phi(E_2^0)\leq \phi(M')=\phi_{max}=\phi(E_1^0)$.

$\neg$ $(E_2^n,M' \in \sigma^{ss}\ \  \mbox{and} \ \ \phi(E_2^n)>\phi(M') )$.  \  \textit{Proof:} If $\phi(M)=\phi_{min}$, then by $\hom(E_2^n,M)\neq 0$ we
get $\phi(E_2^n)=\phi_{min}$. If $\phi(M')=\phi_{max}$, then $\phi(E_2^n)\leq\phi(M')$.

$\neg$ $( M,E_1^n \in \sigma^{ss} \ \  \mbox{and} \ \ \phi(M)>\phi(E_1^n) )$. \  \textit{Proof:} If  $\phi(M)=\phi_{min}$, then from $E_1^n \in \sigma^{ss}$ it follows
$\phi(M)\leq \phi(E_1^n)$.  If  $\phi(M')=\phi_{max}$, then  $\hom(M',E_1^n)\neq 0$ implies $\phi(M)\leq \phi(M')=\phi_{max}=\phi(E_1^n)$.

 The lemma follows.
 \epr
 \begin{coro} \label{M and E_3^m ess} Let $\{S_i\}_{i \in \NN} $ be a sequence of
  pairwise non-isomorphic, semistable objects with $\{S_i\}_{i \in \NN} \subset \{ E_3^j$, $E_4^j \}_{ j \in \NN}$. If any of the two conditions below is satisfied
 \begin{itemize}
    \item[\textbf{(a)}] $M\in \mc P(\phi_{min})$,  $\lim_{i\rightarrow \infty} \phi(S_i)=\phi_{max}$,
    \item[\textbf{(b)}] $M'\in \mc P(\phi_{max})$,  $\lim_{i\rightarrow \infty} \phi(S_i)=\phi_{min}$,
\end{itemize}
then  all the exceptional objects in the set $\{ E_3^j$, $E_4^j \}_{ j \in \NN}$ are semistable.
 \end{coro}
 \bpr Since each  $E \in \{ E_3^j$, $E_4^j \}_{j \in \NN}$ is a trivially coupling object,  it is neither \textbf{B1} nor \textbf{B2} (Corollary
 \ref{from good to good case}).   From  Lemma \ref{M and E_3^m}  we know that  any $ E $ is either semistable or \textbf{Ci}(i=2,3). However, if it is \textbf{Ci}(i=2,3), then:

 \textbf{(a)} By $\hom^*(E,M)=0$ (see table \eqref{Q1 table}), $\phi(M)=\phi_{min}$, and Lemma \ref{lemma for R,E_min} there exists $S \in \sigma^{ss}\cap {\mc A}_{exc}$
 with $\hom^*(S,M)=0$ and $\phi(M)+1\leq \phi(S)+1<\phi_{max}$, which by $\lim_{i\rightarrow \infty} \phi(S_i)=\phi_{max}$ implies that $(M,S,S_i)$ is an exceptional triple for big enough $i$. By Corollary \ref{coro for isom triples} this cannot happen, since $\{S_i\}_{i \in \NN}$ are pairwise non-isomorphic.

 \textbf{(b)} By $\hom^*(M',E)=0$(see table \eqref{Q1 table}), $\phi(M')=\phi_{max}$, and Lemma \ref{lemma S_max,R} there exists
 $S \in \sigma^{ss}\cap {\mc A}_{exc}$ with $\hom^*(M',S)=0$ and $\phi(M')\geq \phi(S)>\phi_{min}+1$, which by
 $\lim_{i\rightarrow \infty} \phi(S_i)=\phi_{min}$ implies that $(S_i,S,M')$ is an exceptional triple for big enough $i$.
 This contradicts  Corollary \ref{coro for isom triples}.
 \epr
The arguments for the proof  of Corollary \ref{M' and E_2^m ess} are  the same, but the role of Lemma  \ref{M and E_3^m} is played by the following Lemma \ref{M' and E_2^m}.
  \begin{lemma} \label{M' and E_2^m} Let  $M'\in \mc P(\phi_{min})$ or $M\in \mc P(\phi_{max})$.
  If for some $m>0$ we have $E_1^m \in \sigma^{ss}$ or $E_2^m \in \sigma^{ss}$, then  there is not a \textbf{C1}
  object in the set $\{ E_1^j$, $E_2^j \}_{j \in \NN}$.
  \end{lemma}
  \bpr Using that for each $j\in \NN$ we have $\hom^*(E_{1/2}^j,M')=0$, $\hom^*(M,E_{1/2}^j)=0$ and   Lemma \ref{comments on C1} (this time the first two rows in the  table) by the same arguments as in Lemma \ref{M and E_3^m} we reduce the proof to the negations:

$\neg$ $\left ( E_2^0,E_1^0 \in \sigma^{ss} \ \  \mbox{and} \ \
\phi(E_2^0)>\phi(E_1^0) \right ) $. \ \ \textit{Proof:} If
$E_2^0,E_1^0 \in \sigma^{ss}$, then by $\hom(E_2^0,E_{1/2}^m)\neq
0$, $\hom(E_{1/2}^m,E_1^0)\neq 0$, $m>0$ it follows that
$\phi(E_2^0)\leq \phi(E_1^0)$.

$\neg$ $\left ( E_3^0,E_1^0 \in \sigma^{ss}\ \  \mbox{and} \ \
\phi(E_3^0)>\phi(E_1^0) \right )$. \ \ \textit{Proof:} Follows
from   $\hom(E_3^0,E_{1/2}^m)\neq 0$ and \\
$\hom(E_{1/2}^m,E_1^0)\neq 0$.

$\neg$ $( E_3^0,E_4^0 \in \sigma^{ss}\ \  \mbox{and} \ \
\phi(E_3^0)>\phi(E_4^0))$. \ \   \textit{Proof:} If
$\phi(M')=\phi_{min}$, then  from  $\hom(E_3^0,M')\neq 0$ it
follows $\phi_{min}=\phi(E_3^0)
 \leq \phi(E_4^0)$. If $\phi(M)=\phi_{max}$, then  $\hom(M,E_4^0)\neq 0$ implies $\phi_{max}=\phi(E_4^0)\geq \phi(E_3^0)$.

 $\neg$ $( M,E_3^n \in \sigma^{ss} \ \  \mbox{and} \ \ \phi(E_3^n)>\phi(M) )$. \ \   \textit{Proof:} If $\phi(M')=\phi_{min}$, then  we use
 $\hom(E_3^n,M')\neq 0$. If $\phi(M)=\phi_{max}$, then $E_3^n \in \sigma^{ss}$
 implies $\phi(E_3^n)\leq \phi(M)$.

$\neg$ $(E_4^n,M' \in \sigma^{ss}\ \  \mbox{and} \ \ \phi(M')>\phi(E_4^n) )$. \ \   \textit{Proof:} If $\phi(M')=\phi_{min}$, then $E_4^n \in \sigma^{ss}$ implies $ \phi(M')\leq
\phi(E_4^n)$. If $\phi(M)=\phi_{max}$, then  the negation follows from  $\hom(M,E_4^n)\neq 0$.

 The lemma follows.
 \epr
\begin{coro} \label{M' and E_2^m ess} Let $\{S_i\}_{i \in \NN} \subset \{ E_1^j$, $E_2^j \}_{j \in \NN}$ be a sequence of  pairwise non-isomorphic, semistable objects. Any of the following two settings:
 \begin{itemize}
    \item[\textbf{(a)}]   $M'\in \mc P(\phi_{min})$,  $\lim_{i\rightarrow \infty} \phi(S_i)=\phi_{max}$,
    \item[\textbf{(b)}] $M\in \mc P(\phi_{max})$,  $\lim_{i\rightarrow \infty} \phi(S_i)=\phi_{min}$,
\end{itemize}
implies that  $\{ E_1^j$, $E_2^j \}_{j \in \NN} \subset \sigma^{ss}$.
 \end{coro} \bpr      The arguments are   the same as  those used  in the proof of Corollary \ref{M and E_3^m ess}, but we use  Lemma \ref{M' and E_2^m} instead of Lemma  \ref{M and E_3^m}.
 \epr

Note  that  the conclusions of   Corollaries  \ref{M' and E_2^m ess} and \ref{M and E_3^m ess}, namely that  $\{E_2^m,E_1^m\}\subset \sigma^{ss}$  and  $\{E_3^m,E_4^m\}\subset \sigma^{ss}$, are components   of the data in the  two non-locally finite cases   \eqref{non-generic}.    In the next  subsection  we derive \eqref{non-generic}  from the assumption $\phi_{max}-\phi_{min}>1$, and Corollaries  \ref{M' and E_2^m ess}, \ref{M and E_3^m ess} will be helpful at some points.

The implications given below  are further  minor steps towards   derivation of the non-locally finite cases  \eqref{non-generic}.  These implications will be used in both Subsection \ref{>1} and Subsection \ref{subsect leq11}.

 \begin{lemma}\label{stabilize the phases} \mbox{}

\textbf{(a)} If $\phi_{max}=\phi(M')$ and $\{E_4^m: m \in \NN\}
\subset \sigma^{ss} $,  then  $\{E_4^m : m \in \NN \}\subset \mc
P(t) $ for some $t\leq \phi_{max}$.

 \textbf{(b)} If $\phi_{max}=\phi(M)$ and $\{E_1^m : m \in \NN\}\subset \sigma^{ss} $,  then  $\{E_1^m : m \in \NN \}\subset \mc P(t) $ for
  some $t\leq \phi_{max}$.

   \textbf{(c)} If $\phi_{min}=\phi(M')$ and $\{E_2^m : m \in \NN\}\subset \sigma^{ss} $,  then  $\{E_2^m : m \in \NN\}\subset \mc P(t) $ for
    some $t\leq \phi_{max}$.

   \textbf{(d)} If $\phi_{min}=\phi(M)$ and  $\{E_3^m : m \in \NN\}\subset \sigma^{ss} $,  then  $\{E_3^m : m \in \NN\}\subset \mc P(t) $ for
    some $t\leq \phi_{max}$.

 \end{lemma}
 \bpr Presenting  the proof  we keep in mind  Remark \ref{inc dec seq}:

 \textbf{(a)}For any $m\in \NN $ we have  $\phi(E_4^{m+1})\leq \phi(E_4^{m}) \leq \phi(M')$.  The triple  $(E_4^{m+1},E_4^m, M')$ has
  $\hom(E_4^m,M')=0$. Hence from Lemma \ref{inequalities} \textbf{(c)} it follows    $\phi(E_4^{m+1})=\phi(E_4^m)$ for each $m \in \NN$.

 \textbf{(b)}We apply  the same arguments as in \textbf{(a)} to the triple $(E_1^{n+1},E_1^{n},M)$ with $\hom(E_1^{n},M)=0$.

 \textbf{(c)} Now $ \phi(M')\leq  \phi(E_2^{n}) \leq \phi(E_2^{n+1})$, $\hom(M',E_2^n)=0$ and we can apply Lemma
 \ref{inequalities} \textbf{(b)}  to the triple $( M',E_2^{n},E_2^{n+1})$, which implies  $\phi(E_2^n)=\phi(E_2^{n+1})$ for each $n\geq 0$.

 \textbf{(d)} We apply  the same arguments as in \textbf{(c)} to the triple $( M,E_3^{n},E_3^{n+1})$ with $\hom(M,E_3^n)=0$.
 \epr

\subsection{The case \texorpdfstring{$\phi_{max}-\phi_{min} > 1$}{\space}} \label{>1}

 In this subsection we show  that the inequality $\phi_{max}-\phi_{min} > 1$  is inconsistent with the assumption that there is not a $\sigma$-exceptional triple.
 The inequality  $\phi_{max}-\phi_{min} > 1$ implies by  Proposition  \ref{from min to max} that (for brevity we denote this product by $\Phi$):
  \begin{gather} \Phi=(\mc P(\phi_{min})\cap {\mc A}_{exc}) \times ( \mc P(\phi_{max})\cap {\mc A}_{exc})    \neq \emptyset. \end{gather}
  If  $(S_{min},S_{max}) \in \Phi$, then $(S_{min},S_{max})$ is an exceptional pair, since $\phi_{max}-\phi_{min} > 1$.  Hence there exists  unique
   $E\in {\mc A}_{exc}$, s. t. $(S_{min},E, S_{max})$ is an exceptional triple. It is  very important for us that \uline{ $E$ must be  necessarily semistable},
    which follows from  \ref{between min and max}.

 F the rest of this subsection we assume that $\phi_{max}-\phi_{min} > 1$. In the end we conclude that $\Phi \neq \emptyset$ contradicts the
 non-existence of a $\sigma$-exceptional triple.

  Since  any $(S_{min},S_{max})\in \Phi$ is an exceptional pair in $\mc A$, it must be some of the pairs listed in Corollary \ref{coro Q1 exceptional pairs}.
  We show case-by-case (in a properly chosen order) that for each pair $(A,B)$ in this list the incidence  $(A,B)\in \Phi$ leads to a contradiction.
  We show first that $(E_1^0,E_3^0)\not \in \Phi$.
\begin{lemma} \label{E_1^0,E_3^0 not in Phi} $(E_1^0,E_3^0)\not \in \Phi$.  \end{lemma}
\bpr Suppose that  $(E_1^0,E_3^0)\in \Phi $. We consider the
triple $(E_1^{0},M,E_3^0)$. From Proposition \ref{between min and
max} it follows that $M\in \sigma^{ss}$, hence
$\phi_{min}=\phi(E_1^{0})\leq \phi(M) \leq
\phi(E_3^0)=\phi_{max}$. One of these inequalities must be proper.
However, by $\hom(E_1^{0},M)=\hom(M,E_3^{0})=0$ and Lemma
\ref{inequalities} \textbf{(b)}, \textbf{(c)} we obtain a
$\sigma$-exceptional triple, which is a contradiction. \epr

We introduce  the following   formal rules, which facilitate the exposition:
 \begin{gather}\label{formal 1}  \bd (A,C)\in \Phi & \rImplies^{(A,B,C)} & \ \mbox{either} \ (B,C)\in \Phi \ \ \mbox{or} \ \ (A,B)\in \Phi  \ed \\
\label{formal 2} \bd (A,C)\in\Phi & \rImplies^{(A,B,C),\ \hom(A,B)=0 } & (A,B)\in \Phi \ed \\
\label{formal 3}  \bd (A,C)\in\Phi & \rImplies^{(A,B,C),\
\hom(B,C)=0 } &  (B,C)\in \Phi \ed. \end{gather} In  \eqref{formal
1}, \eqref{formal 2}, and \eqref{formal 3} the triple $(A,B,C)$ is
the unique exceptional  triple(taken from Lemma \ref{exceptional
colleections}) with first element $A$ and last element $C$.  In
all the three rules  we implicitly   use Proposition \ref{between
min and max}, from which it follows  $B\in\sigma^{ss}$,  and hence
$\phi_{min}=\phi(A) \leq \phi(B)\leq \phi(C)=\phi_{max}$.   The specific arguments
assigned to each  individual rule are:
\begin{itemize}
    \item[\eqref{formal 1}] from Lemma \ref{inequalities} \textbf{(a)} and $\phi_{max}-\phi_{min}>1$ it follows that  either $\phi(A)=\phi(B)=\phi_{min}$
    or $\phi(B)=\phi(C)=\phi_{max}$, whence we reduce to either $(B,C)\in \Phi $ or $(A,B) \in \Phi $;
    \item[\eqref{formal 2}] by Lemma \ref{inequalities} \textbf{(b)} and $\hom(A,B)=0$  we get $\phi(B)=\phi(C)=\phi_{max}$, whence  $(A,B)\in \Phi$;
    \item[\eqref{formal 3}] by Lemma \ref{inequalities} \textbf{(c)} and $\hom(B,C)=0$ we get $\phi(A)=\phi(B)=\phi_{min}$, whence  $(B,C) \in \Phi $.
\end{itemize}
Now we eliminate some  pairs $(X,Y)$ by showing that $(X,Y)\in \Phi$ implies $(E_1^0,E_3^0)\in \Phi$.
\begin{coro} \label{E_4^0,E_3^0 etc}For each $n\in \NN$  any of the pairs $(E_4^0,E_3^0), (E_1^0,E_2^0), (M,E_3^n), (E_4^n,M'), (E_1^n,M)$,
 $(M',E_2^n)$, $(E_1^{n+1},E_4^{n})$, $(E_4^{n},E_1^{n})$, $(E_4^{n+1},E_4^{n})$, $(E_1^{n+1},E_1^{n})$, $(E_2^{n},E_3^{n+1})$, $(E_3^{n},E_2^{n})$,
  $(E_3^{n},E_3^{n+1})$, \\ $(E_2^{n},E_2^{n+1})$ is not in $\Phi$.
\end{coro}
\bpr We keep in mind the formal rules \eqref{formal 1}, \eqref{formal 2}, \eqref{formal 3}. The following expressions and  Lemma \ref{E_1^0,E_3^0 not in Phi} show that each of the listed
pairs is not in  $\Phi$.

$(E_4^0,E_3^0)\in \Phi$ $\bd & \rImplies^{(E_4^{0},E_1^0,E_3^0),\ \hom(E_1^{0},E_3^{0})=0 } &  (E_1^{0},E_3^0)\in\Phi \ed $.

$(E_1^0,E_2^0)\in\Phi$ $\bd & \rImplies^{(E_1^0,E_3^{0},E_2^0),\ \hom(E_1^{0},E_3^{0})=0 } &  (E_1^{0},E_3^0)\in\Phi \ed $.

$(M,E_3^0)\in\Phi$ $ \bd & \rImplies^{(M,E_4^0,E_3^0),\ \hom(E_4^0,E_3^0)=0 } &  (E_4^{0},E_3^0) \in \Phi \ed $.

$(M,E_3^n)\in \Phi$,$n\geq 1$ $ \bd & \rImplies^{(M,E_3^{n-1},E_3^n),\ \hom(M,E_3^{n-1})=0 } &  (M,E_3^{n-1})\in\Phi \ed $    $\bd & \rImplies^{\mbox{  induction} } & \ed$  $(M,E_3^0)$.

$(E_4^0,M')\in\Phi$ $ \bd & \rImplies^{(E_4^0,E_3^0,M'),\ \hom(E_4^0,E_3^0)=0 } &(E_4^{0},E_3^0) \in\Phi\ed $.

$(E_4^n,M')\in \Phi$,$n\geq 1$   $ \bd & \rImplies^{(E_4^n,E_4^{n-1},M'),\ \hom(E_4^{n-1},M')=0 } &  (E_4^{n-1},M')\in\Phi \ed $ $\bd & \rImplies^{\mbox{  induction} } & \ed$ $(E_4^0,M')$.

$(E_1^0,M)\in\Phi$ $ \bd & \rImplies^{(E_1^0,E_2^0,M),\ \hom(E_1^0,E_2^0)=0 } & (E_1^{0},E_2^0)\in\Phi \ed $.

$(E_1^n,M)\in\Phi$,$n\geq 1$   $ \bd & \rImplies^{(E_1^n,E_1^{n-1},M),\ \hom(E_1^{n-1},M)=0 } &  (E_1^{n-1},M)\in\Phi \ed $  $ \bd & \rImplies^{\mbox{  induction} } & \ed   $ $(E_1^0,M)$.

$(M',E_2^0)\in\Phi$ $ \bd & \rImplies^{(M',E_1^0,E_2^0),\ \hom(E_1^0,E_2^0)=0 } & (E_1^{0},E_2^0)\in\Phi \ed $.

$(M',E_2^n)\in\Phi$,$n\geq 1$ $ \bd & \rImplies^{(M',E_2^{n-1},E_2^{n}),\ \hom(M',E_2^{n-1})=0 } &  (M',E_2^{n-1})\in\Phi \ed $  $ \bd & \rImplies^{\mbox{  induction} } & \ed   $  $(M',E_2^0)$.

$(E_1^{n+1},E_4^{n})\in \Phi$,$n\geq 0$ $ \bd & \rImplies^{(E_1^{n+1},M,E_4^n),\ \hom(E_1^{n+1},M)=0 } &  (E_1^{n+1},M) \in \Phi\ed $.

$(E_4^{n},E_1^{n})\in \Phi$,$n\geq 0$ $ \bd & \rImplies^{(E_4^{n},M',E_1^{n}),\ \hom(E_4^{n},M')=0 } &  (E_4^{n},M')\in \Phi \ed $.

$(E_4^{n+1},E_4^{n})\in \Phi$,$n\geq 0$ $\bd & \rImplies^{(E_4^{n+1},E_1^{n+1},E_4^{n})} &  \ed $   either $(E_1^{n+1},E_4^{n})\in \Phi $ or $(E_4^{n+1},E_1^{n+1})\in \Phi $.

$(E_1^{n+1},E_1^{n})\in \Phi$,$n\geq 0$ $\bd & \rImplies^{(E_1^{n+1},E_4^{n},E_1^{n})} &  \ed $   either $(E_4^{n},E_1^{n})\in \Phi $ or $(E_1^{n+1},E_4^{n})\in \Phi $.

$(E_2^{n},E_3^{n+1})\in \Phi$,$n\geq 0$ $ \bd & \rImplies^{(E_2^{n},M,E_3^{n+1}),\ \hom(M,E_3^{n+1})=0 } &  (M,E_3^{n+1})\in \Phi \ed $.

$(E_3^{n},E_2^{n})\in \Phi$,$n\geq 0$ $ \bd & \rImplies^{(E_3^{n},M',E_2^{n}),\ \hom(M',E_2^{n})=0 } &  (M',E_2^{n})\in \Phi \ed $.

$(E_3^{n},E_3^{n+1})\in \Phi$,$n\geq 0$ $\bd & \rImplies^{(E_3^{n},E_2^{n},E_3^{n+1})} &  \ed $   either $(E_2^{n},E_3^{n+1})\in \Phi $ or $(E_3^{n},E_2^{n})\in \Phi $.

$(E_2^{n},E_2^{n+1})\in \Phi$,$n\geq 0$ $\bd &
\rImplies^{(E_2^{n},E_3^{n+1},E_2^{n+1})} &  \ed $   either
$(E_3^{n+1},E_2^{n+1})\in \Phi $ or $(E_2^{n},E_3^{n+1})\in \Phi
$. \epr  We eliminated many pairs by using only Section \ref{two examples},  Proposition \ref{between min and max}, and Lemma \ref{inequalities}.  It remains to consider the incidences:
$(M,E_4^n)$, $ (E_3^n, M')$, $ (M',E_1^n)$, $(E_2^n,M) \in \Phi $ for  $n\geq 0$. From any of these incidences,  with the help of Corollaries \ref{M and
E_3^m ess}, \ref{M' and E_2^m ess} and Lemma \ref{stabilize the phases},    we  will   derive some of  the non-locally finite cases  \eqref{non-generic}, which is excluded by Corollaries \ref{coro for E_3
and E_4}, \ref{coro for E_2 and E_1}. We start with $(M,E_4^n)$.

\begin{lemma}For each $n\geq 0$ we have $(M,E_4^n)\not \in\Phi$. \end{lemma}
\bpr Suppose that $(M,E_4^n)\in\Phi$. In the previous corollary  we showed that  $(E_4^{n+1},E_4^{n})\not \in \Phi $. Now from the implication
$(M,E_4^n)\in\Phi$,$n\geq 0$ $ \bd & \rImplies^{(M,E_4^{n+1},E_4^{n}) } &  \ed $ either $(E_4^{n+1},E_4^{n})\in \Phi $ or $(M,E_4^{n+1})\in \Phi $\\
we deduce that $(M,E_4^{n+1})\in \Phi $, and by induction we obtain $  \phi(E_4^{i})=\phi_{max}  $ for $ i \geq n $.  We are given  also $\phi(M)=\phi_{min}$,
therefore
 we can use Corollary \ref{M and E_3^m ess} \textbf{(a)} to obtain  $  \{E_3^j,E_4^j \}_{j \in \NN} \subset \sigma^{ss}. $
By Remark \ref{inc dec seq} we see
\be\label{M,E_4^n 2} \forall i \geq 0 \ \ \  \phi(E_4^{i})=\phi_{max}. \ee
The next step is to show that
\be \label{M,E_4^n 3} \forall i \geq 0 \ \ \  \phi(E_3^{i})=\phi_{max}-1. \ee
 Since $\hom^1(E_4^0,E_3^0)\neq 0$, we have $\phi(M)= \phi_{min}< \phi_{max}-1 = \phi(E_4^0)-1  \leq  \phi(E_3^0) \leq \phi_{max}$.  Whence:
$$ \phi(M)<\phi(E_3^0)\leq \phi(E_4^0)\leq \phi(E_3^0)+1.$$
If $\phi(E_4^0)<\phi(E_3^0)+1$, then we have $
\phi(M)<\phi(E_3^0)\leq \phi(E_4^0)< \phi(E_3^0)+1$ and Lemma
\ref{inequalities} \textbf{(d)} applied to  the triple
$(M,E_4^0,E_3^0)$ gives us a $\sigma$-exceptional triple.
Therefore $\phi(E_4^0)=\phi(E_3^0)+1$.  We showed above that
$E_3^j $
 is semistable for each $j\in \NN$. From Lemma \ref{stabilize the phases} \textbf{(d)} we get $\phi(E_3^0)=\phi(E_3^i)$ for any $i\geq 0$, thus we
get \eqref{M,E_4^n 3}. However \eqref{M,E_4^n 2} and \eqref{M,E_4^n 3}
contradict Corollary \ref{coro for E_3 and E_4}. \epr

\begin{lemma} For each $n\geq 0$ we have  $(E_3^n,M')\not \in\Phi$. \end{lemma}
\bpr Suppose  that $(E_3^n,M')\in\Phi$. We  obtain  a  contradiction of Corollary \ref{coro for E_3 and E_4} as follows:\\
$(E_3^n,M')\in\Phi$,$n\geq 0$ $ \bd & \rImplies^{(E_3^{n},E_3^{n+1},M') } &  \ed $ either $(E_3^{n+1},M')\in \Phi $ \ or
$(E_3^{n},E_3^{n+1})\in \Phi $ $\bd &\rImplies^{\mbox{Corollary \ref{E_4^0,E_3^0 etc}}} & \ed$ \\ $(E_3^{n+1},M')\in \Phi $ $\bd &\rImplies^{\mbox{ind. }} & \ed$ $ \forall i \geq n \ \
 \phi(E_3^{i})=\phi_{min} $ $\bd &\rImplies^{\mbox{Corollary \ref{M and E_3^m ess} \textbf{(b)} } }  & \ed$ $  \{E_3^j,E_4^j  \}_{j \in \NN} \subset \sigma^{ss}. $
By Remark \ref{inc dec seq} we see that   $\phi(E_3^{i})=\phi_{min}$ for  $  i \geq 0 $.
We show below that  $(E_3^0,M')\in\Phi$ implies that   $\phi(E_4^{i})=\phi_{min}+1 $  for each $i \geq 0 $, which contradicts
  Corollary \ref{coro for E_3 and E_4}.

Indeed, by $\hom^1(E_4^0,E_3^0)\neq 0$   we can write $ \phi_{min}=\phi(E_3^0)\leq  \phi(E_4^0)  \leq  \phi(E_3^0)+1 < \phi_{max}=\phi(M')$.
 The triple $(E_4^0, E_3^0, M')$ has $\hom(E_4^0,E_3^0)=0$, therefore  from  $\phi(E_4^0)  <  \phi(E_3^0)+1$ it follows that for some
 $j\geq 1$ the triple  $(E_4^0,E_3^0,M'[-j])$ is $\sigma$-exceptional. Therefore $\phi(E_4^0)=\phi(E_3^0)+1$.  We showed above that  $\{E_4^j \} \subset \sigma^{ss}$. By  Lemma \ref{stabilize the phases} \textbf{(a)} we conclude  that
 $\phi(E_4^n)=\phi_{min}+1$ for each $n\geq 0$.
\epr

\begin{lemma} For each $n\geq 0$ we have  $(M',E_1^n)\not \in\Phi$. \end{lemma}
\bpr Suppose that $(M',E_1^n)\in\Phi$. We show that  this  contradicts  Corollary \ref{coro for E_2 and E_1} as follows: \\
$(M',E_1^n)\in\Phi$,$n\geq 0$ $ \bd & \rImplies^{(M',E_1^{n+1},E_1^{n}) } &  \ed $ either $(E_1^{n+1},E_1^{n})\in \Phi $ \ or
$(M',E_1^{n+1})\in \Phi $ $\bd &\rImplies^{\mbox{Corollary \ref{E_4^0,E_3^0 etc}}} & \ed$ \\ $(M',E_1^{n+1})\in \Phi $ $\bd &\rImplies^{\mbox{ind. }} & \ed$ $ \forall i \geq n \ \
 \phi(E_1^{i})=\phi_{max} $ $\bd &\rImplies^{\mbox{Corollary \ref{M' and E_2^m ess} \textbf{(a)} } }  & \ed$ $  \{E_1^j,E_2^j \}_{j \in \NN} \subset \sigma^{ss}. $\\
By Remark \ref{inc dec seq} we see that
  $\phi(E_1^{i})=\phi_{max} $ for each $ i \geq 0$. Furthermore, using  $(M',E_1^0)\in\Phi$  we show below that
$ \phi(E_2^{i})=\phi_{max}-1 $ must hold  for $  i \geq 0 $, which  contradicts   Corollary \ref{coro for E_2 and E_1}.

 Indeed, it follows from  $\hom^1(E_1^0,E_2^0)\neq 0$  that  $ \phi_{min}=\phi(M')< \phi(E_1^0)-1\leq \phi(E_2^0)  \leq \phi_{max}=\phi(E_1^0)$.
 If $\phi(E_1^0)  <  \phi(E_2^0)+1$, then $\phi(M')<  \phi(E_2^0)  \leq \phi(E_1^0)<\phi(E_2^0)+1 $, and  the triple
 $( M',E_1^0, E_2^0)$ with $\hom(E_1^0,E_2^0)=0$ gives rise to a $\sigma$-triple  by Lemma \ref{inequalities} \textbf{(d)}. Therefore
 $\phi(E_1^0)=\phi(E_2^0)+1$.  Since  $\{E_2^j\}\subset \sigma^{ss}$,  Lemma \ref{stabilize the phases} \textbf{(c)}
  implies that  $\phi(E_2^n)=\phi(E_2^{0})$ for $n\geq 0$. The lemma is proved.
\epr

\begin{lemma} For each $n\geq 0$ we have  $(E_2^n,M)\not \in\Phi$. \end{lemma}
\bpr Suppose that  $(E_2^n,M)\in\Phi$. We will obtain  a  contradiction of Corollary \ref{coro for E_2 and E_1} as follows:
$(E_2^n,M)\in\Phi$,$n\geq 0$ $ \bd & \rImplies^{(E_2^{n},E_2^{n+1},M) } &  \ed $ either $(E_2^{n+1},M)\in \Phi $ \ or
$(E_2^{n},E_2^{n+1})\in \Phi $ $\bd &\rImplies^{\mbox{Corollary \ref{E_4^0,E_3^0 etc}}} & \ed$ \\ $(E_2^{n+1},M)\in \Phi $ $\bd &\rImplies^{\mbox{ind. }} & \ed$ $ \forall i \geq n \ \
   \phi(E_2^{i})=\phi_{min} $ $\bd &\rImplies^{\mbox{Corollary \ref{M' and E_2^m ess} \textbf{(b)} } }  & \ed$ $  \{E_1^j,E_2^j \}_{j \in \NN} \subset \sigma^{ss}. $\\
By Remark \ref{inc dec seq} we conclude that
  $\phi(E_2^{i})=\phi_{min} $  for $  i \geq 0$.  We show below that $(E_2^0,M)\in \Phi$  implies that $ \phi(E_1^{i})=\phi_{min}+1 $ for each $ i \geq 0 $, which   contradicts     Corollary \ref{coro for E_2 and E_1}.

Indeed, it follows from $\hom^1(E_1^0,E_2^0)\neq 0$ and $(E_1^0,M)\not \in \Phi$(see Corollary \ref{E_4^0,E_3^0 etc})  that  $ \phi_{min}=\phi(E_2^0)<  \phi(E_1^0)\leq
\phi(E_2^0) +1 < \phi_{max}=\phi(M)$.
 If $\phi(E_1^0)  <  \phi(E_2^0)+1$, then  for some $j\geq 1$   the triple  $( E_1^0, E_2^0,M[-j])$  is   $\sigma$-exceptional, since   $( E_1^0, E_2^0,M)$ is exceptional and $\hom(E_1^0,E_2^0)=0$. Therefore $\phi(E_1^0)=\phi(E_2^0)+1<\phi(M)$.   Lemma \ref{stabilize the phases} \textbf{(b)}
  gives us  $\phi(E_1^n)=\phi(E_1^{0})=\phi_{min}+1$ for $n\geq 0$.
\epr
Therefore, we reduce to $\phi_{max}-\phi_{min}\leq
1$, which will be assumed until the end of the proof.
\subsection{The case \texorpdfstring{$\phi_{max}-\phi_{min} \leq 1$}{\space}}\label{subsect leq11}
From this inequality we obtain a contradiction here again, by deriving    the non-locally finite cases  \eqref{non-generic}.  We  show first in a series of lemmas that $\mc A_{exc}\subset \mc P(\phi_{min})\cup \mc P(\phi_{max})$, $\phi_{\max}-\phi_{min}=1$.   Lemma \ref{from C2 C3 to phi_max - phi_min > 1} and Corollary \ref{lemma for final good case objects in main th for Q1} imply immediately
\begin{lemma} \label{max-min leq 1 fund} Any $E\in {\mc A}_{exc}$ is either semistable or irregular or a
final \textbf{C1} object.
\end{lemma}

 Any $X\in \{E_i^j:j\in \NN, 1\leq i\leq 4\}$ is a trivially
coupling object,
 hence by Lemma \ref{from good to good case} we have only two possibilities: $X$ is semistable or $X$ is
  a final \textbf{C1} object (cannot be irregular).
\begin{coro} \label{M is not C1}
The objects $E_1^0$, $E_3^0$ are semistable, and $M$ is either irregular or semistable.
\end{coro}
\bpr  The objects $E_1^0$, $E_3^0$, $M$ cannot be  \textbf{C1} by   Lemma \ref{C1 cases for M' etc}. \epr
\begin{lemma}The object  $E_2^0$ is semistable.
\end{lemma}
\bpr Suppose that  $E_2^0$ is not semistable.

 Therefore $E_2^0$ must be \textbf{C1}, and we have $\bd E_2^0 & \rDotsto^{\textbf{C1}} &(X,Y)\ed $ for some exceptional pair $(X,Y)$. By Lemma \ref{C1 cases for M' etc}, we see that  $(X,Y)=(M,E_3^0)$.
 Since $E_2^0$ is final, we can write  \ben M, E_3^0 \in \sigma^{ss} \qquad \phi(M) < \phi(E_3^0).\een From the triple $(E_1^0,M,E_3^0)$, which satisfies
  $\hom(E_1^0,M)=\hom(E_1^0,E_3^0)=\hom(M,E_3^0)=0$, and Lemma \ref{inequalities} (f) we see that $\phi(E_1^0)=\phi(M)+1$
  (recall that $E_1^0$ is semistable). From $\phi_{max}-\phi_{min} \leq 1$ it is clear that
\ben \phi(M)=\phi_{min}, \qquad
\phi(E_1^0)=\phi_{max}=\phi(M)+1.\een The obtained relations imply
that  $E_4^0$ is semistable. Indeed, if $E_4^0$ is not semistable,
then it must be final \textbf{C1}, hence  by  Lemma \ref{C1 cases
for M' etc} we have  $\bd E_4^0 & \rDotsto^{\textbf{C1}}
&(E_1^0,M)\ed $, which in turn implies $\phi(M)>\phi(E_1^0)$
contradicting   $\phi(E_1^0)=\phi(M)+1$. Therefore $E_4^0$ is
semistable.  Now consider the triple $(M,E_4^0,E_3^0)$ with
$\hom(E_4^0,E_3^0)=0$. We have  $\phi(M)\leq \phi(E_4^0)\leq
\phi(M)+1$, $\phi(M) < \phi(E_3^0)\leq \phi(M)+1$. If
$\phi(M)<\phi(E_4^0)$, then $\phi(M)-1<\phi(E_4^0[-1])\leq \phi(M)
$, $\phi(M)-1<\phi(E_3^0[-1])\leq \phi(M) $ and
$(M,E_4^0[-1],E_3^0[-1])$ is a $\sigma$-exceptional triple. So
far,  assuming  that $E_2^0$ is not semistable we get: \ben
\phi_{min}=\phi(M)=\phi(E_4^0)<\phi(E_3^0)\leq
\phi(E_1^0)=\phi(E_4^0)+1. \een Therefore  $\phi(E_4^0)-1 <
\phi(E_3^0[-1])\leq  \phi(E_1^0[-1])=\phi(E_4^0)$ and then the
triple $(E_4^0, E_1^0[-1], E_3^0[-1])$ is  a $\sigma$-exceptional
triple (since $\hom(E_1^0,E_3^0)=0$). This triple contradicts our
assumption on $\sigma$. \epr
\begin{lemma} The object $E_4^0$ is semistable.
\end{lemma}
\bpr Suppose that $E_4^0$ is not semistable. Hence it is final \textbf{C1}, and by Lemma \ref{C1 cases for M' etc} we have\\
$\bd E_4^0 & \rDotsto^{\textbf{C1}} &(E_1^0,M)\ed $. Since $E_4^0$
is final, it follows:
 \ben M, E_1^0 \in \sigma^{ss} \qquad \phi(M) > \phi(E_1^0).\een
 From the simple objects triple $(E_1^0,M,E_3^0)$    and Lemma \ref{inequalities} \textbf{(f)}
 it follows that $\phi(M)=\phi(E_3^0)+1$ (recall that $E_3^0$ is semistable),  hence by  $\phi_{max}-\phi_{min} \leq 1$:
\ben \phi(E_3^0)=\phi_{min}, \qquad \phi(M)=\phi_{max}=\phi_{min}+1.\een
 Now we have   $\{E_1^0,E_2^0,E_3^0, M \} \subset \sigma^{ss}$.
 From Corollary \ref{ineq M} it follows $\phi(E_2^0)\leq \phi(E_1^0)$. Whence, we have $\phi_{min}=\phi(M)-1\leq \phi(E_2^0)\leq \phi(E_1^0)<\phi(M)$.
The triple $(E_1^0,E_2^0,M)$ has  $\hom(E_1^0,E_2^0)=0$. We rewrite the last inequalities as follows
 $\phi(M[-1])\leq \phi(E_2^0)\leq \phi(E_1^0)<\phi(M[-1])+1$ and obtain a $\sigma$-exceptional triple $(E_1^0,E_2^0,M[-1])$, which is a contradiction.
 The lemma follows.
 \epr
Now, using that $\{E_i^0\}_{i=1}^4\subset \sigma^{ss}$, we show
that  $M$, $M'$ cannot be irregular.
\begin{coro} \label{no B2} There does not exist a \textbf{B2} object.
\end{coro}
\bpr Suppose that $E\in \mc A$ is a \textbf{B2} object. Since the
only Ext-nontrivial couple is  $\{M,M'\}$,  we have $E\in
\{M,M'\}$ and  we can write \be  \mk{alg}(E)=\begin{diagram}[1em]
A & \rTo      &     &       &   E \\
  & \luDashto &     & \ldTo &       \\
  &           & B[1] &       &
\end{diagram} \ \ \qquad\begin{array}{c} \{E,\Gamma \}=\{ M,M'\} \ \mbox{for some} \  \Gamma \in Ind(B), \\ \phi_-(A)>\phi(B)+1=\phi(\Gamma)+1.  \end{array}  \ee
 From
$\hom(M,E_4^0)\neq 0$, $\hom(M',E_1^0)\neq 0$, and $\{E_1^0,
E_4^0\}\subset\sigma^{ss}$ (shown in the  preceding lemmas)  it
follows that there exists $X \in \sigma^{ss}\cap {\mc A}_{exc}$
with $\hom(E,X)\neq 0$, hence $\hom(A,X)\neq 0$ and $\phi_-(A)\leq
\phi(X)$. Whence, we obtain $\phi(X)\geq \phi_-(A) >
\phi(\Gamma)+1$ with $X,\Gamma \in \sigma^{ss}\cap {\mc A}_{exc}$,
which  contradicts  the inequality $\phi_{max}-\phi_{min}\leq 1$.
\epr
\begin{lemma}\label{no bad} There does not exist a $\sigma$-irregular object.
\end{lemma}
\bpr  By Corollary \ref{no B2} we have to show that neither $M$ nor $M'$ can be \textbf{B1}.

Suppose that  $E \in \{M,M'\}$  is \textbf{B1}, then we can write:
\begin{gather}  \mk{alg}(E)=\begin{diagram}[1em]
A_1 \oplus A_2[-1] & \rTo      &     &       &   E \\
  & \luDashto &     & \ldTo &       \\
  &           & B &       &
\end{diagram} \ \ \qquad \begin{array}{c} \{E,\Gamma \}=\{ M,M'\} \ \mbox{for some} \    \Gamma \in Ind(A_2), \\  \hom^1(A_1,A_1)=\hom^1(A_2,A_2)=0 \\
 \phi_-(A_1\oplus A_2[-1])\geq \phi(B)=\phi_-(E).  \end{array} \nonumber  \end{gather}
 We show first that  each $Y\in Ind(A_2)$ must be semistable with $\phi(Y)=\phi(B)+1$, which implies
\be \label{no bad 1} A_2 \in \sigma^{ss}, \ \
\phi(A_2[-1])=\phi(B).\ee To that end we observe  that there
exists $X \in \sigma^{ss}\cap {\mc A}_{exc}$ with  $\hom(X,B)\neq
0$, and hence \be \label{no bad 2} \phi(X)\leq \phi(B) \ \ \ \ X
\in \sigma^{ss}\cap {\mc A}_{exc} .\ee Indeed, if we find $X \in
\sigma^{ss}\cap {\mc A}_{exc}$ with  $\hom^*(X,E)=0$ and
$\hom(X,\Gamma)\neq 0$, then from the triangle $\mk{alg}(E)$ it
follows that $\hom(X,B)\cong \hom(X,A_1[1]\oplus A_2)\neq 0$ (the
latter  does not vanish by $\Gamma \in Ind(A_2)$ and
$\hom(X,\Gamma)\neq 0$). Looking at table \eqref{Q1 table} we see
that  $\hom^*(E_2^0,M')=0, \hom(E_2^0,M)\neq 0$,
$\hom^*(E_3^0,M)=0,    \hom(E_3^0,M')\neq 0$,  therefore
\be \label{no bad 3} \begin{array}{c c c} X=E_2^0 & \mbox{if} & E=M'\\
X=E_3^0 & \mbox{if} & E=M.
        \end{array} \ee
        Let us take  any  $Y\in Ind(A_2)$.  From Lemma \ref{no two bad cases} \textbf{(c)} it follows that $Y$ cannot be $\sigma$-irregular. Hence it is either semistable or a
        final \textbf{C1} object.
         If $Y$ is \textbf{C1}, then   $\bd Y & \rDotsto^{\textbf{C1}} & (Z,W) \ed$ for some
        $Z,W \in \sigma^{ss}\cap {\mc A}_{exc} $, and we can write $\phi(W)>\phi(Z)=\phi_-(Y) \geq \phi_-(A_2)\geq \phi(B)+1\geq  \phi(X)+1$, which contradicts $\phi_{max}-\phi_{min}\leq 1$.
        If $Y$ is semistable, then by $\phi_{max}-\phi_{min}\leq 1$ it follows $\phi(Y)\leq \phi(X)+1$, which, together with   $\phi(Y)\geq \phi(B)+1 \geq \phi(X)+1$, implies  $\phi(Y)=\phi(B)+1 = \phi(X)+1$. Whence,  we proved \eqref{no bad 1}. Furthermore, we see that \eqref{no bad 2} must be equality.

        Being a \textbf{B1} object,  $E$ is not  semistable. From the triangle $\mk{alg}(E)$, the equality $\phi(A_2[-1])=\phi(B)$
        and the fact that $\mc P(t)$ is an extension closed subcategory of $\mc T$ it follows that $A_1\neq 0$ and $\phi_+(A_1)>\phi(B)$.
        From \textbf{B1.1}  we know that    $A_1$ is a proper $\mc A$-subobject of $E$. Since $M$ is simple in $\mc A$, it follows
        that $E$ cannot be $M$. Whence, $E$ must be $M'$ and then $X=E_2^0$ (see \eqref{no bad 3}). The only proper
         subobject of $M'$ in $\mc A$ up to isomorphism is  $E_3^0$ and we know that it is semistable. Whence,  we arrive at
         $\phi(E_3^0)>\phi(B)=\phi(X)=\phi(E_2^0)$, It follows that $\hom(E_3^0,E_2^0)=0$, which contradicts table \eqref{Q1 table}.
\epr
 \begin{coro} \label{M,M' are semistable}The objects $M$, $M'$ are   semistable and
 \be \label{the inequalities in leq 1}  \phi(E_2^0)\leq \phi(E_1^0) \qquad \phi(E_3^0)\leq \phi(E_1^0)  \qquad \phi(E_3^0)\leq \phi(E_4^0).\ee
 \end{coro}
 \bpr The semistability of  $M$   follows from Corollary \ref{M is not C1} and Lemma \ref{no bad}. Then from Corollary \ref{ineq M}
 we get $\phi(E_2^0)\leq \phi(E_1^0)$ and   $\phi(E_3^0)\leq \phi(E_1^0)$.

 Using $\phi_{max}-\phi_{min}\leq 1$ we showed so far that the cases \textbf{C2, C3, B1, B2 } can not appear.  Therefore we have only two options for $M'$: either semistable or final \textbf{C1}.

 Suppose that $M'$ is  final \textbf{C1}.   Lemma \ref{C1 cases for M' etc} implies that $\bd M' & \rDotsto^{\textbf{C1}} &(E_1^0,E_3^0)\ed $. Therefore
 $\phi(E_3^0)>\phi(E_1^0)= \phi_-(M').$
 However we showed already that  $ \phi(E_3^0)\leq \phi(E_1^0)$. Hence, $M'$ must be also semistable. Now Corollary \ref{ineq M'}
implies  $ \phi(E_3^0)\leq  \phi(E_4^0)$.
 \epr
So far, we showed that the low dimensional exceptional objects  $\{E_i^0\}_{i=1}^4$, $M$, $M'$ are semistable. The following implications, due  to table  \eqref{C1 table} in Lemma \ref{comments on C1},  will help us to show that $\mc A_{exc}\subset \sigma^{ss}$.
 \begin{coro} \label{if Ci ineq} Let $R\in \{E_i^m: m \geq  1, 1\leq i \leq 4  \}$ and let $R$ be non-semistable.

\begin{itemize}
    \item[\textbf{(a)}] If $R=E_1^m$, then $\phi(E_4^n)<\phi(M')$ for some $n<m$, and  hence $\phi(M)<\phi(M')$
    \item[\textbf{(b)}] If $R=E_2^m$, then $\phi(M)<\phi(E_3^n)$ for some $n\leq m$,  and  hence $\phi(M)<\phi(M')$
    \item[\textbf{(c)}] If $R=E_3^m$, then $\phi(M')<\phi(E_2^n)$ for some $n <  m$, and  hence $\phi(M')<\phi(M)$
    \item[\textbf{(d)}] If $R=E_4^m$, then  $\phi(E_1^n)<\phi(M)$ for some $n \leq  m$, and hence $\phi(M')<\phi(M)$.
\end{itemize}

 \end{coro}
 \bpr Now we have $M$, $M'\in \sigma^{ss}$ and   any non-semistable $R\in {\mc A}_{exc}$ is a final \textbf{C1} object. Note also that for each $n\in \NN$ we have  $\hom(M,E_4^n)\neq 0$, $\hom(E_3^n,M')\neq 0$, $\hom(E_2^n,M)\neq 0$,  $\hom(M',E_1^n)\neq 0$,  which implies
 $\phi(M)\leq \phi_+(E_4^n)$,  $\phi_-(E_3^n)\leq \phi(M')$, $\phi_-(E_2^n)\leq \phi(M)$, $\phi(M')\leq \phi_+(E_1^n)$. Due to  Lemma \ref{coro for final good case} and the inequalities \eqref{the inequalities in leq 1} in Corollory \ref{M,M' are semistable}, we can   remove the pairs $(E_1^0,E_2^0)$, $(E_4^0,E_3^0)$, $(E_1^0,E_3^0)$ from table
  \eqref{C1 table} in Lemma \ref{comments on C1}.    If   $E_i^m\not \in \sigma^{ss}$ for some $m \geq  1$, $1\leq i \leq 4$, then  $E_i^m$ is a final \textbf{C1} object  and  the corollary follows from  table \eqref{C1 table} in Lemma \ref{comments on C1}.  \epr
 Knowing that the triple  $(E_1^0,M,E_3^0)$ of the simple objects is semistable, we  obtain that one of  three equalities below must hold, which implies $\phi_{max}-\phi_{min}=1$.
 \begin{lemma} \label{phi_m-phi_min=1}  There is an equality $\phi_{max}-\phi_{min}=1$.   One of the following  equalities must hold:
\be \label{three poss for all stable}  \phi(E_1^0)=\phi(M)+1, \ \ \  \phi(E_1^0) =  \phi(E_3^0)+1,\ \ \  \phi(M) =  \phi(E_3^0)+1.\ee
 \end{lemma}
\bpr From $\hom(E_1^0,M[1])\neq 0$, $\hom(E_1^0,E_3^0[1])\neq 0$,
$\hom(M,E_3^0[1])\neq 0$ we have $\phi(E_1^0)\leq \phi(M)+1$,
$\phi(E_1^0)\leq \phi(E_3^0)+1$, $\phi(M)\leq \phi(E_3^0)+1$.
Applying \textbf{(f)} of  Lemma  \ref{inequalities} to the  triple
$(E_1^0, M, E_3^0)$, we see that one of the equalities
\eqref{three poss for all stable} holds. Hence $
\phi_{max}-\phi_{min} \geq 1 $  and the lemma follows. \epr
\begin{coro} \label{min max of M,M' 0}   $\phi(M)\in \{ \phi_{min}, \phi_{max}\}$.
\end{coro}
 \bpr    Suppose that  $\phi_{min}< \phi(M)<\phi_{max}$. By Lemma \ref{phi_m-phi_min=1}  we get  $ \phi_{min}=\phi(E_3^0)$ and $ \phi(E_3^0)+1=\phi(E_1^0)=\phi_{max}. $
 Therefore,  we can write $\phi(E_3^0)< \phi(M)< \phi(E_3^0)+1$ and $\phi(E_3^0)\leq \phi(E_2^0)\leq \phi(E_3^0)+1$. Now by
  combining  \eqref{phase formula}  and the equality $ Z(E_2^0)= Z(E_3^0)+Z(M)$ (see Lemma \eqref{exceptional objects in Q1})  we  obtain:
\be  \phi_{min}=\phi(E_3^0)< \phi(E_2^0)< \phi(E_3^0)+1=\phi(E_1^0)=\phi_{max}. \ee

By semistability of $M'$ we have either  $\phi_{min}=\phi(E_3^0)< \phi(M')$ or  $\phi_{min}=\phi(E_3^0)=\phi(M')$. We  aim at a contradiction\footnote{of the assumption that there is not a $\sigma$-exceptional triple}  by using  either the triple $(E_3^0,M',E_2^0)$  with  $\hom(M', E_2^0)=0$ or the triple $(M',E_1^0,E_2^0)$ with
 $\hom(E_1^0, E_2^0)=0$.  If $\phi_{min}=\phi(E_3^0)< \phi(M')$, then we have  $\phi(E_3^0)< \phi(M')\leq \phi(E_3^0)+1$,
 $\phi(E_3^0)< \phi(E_2^0)< \phi(E_3^0)+1$, hence the triple $(E_3^0,M'[-1],E_2^0[-1])$ is  $\sigma$-exceptional.
 If $\phi_{min}=\phi(E_3^0) =\phi(M')$, then we have  $\phi(E_1^0)= \phi(M')+1$,  $\phi(M')< \phi(E_2^0)< \phi(M')+1$,  hence
 the triple $(M',E_1^0[-1],E_2^0[-1])$ is  $\sigma$-exceptional.
\epr
\begin{coro}  \label{min max of M,M'}   We have  $ \{\phi(M),\phi(M'), \phi(E_j^0) \}\subset \{ \phi_{min}, \phi_{max}\}$ for  $j =1,2,3,4$.
\end{coro}
\bpr Now we have  $ \{\phi(M),\phi(M'), \phi(E_j^0) \}\subset \sigma^{ss}$ and   $\phi(M)\in \{ \phi_{min}, \phi_{max}\}$.
It is enough to show $\{\phi(E_1^0), \phi(E_3^0)\} \subset\{ \phi_{min}, \phi_{max} \}$, because  then by formula  \eqref{phase formula}, the equalities  $Z(M')=Z(E_1^0)+Z(E_3^0)$, $Z(E_2^0)=Z(M)+Z(E_3^0)$, $Z(E_4^0)=Z(M)+Z(E_1^0)$, and the inequalities
  $\phi_{min}\leq \phi(M'), \phi(E_2^0), \phi(E_4^0) \leq \phi_{max}$ it follows that  $\{ \phi(M'), \phi(E_2^0), \phi(E_4^0)\}\subset  \{ \phi_{min}, \phi_{max}\}$.

If $\phi(M)=\phi_{min}$, then by $\hom(E_2^0,M)\neq 0$ it follows that $\phi(E_2^0)=\phi_{min}$, and
by Lemma \ref{phi_m-phi_min=1} it follows that  $\phi(E_1^0) = \phi_{max}$. Expanding the equality  $Z(E_2^0)=Z(M)+Z(E_3^0)$  by  formula \eqref{phase formula}, and using
  $\phi(M)=\phi(E_2^0)=\phi_{min}$,    $\phi_{min}\leq \phi(E_3^0) \leq \phi_{max}$,  we conclude $\phi(E_3^0) \in \{ \phi_{min}, \phi_{max}\}$.

If $\phi(M)=\phi_{max}$, then by $\hom(M,E_4^0)\neq 0$ it follows $\phi(E_4^0)=\phi_{max}$, and  by Lemma \ref{phi_m-phi_min=1} it follows $\phi(E_3^0) = \phi_{min}$.
 Finally,   $\phi(E_1^0) \in \{ \phi_{min}, \phi_{max}\}$  follows from  $\phi(M)=\phi(E_4^0)=\phi_{max}$, $Z(E_4^0)=Z(M)+Z(E_1^0)$, and
formula \eqref{phase formula}.
The corollary is proved.
\epr
 The proofs of semistability for  $E_1^m$ and  $ E_2^m$  share some steps because   the non-semistability of  any of them implies $\phi(M)<\phi(M')$(Corollary \ref{if Ci ineq} \textbf{(a)}, \textbf{(b)}). Similarly,  the starting argument in the proof of Lemma \ref{E_3^m,E_4^m semistable} is that  the  non-semistability of   $E_3^m$ or $E_4^m $ implies $\phi(M')<\phi(M)$.

\begin{lemma} \label{E_1^m,E_2^m semistable} All objects  in   $\{E_1^m, E_2^m  \}_{m \in \NN}$ are semistable.
\end{lemma}
 \bpr Suppose that $E_1^m$ is not semistable for some $m \in \NN$. Corollary \ref{if Ci ineq} \textbf{(a)}  shows that $E_4^n \in \sigma^{ss}$, $\phi(E_4^n)<\phi(M')$  for some
 $n\in \NN$, and $\phi(M)<\phi(M')$.    The latter inequality implies, due to  Corollary \ref{if Ci ineq} \textbf{(c)} and \textbf{(d)},   that $\{E_4^m, E_3^m\}_{m \in \NN}\subset \sigma^{ss}$, and, due to  Corollary \ref{min max of M,M'}, it implies
 \be \label{E_1^m,E_2^m semistable 1}  \phi_{min}=\phi(M), \ \ \ \phi(M')=\phi_{max}=\phi_{min}+1.\ee By  Lemma \ref{stabilize the phases} \textbf{(a)}  we can write
 $\phi(E_4^0)=\phi(E_4^n)<\phi(M')$ and combining with  Corollary  \ref{M,M' are semistable} we arrive at $\phi_{min}=\phi(M')-1 \leq \phi(E_3^0)\leq \phi(E_4^0)<\phi(M')$,
  hence the triple $(E_4^0,E_3^0,M'[-1])$ with $\hom(E_4^0,E_3^0)=0$ is a $\sigma$-exceptional triple. Therefore $\{E_1^m\}_{m \in \NN}\subset \sigma^{ss}$.

Next, suppose that    $E_2^m$ is not semistable  for some $m \in
\NN$. Then  by Corollary \ref{if Ci ineq} \textbf{(b)} we have $E_3^n \in
\sigma^{ss}$, $\phi(M)<\phi(E_3^n)$ for some
 $n\in \NN $,  and   $\phi(M)<\phi(M')$. Now by the same arguments as above we get
 \eqref{E_1^m,E_2^m semistable 1} and $\{E_4^m, E_3^m\}_{m \in \NN} \subset \sigma^{ss}$.
  By Lemma \ref{stabilize the phases} \textbf{(d)} we can write  $\phi(E_3^0)=\phi(E_3^n)>\phi(M)$.
  Combining with Corollary  \ref{M,M' are semistable} we arrive at $\phi_{min}=\phi(M) < \phi(E_3^0)\leq \phi(E_4^0)\leq \phi(M)+1$.
  These  inequalities and the exceptional triple  $(M,E_4^0,E_3^0)$ with $\hom(E_4^0,E_3^0)=0$ provide  a  $\sigma $-exceptional triple $(M,E_4^0[-1],E_3^0[-1])$.
   The lemma follows.
   \epr
 \begin{lemma} \label{E_3^m,E_4^m semistable} All objects  in   $\{E_3^m, E_4^m  \}_{m \in \NN}$ are semistable.
\end{lemma}
\bpr Suppose  that  $E_4^m$ or $E_3^m$ is not semistable for some
$m\in \NN$. By Lemma \ref{if Ci ineq} we get $\phi(M')<\phi(M)$.
Since $\{\phi(M), \phi(M')\}\subset \{\phi_{min},\phi_{max}\}$
(Corollary \ref{min max of M,M'}), we  find that:
 \ben \label{E_3^m,E_4^m semistable 1}  \phi_{min}=\phi(M'), \ \ \ \phi(M)=\phi_{max}=\phi_{min}+1.\een
We have also $\{E_1^m,E_2^m\}_{m \in \NN} \subset \sigma^{ss}$.
Thus,  \textbf{(b)} and \textbf{(c)} in Lemma \ref{stabilize the
phases} can be used to obtain:  \be \label{E_3^m,E_4^m semistable
2} \forall m \in \NN \ \ \ \ \ \ \ \ \ \ \
\phi(E_1^m)=\phi(E_1^0), \ \ \phi(E_2^m)=\phi(E_2^0).\ee From
$\hom(M,E_4^0) \neq 0$, and $\hom(E_4^0, E_1^0)\neq 0$ (note that
$\hom(E_4^0, E_1^m)= 0$ for $m\geq 1$) it follows
$\phi(M)=\phi_{max}=\phi(E_1^0)$. On the other hand, from the
triple $(M',E_1^0,E_2^0)$ with $\hom(E_1^0,E_2^0)=0$ it follows
that $\phi(E_2^0)=\phi(M')=\phi_{min}$ (otherwise
$(M',E_1^0[-1],E_2^0[-1])$ would be a $\sigma$-exceptional
triple).  Using \eqref{E_3^m,E_4^m semistable 2} we obtain \ben
\forall m \in \NN   \ \ \ \ \ \ \ \ \ \  \phi_{max}=\phi(M)
=\phi(E_1^m), \ \ \phi_{min}=\phi(M')=\phi(E_2^m).\een However,
due to \textbf{(c)} and \textbf{(d)} in Corollary \ref{if Ci
ineq}, these equalities contradict the assumption that  $E_3^m$ or
$E_4^m$ is not semistable for some $m$. The lemma follows. \epr
 \begin{coro}\label{all semistable} All exceptional objects are semistable and their  phases are in  $\{\phi_{min},\phi_{max}\}$.
 \end{coro}
 \bpr  We have already proved that the exceptional objects are semistable. As in subsection \ref{two limit points}, we denote $\delta_Z=Z(M)+Z(E_1^0)+ Z(E_3^0)$.
By  Bridgleand's  axiom \eqref{phase formula} we can rewrite  \eqref{Z(E_i^m)} as follows:
\begin{gather}  \label{all ex ar stable 6}  r(E_j^m) \exp(\ri \pi \phi(E_j^m))=m \delta_Z+ r(E_j^0) \exp(\ri \pi \phi(E_j^0))  \qquad m \in \NN,  j=1,2,3,4.
\end{gather}
In Corollary \ref{min max of M,M'} we have  $\{ \phi(M),
\phi(E_1^0),  \phi(E_3^0) \} \subset  \{ \phi_{min},
\phi_{max}\}$, therefore we can write $ \delta_Z = \Delta \exp(\ri
\pi \gamma) $ with   $\Delta \geq 0$ and $\gamma \in \{
\phi_{min}, \phi_{max}\}$.

Now \eqref{all ex ar stable 6} restricts all the phases in the set
$\{\phi_{min},\phi_{max}\}$, since $\phi_{min} \leq \phi(E_j^m)
\leq \phi_{max}=\phi_{min}+1$ for any $1\leq j \leq 4$, $m\in \NN
$ .
 \epr
We are already close to  \eqref{non-generic}. To derive completely some of the non-locally finite cases in \eqref{non-generic} we consider each   of the three equalities \eqref{three poss for all
stable}. We showed that one of them  holds.
\subsubsection{\textbf{If \texorpdfstring{$\phi(E_1^0)=\phi(M)+1$}{\space}}} Then $\phi_{min}=\phi(M)$ and $\phi_{max}=\phi(E_1^0)$.

Since $\hom(E_2^m,M)\neq 0$, we have   $\phi(E_2^m) \leq \phi(M)=\phi_{min}$ for  $m\in \NN$. Hence
$ \{ E_2^{m} \} \subset \mc P(\phi_{min})$.  We will show below that $\{E_1^m\}\subset \mc P(\phi_{max})$ and so we obtain
the first case in \eqref{non-generic}.

The sequence $\{\phi(E_1^m)\}_{m\in \NN}$ is non-increasing (see
Remark \ref{inc dec seq}) and has at most  two values. The first
value is $\phi(E_1^0)=\phi_{max}=\phi(M)+1$. Suppose that
$\phi(E_1^l)=\phi(M)$ for some $l>0$. We can assume that $l$ is
minimal, so $\phi(E_1^{l-1})=\phi(M)+1$. In table \eqref{Q1 table}
we see that $\hom(M',E_1^l)\neq 0$,  hence  $ \phi(M') \leq
\phi(M)=\phi_{min}$, i. e. $ \phi(M') = \phi(M)=\phi_{min} . $ We
have the triple $(E_1^l, M, E_4^{l-1})$ with $\hom(E_1^l, M)=0$
and  $\phi(E_1^{l})=\phi(M)$.  It follows that
$\phi(E_4^{l-1})=\phi(M)$, otherwise      Lemma \ref{inequalities}
\textbf{(b)}
  produces a $\sigma$-triple. However, now the exceptional triple  $(E_4^{l-1}, M', E_1^{l-1})$ with $\hom(E_4^{l-1},M')=0$ satisfies    $\phi(M)=\phi(E_4^{l-1})=  \phi(M')<\phi(E_1^{l-1})=\phi(M)+1$ and    Lemma \ref{inequalities} \textbf{(b)} gives a contradiction.

Whence, the equality
$\phi(E_1^0)=\phi(M)+1$ implies the first case in \eqref{non-generic}, which contradicts Corollary \ref{coro for E_2 and E_1}. Therefore, for the rest of the proof we can use  the strict  inequality: \be \label{all
are stable eq 5} \phi(E_1^0)<\phi(M)+1.\ee
\subsubsection{\textbf{\textbf{If $\phi(E_1^0)=\phi(E_3^0)+1$ or $\phi(M)=\phi(E_3^0)+1$.}}}  In both cases $\phi_{min}=\phi(E_3^0)$, $\phi_{max}=\phi(E_3^0)+1$.

We note first that  $\hom(M,E_4^m)\neq 0$ and $\hom(E_4^m,E_1^m)\neq 0$   for each integer $m$, hence \be \label{phi(M)leq phi(E_4^m)leq...} \phi(M)\leq \phi(E_4^m) \leq \phi(E_1^m)\leq \phi_{max}=\phi(E_3^0)+1  \qquad m \in \NN.\ee
 Threfore,  it is enough to consider the case $\phi(E_1^0)=\phi(E_3^0)+1$.
The latter equality  and  \eqref{all are stable eq 5}  imply  $\phi_{max}=\phi(E_1^0)$ and  $\phi(E_3^0)=\phi_{min} < \phi(M) $. It follows that
$ \phi(M)= \phi_{max}. $  Now  \eqref{phi(M)leq phi(E_4^m)leq...} implies $\{E_4^m\}_{m\in \NN}\subset \mc P(\phi_{max})$. We will show that $\{E_3^m\}_{m\in \NN}\subset \mc P(\phi_{min})$ and so we obtain the second case in \eqref{non-generic}.

Now we have $\phi(E_3^0)=\phi_{min}$. Suppose that $\phi(E_3^l)=\phi_{max}$ for some $l>0$. Choosing the minimal $l$ with this property, we have $\phi(E_3^{l-1})=\phi_{min}$. By $\hom(E_3^l, M')\neq 0$ we get $\phi(M')=\phi_{max}=\phi(M)$. It follows that $\phi(E_2^{l-1})=\phi_{min}$, because otherwise $(E_3^{l-1},M'[-1],E_2^{l-1}[-1])$ is a $\sigma$-triple, due to
$\hom(M',E_2^{l-1})=0$. However, now   $(E_2^{l-1},M[-1],E_3^{l}[-1])$ is a  $\sigma$-exceptional triple,
due to $\hom(M,E_3^{l})=0$.

Whence, any of the equalities
$\phi(E_1^0)=\phi(E_3^0)+1$ and $\phi(M)=\phi(E_3^0)+1$ implies  \eqref{non-generic}, which is the desired contradiction. Theorem \ref{main theorem for Q_1} is proved.

\appendix

\section{} \label{table with matrices}
In the table below we present the dimensions of some vector spaces
of matrices. We skip the computations.  For $m,n \geq 1$ we denote by ${\mc
M}_{k}(m,n)$ the vector space of $m\times n$ matrices over the
field $k$. The notations $\pi_\pm^m, j_\pm^m$ for $m\in \NN$ are explained before Proposition \ref{exceptional objects in Q1}.
\small
\begin{gather} \label{vect space table} \begin{array}{| c | c | c |}
  \hline
                &   V             &   \dim_k(V)     \\ \hline
  1\leq n < m        & \left \{(X,Y)\in {\mc M}_{k}(n+1,m+1)\times {\mc M}_{k}(n,m) : X\circ j_+^m=j_+^n \circ Y , \ \ X\circ j_-^m=j_-^n \circ Y   \right  \}  & 0      \\ \hline
   1\leq m \leq n        &  \line(1,0){6} \  \line(0,1){0.3} \  \line(0,1){0.3} \ \line(1,0){6}   & 1+n-m      \\ \hline
   1\leq m < n        & \left \{(X,Y)\in {\mc M}_{k}(n,m)\times {\mc M}_{k}(n+1,m+1) : X\circ \pi_+^m=\pi_+^n \circ Y, \ \   X\circ \pi_-^m=\pi_-^n \circ Y  \right  \}  & 0      \\ \hline
    1\leq n \leq m        &  \line(1,0){6} \ \line(0,1){0.3} \ \line(0,1){0.3} \ \line(1,0){6}   & 1+m-n      \\ \hline
     1\leq m, 1 \leq n        & \left \{(X,Y)\in {\mc M}_{k}(n+1,m)\times {\mc M}_{k}(n,m+1): X\circ \pi_+^m=j_+^n \circ Y,  \ \ X\circ \pi_-^m=j_-^n \circ Y   \right  \}  & 0      \\ \hline
       1\leq m, 1 \leq n        & \left \{X\in {\mc M}_{k}(n,m)  : j_+^n \circ X\circ \pi_-^m=j_-^n \circ X \circ \pi_+^m    \right  \}  & 0      \\ \hline
         1\leq m\leq n        & \left \{X\in {\mc M}_{k}(n+1,m)  : \pi_-^n \circ X\circ \pi_+^m=\pi_+^n \circ X\circ \pi_-^m   \right  \}  & 0      \\ \hline
            0\leq n < m        &  \line(1,0){6} \ \line(0,1){0.3} \ \line(0,1){0.3} \ \line(1,0){6}   & m-n      \\ \hline
            1\leq n\leq m        & \left \{X\in {\mc M}_{k}(n,m+1)  : j_-^n \circ X\circ j_+^m=j_+^n \circ X\circ j_-^m   \right  \}  & 0      \\ \hline
             0\leq m < n        &  \line(1,0){6} \ \line(0,1){0.3} \ \line(0,1){0.3} \ \line(1,0){6}   & n-m      \\ \hline
                  1\leq m, 1 \leq n        & \left \{(X,Y)\in {\mc M}_{k}(m,n+1)\times {\mc M}_{k}(m+1,n) : X\circ j_+^n=\pi_+^m \circ Y,  \ \ X\circ j_-^n=\pi_-^m \circ Y   \right  \}  & m+n      \\ \hline
                    0\leq n < m        & \left \{(X,Y)\in {\mc M}_{k}(n+1,m)^2 : \pi_+^n\circ X =\pi_-^n \circ Y,  \ \ X\circ \pi_+^m=Y\circ \pi_-^m    \right  \}  & m-n-1      \\ \hline
                     0\leq m < n        & \left \{(X,Y)\in {\mc M}_{k}(n,m+1)^2 : j_+^n\circ X =j_-^n \circ Y,  \ \ X\circ j_+^m=Y\circ j_-^m    \right  \}  & n-m-1      \\ \hline
                       0\leq m,0 \leq n        & \left \{(X,Y)\in {\mc M}_{k}(n+1,m+1)^2 : \pi_+^n\circ X =\pi_-^n \circ Y,  \ \ X\circ j_+^m=Y\circ j_-^m    \right  \}  & n+m+2      \\ \hline
   0\leq m, 0 \leq n        & \left \{X\in {\mc M}_{k}(n+1,m+1) : \pi_-^n \circ X\circ j_+^m=\pi_+^n \circ X \circ j_-^m    \right  \}  & m+n+1      \\ \hline
        \end{array}
    \end{gather}
\normalsize

\section{The Kronecker quiver} \label{Kronecker}

\subsection{There are not Ext-nontrivial couples in \texorpdfstring{$Rep_k(K(l))$}{\space}}  \label{there are no Ext-nontrivial...}

The quiver  with two vertices and  $l\geq 2$ parallel  arrows  will be denoted by  $K(l)$. Here we  revisit    \cite[Lemma 4.1]{Macri}. This  lemma implies the title of this subsection.

 Following the notations of \cite{Macri}, let $s_0$ and $s_1$ be  the exceptional objects in $D^b(K(l))$, such that $s_0[1]$ is the simple representation with $k$ at the source, and $s_1$ is the simple representation with $k$ at the sink, and then define  $s_i$   for each $i\in \ZZ$ as follows:
 \begin{gather} s_{-i}=L_{s_{-i+1}}(s_{-i+2}) , \ \  s_{i+1}=R_{s_{i}}(s_{i-1})  \qquad i\geq 1. \end{gather}
 The Braid group $B_2$ is isomorphic to $\ZZ$. By the transitivity of the action of $B_2$ on the set of  full exceptional collections, shown in   \cite{WCB1}, it follows that, up to shifts, the complete list of  the exceptional  pairs in $Rep_k(K(l))$ is $\{(s_i,s_{i+1})\}_{i\in\ZZ}$.    \textit{Lemma 4.1 in \cite{Macri}  says that $s_{\leq 0}[1], s_{\geq 1} \in Rep_k(K(l))$, and:} \be \label{lemma 4.1 in Macri}  p\neq 0 \Rightarrow  \hom^p(s_i,s_j)= 0 ; \qquad \quad  p\neq 1 \Rightarrow  \hom^p(s_j,s_i)= 0; \qquad \qquad  i<j. \ee
  Now $\{s_{-i}[1]\}_{i\geq 0}\cup \{s_{i}\}_{i\geq 1}$ is the complete list of exceptional objects of $Rep_k(K(l))$, and from  the vanishings \eqref{lemma 4.1 in Macri} it follows  that for any couple $\{X,Y\}$ in this list  $\hom^1(X,Y)\neq 0$ implies  $\hom^1(Y,X)= 0$. Thus, there are not Ext-nontrivial couples in  $Rep_k(K(l))$.

    One can show that     the  following inequalities hold  for each $i\in \ZZ$:
 \begin{gather}\label{nonzerohominapair} l=\hom(s_i,s_{i+1})<\hom(s_i,s_{i+2})<\dots; \quad
 0= \hom^1(s_i,s_{i-1})<\hom^1(s_i,s_{i-2})<\dots, \\
 \dim_k(s_1)=\dim_k(s_0[1])<\dim_k(s_2)=\dim_k(s_{-1}[1])<\dots.\end{gather}
  which implies that $\{s_{-i}[1]\}_{i\geq 0}\cup \{s_{i}\}_{i\geq 1}$ are pairwise non-isomorphic. Whence, in this case the action of the Braid group is  free (compare with Remark \ref{braid}).

\subsection{\texorpdfstring{$\sigma$}{\space}-exceptional pairs in \texorpdfstring{$D^b(K(l))$}{\space}} \label{exceptional pairs} \mbox{}\\

The full exceptional collections in $D^b(K(l))$ have length two, so the analogue of Theorem \ref{main theorem for Q_1} is:
\begin{lemma}\label{main theorem for K(l)}
For each $\sigma \in \st(D^b(K(l)))$ there exists a
$\sigma$-exceptional pair.
\end{lemma}
The statement of    \cite[Lemma 4.2]{Macri} is equivalent to the statement of Lemma \ref{main theorem for K(l)}.
For the sake of completeness we give a proof of Lemma \ref{main theorem for K(l)} here.

Denote, for brevity $\mc A=Rep_k(K(l))$, and take any $\sigma=(\mc
P, Z) \in \st(D^b(\mc A))$.  There are not Ext-nontrivial couples
in $\mc A$ and    the exceptional pairs of $D^b(\mc A)$, up to
shifts, are  a sequence  $\{(s_i,s_{i+1})\}_{i\in \ZZ}$, where
$\{s_{-i}[1]\}_{i \geq 0} \cup \{s_{i}\}_{i \geq 1}\subset \mc
A$(see Appendix \ref{there are no Ext-nontrivial...}).   By Remark
\ref{if there are not Ext-nontrivila pairs} we reduce the proof
immediately  to the case, where all the  exceptional objects are
semistable.    In \eqref{nonzerohominapair} we have $\{
\hom(s_i,s_{i+1})\neq 0 \}_{i\in \ZZ}$, hence $\{
\phi(s_i)\leq\phi(s_{i+1})\}_{i\in \ZZ}$. If $\phi(s_i)
<\phi(s_{i+1})$ for some  $i\in \ZZ$, then there exists   $j\geq
1$ with $ \phi(s_{i+1}[-j])\leq \phi(s_i)<\phi(s_{i+1}[-j])+1$,
and hence,  due to \eqref{lemma 4.1 in Macri}, the pair
$(s_i,s_{i+1}[-j])$ is  $\sigma$-exceptional. Thus, we reduce to
the case, where all $\{ s_i \}_{i\in \ZZ}$ have the same
phase,\footnote{In the end of Appendix \ref{there are no
Ext-nontrivial...} we pointed out that  $\{s_i\}_{i\geq 1}$ are
pairwise non-isomorphic.} say $t\in \RR$:  \be \label{s_0,s_1 in
P(t)}  \{ s_i \}_{i\in \ZZ} \subset \mc P(t). \ee We  show now
that the obtained inclusion contradicts the locally finiteness of
$\sigma$, i. e. \eqref{s_0,s_1 in P(t)} is a non-locally finite
case.

Since all  the exceptional pairs in $\mc A$  are  $\{(s_{i-1}[1],s_{i}[1])\}_{i\leq -1} \cup \{ (s_0[1],s_1) \} \cup \{(s_{i},s_{i+1})\}_{i\geq 1}$, it follows from \eqref{s_0,s_1 in P(t)} that:
\be \label{cont} \mbox{\textit{For each exceptional pair $(S,E)$
with  $S,E \in \mc A$ we have $\phi(S)\geq \phi(E)$.} } \ee

 We will obtain a contradiction by constructing an exceptional pair $(S,E)$ in $\mc A$ with $\phi(S)<\phi(E)$.
  Recall that  $Z$ is the central charge of $\sigma$.  By \eqref{s_0,s_1 in P(t)} and  \eqref{phase formula}  we have  $\{ Z(s_1), Z(s_0[1])=-Z(s_0)\} \subset \RR \exp(\ri \pi t)$.  Since\footnote{This isomorphism is determined  by  assigning to $[X] \in K(D^b(\mc A))$, for $X\in \mc A $,  the dimension vector $\ul{\dim}(X)\in \ZZ^2$.} $K(D^b(\mc A))\cong \ZZ^2 $ and the simple objects $s_0[1]$, $s_1$ form  a basis of $K(D^b(\mc A))$, it follows that $\im(Z) \subset \RR \exp(\ri \pi t)$. Now using \eqref{phase formula} again, we concude that $\mc P(x)$ is trivial  for $x\in(t-1,t)$, therefore $\mc P(t-1,t]=\mc P(t)$. From the very foundation  \cite{Bridg1} given by  T. Bridgeland,    we know that $\mc P(t-1,t]$ is a heart of a  bounded  $t$-structure of $D^b(\mc A)$, so $\mc P(t)$ is a heart as well. Due to  this  property of $\mc P(t)$, it is also  well known that   $\bd K(\mc P(t)) & \rTo^{K(\mc P(t)\subset D^b(\mc A))} & K(D^b(\mc A))\ed $ is an isomorphism, so $ K(\mc P(t))\cong \ZZ^2$. The  locally finiteness of $\sigma$ implies that $\mc P(t)$ is an abelian category of finite length, which in turn, combined with  $ K(\mc P(t))\cong \ZZ^2$,  implies that $\mc P(t)$ has exactly two simple objects, say $X,Y \in \mc P(t)$. It follows by Lemma \ref{P(t) is a thick subcategory}, that $\{X,Y\}$ are indecomposable in $D^b(\mc A)$, therefore $X=X'[i]$, $Y=Y'[j] $ for some $i,j \in \ZZ$ and $X', Y' \in \mc A$.  Viewing $\mc A$ as   the extension closure of  $s_0[1]$, $s_1$, we see that  $X', Y' \in \mc A\subset \mc P[t,t+1]$. Now from $\{ X'[i],Y'[j]\} \subset \mc P(t)$ it follows  that either $\phi(X')=t$, $i=0$  or $\phi(X')=t+1$, $i=-1$, and the same holds for $Y',j$. If either $i=i'=-1$ or $i=i'=0$, then  $\hom(s_1,X)=\hom(s_1,Y)=0$ or  $\hom(X,s_0)=\hom(Y,s_0)=0$, which contradicts the
   existence  of a Jordan-H\"{o}lder filtration of $s_0,s_1 \in \mc P(t)$  via the simples $X,Y$ of $\mc P(t)$. Thus,  we arrive at:
\be \label{formula for Kronecker} X=X',  \qquad  Y=Y'[-1],\qquad X',Y'  \in \mc A\qquad \phi(X')=t, \ \phi(Y')=t+1.\ee
By $\phi(Y')>\phi(X')$ it follows $\hom(Y',X')=0$. Since $Y'[-1]$, $X'$ are non-isomorphic simple objects in the abelian category $\mc P(t)$, it follows that $\hom(Y'[-1],X')=0$ as well, hence $ \hom^*(Y',X')=0 $.

 The pair  $(X',Y')$ in $\mc A$ has $\phi(X')<\phi(Y')$ and $ \hom^*(Y',X')=0 $, and it almost contradicts  \eqref{cont}, but  we have no arguments for the vanishings ${\rm Ext}^1(X',X')=0$ and  ${\rm Ext}^1(Y',Y')=0$.

Keeping in mind the comments in the beginning of Subsection
\ref{comments on stab cond}, we can view   $\mc P(t)$ as the
extension closure in $D^b(\mc A)$ of the  set $\{Y'[-1], X'\}$.
Denoting  the extension closures of $X'$ and $Y'$  by $\mc X$ and
$\mc Y$, respectively, it is clear that\textit{ $\mc P(t)$ is the
extension closure of  $ \mc Y[-1] \cup \mc X$ and }\begin{gather}
\label{hom^*()=0}  [\mc X]=\NN [X'], \ \  [\mc Y]=\NN [Y'], \qquad
\hom^*(\mc Y,\mc X )=0, \qquad \mc X \subset \mc A\cap \mc P(t), \
\  \mc Y \subset \mc A\cap \mc P(t+1), \end{gather}  where the
first two equalities are between subsets of $K(D^b(\mc A))$. Using
$\hom^*(\mc Y,\mc X )=0$ and that $\mc P(t)$ is the extension
closure of  $ \mc Y[-1] \cup \mc X$, as in the case of
semi-orthogonal decompositions,   one can show  that for each
$X\in \mc P(t)$ there exists  a triangle $ \bd A[-1]& \rTo & X &
\rTo & B &  \rTo & A \ed$ with  $ A \in \mc Y$,$ B \in \mc X$ and
$\hom^*(A,B)=0$.  Since $s_j \in {\mc A}_{exc}\cap \mc P(t)$ for
$j\geq 1$, the corresponding triangle for $s_j$ is:
\begin{gather} \label{triangle for s_j} \bd[width=1em] s_j& \rTo & B & \rTo & A &  \rTo & s_j[1] \ed, \qquad  \hom^*(A,B)=0, \   A\in \mc Y, B \in \mc X , s_j  \in {\mc A}_{exc}. \end{gather}
To prove Lemma \ref{main theorem for K(l)}, we show first that  we
can assume   $A\neq 0$. After that we  recall some of the
arguments used  in  Subection \ref{if B_0 neq 0} for obtaining the
properties \textbf{C2.1}     in the triangle \eqref{triangle woth
f_0}.  These arguments lead to the vanishings
$\hom^1(B,B)=\hom^1(A,A)=0$.  Taking any $S \in Ind(B)$, $E \in
Ind(A)$, we obtain an exceptional pair $(S,E)$ in $\mc A$ with
$\phi(S)<\phi(E)$, which contradicts  \eqref{cont}.

  Suppose that $A=0$. Then $s_j \cong B\in \mc X$ and by \eqref{hom^*()=0} we have  $\ul{\dim}(s_j)=p \   \ul{\dim}(X')$ for some $p\in \NN$. Since  $s_j$ is  exceptional   and $X'$ is  indecomposable, then $\scal{\ul{\dim}(s_j),\ul{\dim}(s_j)}=1$(see \eqref{euler}) and  $\scal{\ul{\dim}(X'),\ul{\dim}(X')}\leq 1$ (see \cite[p. 58]{Kac}).\footnote{where $\scal{,}$ is the Euler form of $K(l)$.} It follows that $\ul{\dim}(s_j)=\ul{\dim}(X')$, $\scal{\ul{\dim}(X'),\ul{\dim}(X')}= 1$. Recall that $X'$ is simple in $\mc P(t)$, which implies $\hom(X',X')=1$. Now formula \eqref{euler} shows  that  $X'$ is an exceptional object, and hence  $\ul{\dim}(X')=\ul{\dim}(s_j)$ implies  that $X'\cong s_j$.\footnote{  There is at most one representation without self-extensions of a given dimension vector(\cite[p. 13]{WCB2}).}  Thus, $A=0$ implies $X'\cong s_j$.
   It follows,   since $\{s_i\}_{i\geq 1}$ are pairwise non-isomorphic,
  that  in \eqref{triangle for s_j} the object $A$ can vanish for at most  one integer $j\geq 1$.  Hence,  we can take $j\geq 1$ so that $A \neq 0$. 

Since $\Hom^1(A,B)=\Hom^2(A,s_j)=0$, by applying $\Hom(A,\_)$ to
\eqref{triangle for s_j}   we obtain  $\Hom^1(A,A)=0$.  Because we have $\hom^*(A,B)$, it follows that $\{ \hom^1(\Gamma,s_j)\neq
0\}_{\Gamma \in Ind(A)}$.\footnote{see the last paragraph of the proof of Lemma \ref{for pro C1.4 C2.4} with $E$ replaced by $s_j$, $A_2$  by $A$, $B_0$ by $B$,
and letting $A_1=0$}      Since there are not Ext-nontrivial
couples in $\mc A$, we obtain $\{ \hom^1(s_j,\Gamma)=
0\}_{\Gamma \in Ind(A)}$, hence $\hom^1(s_j,A)= 0$. Now the triangle
\eqref{triangle for s_j}  and $\Hom(s_j,\_)$ imply
$\hom^1(s_j,B)=0$. Finally, the same triangle and  $\Hom(\_,B[1])$
imply  $\Hom(B,B[1])=0$. Lemma \ref{main theorem for K(l)} is
proved.

\end{document}